\documentclass[twoside,final]{interact}
\usepackage{url}
\usepackage{color}
\usepackage{tabularx}
\usepackage{enumerate}
\newcommand{\colblue}[1]{{#1}}
\newcommand{\colred}[1]{{#1}}
\newcommand{\MyParagraphAbstand}[1]{\medskip\noindent\textbf{#1.}\;}
\newcommand{\RR}{\mathbb{R}}
\newcommand{\NN}{\mathbb{N}}
\newcommand{\ZZ}{\mathbb{Z}}
\newcommand{\Sone}{\mathbb{S}_1}
\newcommand{\ee}{\mathrm{e}}
\newcommand{\dd}{\mathrm{d}}
\newcommand{\ii}{\mathrm{i}}
\newcommand{\nA}{\mathcal{A}}
\newcommand{\nM}{\mathcal{M}}
\newcommand{\MySin}{\sin\!}
\newcommand{\MyCos}{\cos\!}
\newcommand{\absbig}[1]{\big|#1\big|}
\newcommand{\normbig}[1]{\big\|#1\big\|}
\newcommand{\CompProduct}{\,.\!*}
\begin{document}
\articletype{Contribution to Special Issue \\
Geometric Numerical Integration, 25 Years Later}
\title{On the reliable and efficient numerical integration of the \\ Kuramoto model and related dynamical systems on graphs\thanks{Dedicated to Jes{\'u}s Mar{\'i}a Sanz-Serna due to his seminal contributions in the area of Geometric numerical integration.}}
\author{\name{T.~B{\"o}hle\textsuperscript{a}, C.~Kuehn\textsuperscript{a}, M.~Thalhammer\textsuperscript{b}\thanks{M.~Thalhammer. Email: mechthild.thalhammer@uibk.ac.at}}
\affil{\textsuperscript{a}Technische Universit{\"a}t M{\"u}nchen, Fakult{\"a}t f{\"u}r Mathematik, \\ Boltzmannstra{\ss}e~3, 85748~Garching, Germany. \\ \textsuperscript{b}Leopold--Franzens Universit\"{a}t Innsbruck, Institut f{\"u}r Mathematik, \\ Technikerstra{\ss}e 13/7, 6020 Innsbruck, Austria.}}
\maketitle
\begin{abstract}
In this work, a novel approach for the reliable and efficient numerical integration of the Kuramoto model on graphs is studied.
For this purpose, the notion of order parameters is revisited for the classical Kuramoto model describing all-to-all interactions of a set of oscillators. 
First numerical experiments confirm that the precomputation of certain sums significantly reduces the computational cost for the evaluation of the right-hand side and hence enables the simulation of high-dimensional systems.
In order to design numerical integration methods that are favourable in the context of related dynamical systems on network graphs, the concept of localised order parameters is proposed. 
In addition, the detection of communities for a complex graph and the transformation of the underlying adjacency matrix to block structure is an essential component for further improvement.
It is demonstrated that for a submatrix comprising relatively few coefficients equal to zero, the precomputation of sums is advantageous, whereas straightforward summation is appropriate in the complementary case.
Concluding theoretical considerations and numerical comparisons show that the strategy of combining effective community detection algorithms with the localisation of order parameters potentially reduces the computation time by several orders of magnitude.
\end{abstract}
\section{Introduction}
In the present work, we propose a novel approach for the reliable and efficient numerical integration of nonlinear dynamical systems on network graphs and provide various numerical comparisons confirming its potential.
In essence, our objective is to combine effective algorithms for the detection of communities with the concept of localised order parameters based on the precomputation of certain sums.
For the sake of concreteness, we focus on Kuramoto-type models for a large set of individual oscillators \colred{and graphs that are determined by adjacency matrices with coefficients equal to one and zero, respectively.}
Generally, the interactions between the oscillators are described by a system of nonlinear ordinary differential equations involving the sines of the associated phases.
The classical Kuramoto--Daido model~\cite{Daido1996,Kuramoto1984} reflects the case of all-to-all coupling and thus corresponds to a complete graph.
Realistic extended models incorporate complex networks~\cite{RodriguesEtAl2}.
Yet, the efficient numerical integration of a Kuramoto-type model on a graph comprising a high number of nodes is a relevant issue, and even the evaluation of the vector field defining the right-hand side of the system poses a major challenge.
This problem actually permeates the simulation of complex dynamics on networks~\cite{PorterGleeson}.
In order to exemplify our strategy and to illustrate its capability, we initially consider the classical Kuramoto--Daido model and revisit the well-known notion of order parameters~\cite{ArenasEtAl,Daido1996,Kuramoto1984,PikovskyRosenblumKurths,Strogatz2000,Strogatz1991}.
An essential component of our procedure for Kuramoto-type models on graphs is the detection of communities~\cite{GirvanNewman,PorterOnnelaMucha}, since this permits a transformation of the associated adjacency matrix to block structure.
We demonstrate that for a submatrix comprising relatively few coefficients equal to zero the precomputation of sums is indeed advantageous and that straightforward summation is suitable in the complementary case.
We conclude with theoretical considerations and numerical comparisons, which show that our approach ensures a significant reduction of the required memory capacity as well as the computational cost for the evaluation of the right-hand side.
As a consequence, long-term simulations of systems involving a large number of oscillators by geometric integrators~\cite{BlanesCasas,IserlesQuispel,SanzSerna,SanzSernaCalvo} are within reach. 

\MyParagraphAbstand{Scope of the model}
The starting point of our investigations is the Kuramoto--Daido model~\cite{Daido1996,Kuramoto1984}.
We henceforth refer to it as (classical) Kuramoto model.
This system of coupled nonlinear ordinary differential equations is a fundamental mathematical model for the dynamical behaviour of a set of weakly coupled, nearly identical oscillators and specifies the time evolution of the associated phases.
Despite its simple structure, the Kuramoto model exhibits fascinating phenomena such as synchronisation and phase locking~\cite{Strogatz2000}.
Originally introduced to describe processes in chemistry and biology~\cite{Kuramoto1984,Winfree1967}, it was found to have various applications in other fields such as physics, neuroscience, and engineering~\cite{AcebronEtAl,DoerflerBullo,Golomb2001}.
From a theoretical perspective, the Kuramoto model has deep connections to effects present in Hamiltonian systems, particularly to Landau damping~\cite{Dietert,FernandezGerardVaretGiacomin} and bifurcations from essential spectra~\cite{Chiba}.
A variety of Kuramoto-type models have a gradient flow structure entailing further interesting mathematical cross-connections~\cite{HaKimRyoo}.

\MyParagraphAbstand{Kuramoto model}
The main terms defining the right-hand side of the classical Kuramoto model for~$M$ individual oscillators have the form 
\begin{equation}
\label{eq:Sum1}
\tfrac{1}{M} \, \sum_{\ell=1}^{M} \MySin\big(\vartheta_{\ell}(t) - \vartheta_m(t)\big)\,, \quad t \in [0, T]\,, \quad m \in \{1, 2, \dots, M\}\,.
\end{equation}
Here, $\vartheta_m$ denotes the time-dependent phase of the $m$-th oscillator, which takes values in the circle $\Sone = \RR/(2\pi\ZZ)$.
Henceforth, we employ the convenient vector notation
\begin{equation*}
\vartheta = (\vartheta_1, \dots, \vartheta_M)^T: [0, T] \longrightarrow \Sone^M\,.
\end{equation*}
With regard to the numerical simulation of a high number of oscillators, an essential requirement is the efficient evaluation of these sums at certain time grid points.
Thereby, one issue is the limited memory capacity.
As an example, we mention the currently widely used software \textsc{Matlab}, which has a maximum array size preference of about~$10^{10}$ (74.5\,\textsc{GB}) corresponding to a square matrix of dimension~$10^5$.
Consequently, it is desirable to avoid the creation of the matrix
\begin{equation*}
\big(\vartheta_{\ell}(t) - \vartheta_m(t)\big)_{\ell,m \in \{1, 2, \dots, M\}}\,, \quad t \in [0, T]\,,
\end{equation*}
since this would restrict the dimension of the system considerably.
An important alternative for the efficient evaluation of~\eqref{eq:Sum1} relies on the macroscopic order parameters, which are given through \begin{equation*}
r\big(\vartheta(t)\big) \, \ee^{\ii \, \psi(\vartheta(t))} = \tfrac{1}{M} \sum_{m=1}^{M} \ee^{\ii \, \vartheta_m(t)}\,, \quad t \in [0, T]\,.
\end{equation*}
More precisely, applying the addition theorem for the sine function to~\eqref{eq:Sum1}, we obtain the following reformulation
\begin{equation}
\label{eq:Sum2}
\begin{gathered}
S_M\big(\vartheta(t)\big) = \tfrac{1}{M} \sum_{m=1}^{M} \MySin\big(\vartheta_m(t)\big)\,, \quad 
C_M\big(\vartheta(t)\big) = \tfrac{1}{M} \sum_{m=1}^{M} \MyCos\big(\vartheta_m(t)\big)\,, \\
\tfrac{1}{M} \sum_{\ell=1}^{M} \MySin\big(\vartheta_{\ell}(t) - \vartheta_m(t)\big)
= S_M\big(\vartheta(t)\big) \, \MyCos\big(\vartheta_m(t)\big)
- C_M\big(\vartheta(t)\big) \, \MySin\big(\vartheta_m(t)\big)\,, \\
t \in [0, T]\,, \quad m \in \{1, 2, \dots, M\}\,.
\end{gathered}
\end{equation}
An evident though crucial observation is that precomputing the sums~$S_M(\vartheta(t))$ and~$C_M(\vartheta(t))$ permits to evaluate
$S_M(\vartheta(t)) \, \cos(\vartheta_m(t)) - C_M(\vartheta(t)) \, \sin(\vartheta_m(t))$ for $m \in \{1, 2, \dots, M\}$ in an efficient manner.
Numerical comparisons described in Section~\ref{sec:SectionKuramoto} confirm that this approach reduces the number of function evaluations and thus the computation time considerably.

\MyParagraphAbstand{Kuramoto model on graphs}
As indicated above, it is of high relevance to study extensions of the classical Kuramoto model in the context of dynamical networks~\cite{RodriguesEtAl2}.
We focus on situations, where the sum over all phases is replaced by a sum over certain phases.
That is, the right-hand side of the system involves terms of the form 
\begin{equation*}
\sum_{\ell=1}^{M} A_{m \ell} \, \MySin\big(\vartheta_{\ell}(t) - \vartheta_m(t)\big)\,, \quad A_{m \ell} \in \{0, 1\}\,, \quad t \in [0, T]\,, \quad
\ell, m \in \{1, 2, \dots, M\}\,.
\end{equation*}
The associated matrix
\begin{equation*}
A = \big(A_{m \ell}\big)_{\ell,m \in \{1, \dots, M\}} \in \RR^{M \times M}
\end{equation*}
has the natural interpretation as the adjacency matrix of a graph.
In Section~\ref{sec:SectionKuramotoGraph}, we consider the complementary cases of sparse and dense adjacency matrices.
Numerical tests for randomly generated matrices show that the precomputation of sums is advantageous whenever the number of coefficients equal to one is larger than the number of coefficients equal to zero, whereas straightforward summation is appropriate otherwise.
Subsequently, we examine algorithms detecting communities in a graph, which yield as outputs partitions of the nodes, and extend our approach to block matrices.
In our considerations, we do not presume a symmetric adjacency matrix, which would lead to a gradient system. 

\MyParagraphAbstand{Generalisations}
\colred{
We point out that our approach applies to the more general case, where the coefficients of the adjacency matrix take values in a finite set, for instance
\begin{equation*}
A_{m \ell} \in \{0, 1, \dots, J\}\,, \quad \ell, m \in \{1, 2, \dots, M\}\,.
\end{equation*}
Provided that a suitable permutation permits the transformation to a block matrix such that a single value is prevalent in each block and the total number of blocks is relatively small, we may expect a significant gain in efficiency by the precomputation of sums.}
Moreover, it is straightforward to generalise our approach to systems that in addition involve multiple sums such as 
\begin{equation*}
\sum_{j,k,\ell=1}^{M} A_{m j k \ell}\MySin\big(\vartheta_j(t) + \vartheta_k(t) - \vartheta_{\ell}(t) - \vartheta_m(t)\big)\,, \quad t \in [0, T]\,, \quad m \in \{1, 2, \dots, M\}\,.
\end{equation*}
Due to the fact that higher-order Kuramoto models describe interactions beyond usual graph structures, they have recently raised remarkable interest, see~\cite{BickBoehleKuehn,SkardalArenas} and references given therein.
However, as their incorporation requires laborious notation and would obstruct a comprehensible presentation of the key idea, we do not include detailed calculations here.

\MyParagraphAbstand{Remarks}
We note that suitable reformulations of the classical Kuramoto model and Kuramoto-type models on graphs can also be based on complex exponentials.
In view of expedient practical implementations, however, it is preferable to use real-valued quantities.
On account of cross-references in captions, figures are included in the end of this manuscript. 
\section{Kuramoto model}
\label{sec:SectionKuramoto}
In this section, we state the classical Kuramoto model and study different viewpoints on its reliable and efficient numerical integration.

\MyParagraphAbstand{Original formulation}
We consider a set of $M$~limit-cycle oscillators with time-dependent phases 
\begin{subequations}
\label{eq:Kuramoto}  
\begin{equation}
\vartheta_m: [0, T] \longrightarrow \Sone\,, \quad m \in \{1, 2, \dots, M\}\,.
\end{equation}
In the absence of external driving or damping forces, respectively, the oscillators have the intrinsic frequencies
\begin{equation}
\omega_m \in \RR\,, \quad m \in \{1, 2, \dots, M\}\,.
\end{equation}
The pairwise interactions between the oscillators are described by the following system of nonlinear ordinary differential equations
\begin{equation}
\begin{cases}
\displaystyle \vartheta_m'(t) = \omega_m + \tfrac{K}{M} \sum_{\ell=1}^{M} \MySin\big(\vartheta_{\ell}(t) - \vartheta_m(t)\big)\,, \\
\vartheta_m(0) \text{ given}\,, \quad t \in (0, T)\,, \quad m \in \{1, 2, \dots, M\}\,,
\end{cases}
\end{equation}
where $K > 0$ denotes the coupling constant.

\MyParagraphAbstand{Special choices}
In our numerical tests, we consider intrinsic frequencies defined by a single real number $\omega_0 \in \RR$ and initial phases of the form
\begin{equation}
\omega_m = 1 + \omega_0 \, \tfrac{(2 m - M - 1)}{M - 1}\,, \quad 
\vartheta_m(0) = \tfrac{2 \pi m}{M}\,, \quad m \in \{1, 2, \dots, M\}\,.
\end{equation}
\end{subequations}

\MyParagraphAbstand{Reformulation}
Regarding the efficient numerical integration of the classical Kuramoto model, it is to the best advantage to employ a reformulation that relies on elementary addition theorems for trigonometric functions  
\begin{equation*}
\begin{gathered}
\MySin\big(\vartheta_{\ell}(t) - \vartheta_m(t)\big)
= \MySin\big(\vartheta_{\ell}(t)\big) \MyCos\big(\vartheta_m(t)\big) - \MyCos\big(\vartheta_{\ell}(t)\big) \MySin\big(\vartheta_m(t)\big)\,, \\
\quad \ell, m \in \{1, 2, \dots, M\}\,.
\end{gathered}
\end{equation*}
Introducing the abbreviations
\begin{subequations}
\label{eq:KuramotoReformulation}
\begin{equation}
\begin{gathered}
S_M\big(\vartheta(t)\big) = S_M\big(\vartheta_1(t), \vartheta_2(t), \dots, \vartheta_M(t)\big) = \tfrac{1}{M} \sum_{m=1}^{M} \MySin\big(\vartheta_m(t)\big) \in \RR\,, \\
C_M\big(\vartheta(t)\big) = C_M\big(\vartheta_1(t), \vartheta_2(t), \dots, \vartheta_M(t)\big) = \tfrac{1}{M} \sum_{m=1}^{M} \MyCos\big(\vartheta_m(t)\big) \in \RR\,, 
\end{gathered}
\end{equation}
the governing equations~\eqref{eq:Kuramoto} read as 
\begin{equation}
\begin{cases}
\vartheta_m'(t) = \omega_m + K \, \Big(S_M\big(\vartheta(t)\big) \MyCos\big(\vartheta_m(t)\big) 
- C_M\big(\vartheta(t)\big) \MySin\big(\vartheta_m(t)\big)\Big)\,, \\
\vartheta_m(0) \text{ given}\,, \quad t \in (0, T)\,, \quad m \in \{1, 2, \dots, M\}\,.
\end{cases}
\end{equation}
Employing the standard vector notation
\begin{equation}
\begin{gathered}
\omega = \begin{pmatrix} \omega_1 \\ \omega_2 \\ \vdots \\ \omega_M \end{pmatrix} \in \RR^M\,, \quad
\vartheta(t) = \begin{pmatrix} \vartheta_1(t) \\ \vartheta_2(t) \\ \vdots \\ \vartheta_M(t) \end{pmatrix} \in \Sone^M\,, \quad t \in [0, T]\,,
\end{gathered}
\end{equation}
and setting accordingly 
\begin{equation}
\begin{gathered}
\MySin\big(\vartheta(t)\big) = \begin{pmatrix} \MySin\big(\vartheta_1(t)\big) \\ \MySin\big(\vartheta_2(t)\big) \\ \vdots \\ \MySin\big(\vartheta_M(t)\big) \end{pmatrix}\,, \quad
\MyCos\big(\vartheta(t)\big) = \begin{pmatrix} \MyCos\big(\vartheta_1(t)\big) \\ \MyCos\big(\vartheta_2(t)\big) \\ \vdots \\ \MyCos\big(\vartheta_M(t)\big) \end{pmatrix}\,, \\
F\big(\vartheta(t)\big) = \omega + K \, \Big(S_M\big(\vartheta(t)\big) \MyCos\big(\vartheta(t)\big) - C_M\big(\vartheta(t)\big) \MySin\big(\vartheta(t)\big)\Big)\,, \quad t \in [0, T]\,, 
\end{gathered}
\end{equation}
the system takes the compact form  
\begin{equation}
\begin{cases}
\vartheta'(t) = F\big(\vartheta(t)\big)\,, \quad t \in (0, T)\,, \\
\vartheta(0) \text{ given}\,.
\end{cases}
\end{equation}
\end{subequations}

\MyParagraphAbstand{Potential}
The classical Kuramoto model~\eqref{eq:Kuramoto}--\eqref{eq:KuramotoReformulation} has the intrinsic structure of a gradient system. 
That is, the right-hand side is given by the gradient of a real-valued potential function
\begin{equation*}
\begin{split}
- \, \nabla V(\vartheta)
&= \bigg(\omega_m + \tfrac{K}{M} \sum_{\ell=1}^{M} \sin(\vartheta_{\ell} - \vartheta_m)\bigg)_{m \in \{1, 2, \dots, M\}} \\
&= \omega + K \, \big(S_M(\vartheta) \cos(\vartheta) - C_M(\vartheta) \sin(\vartheta)\big)\,, \quad
\vartheta = (\vartheta_1, \dots, \vartheta_M)^T \in \Sone^M\,,
\end{split}
\end{equation*}
such that the governing equations rewrite as 
\begin{equation*}
\begin{cases}
\vartheta'(t) = - \nabla V\big(\vartheta(t)\big)\,, \quad t \in (0, T)\,, \\
\vartheta(0) \text{ given}\,.
\end{cases}
\end{equation*}
This in particular implies that the values of the potential decrease when time evolves 
\begin{equation*}
\begin{gathered}
\tfrac{\dd}{\dd t} \, V\big(\vartheta(t)\big)
= \Big(\nabla V\big(\vartheta(t)\big)\Big)^T \vartheta'(t)
= - \, \normbig{\nabla V\big(\vartheta(t)\big)}^2 \leq 0\,, \\
V\big(\vartheta(t)\big) \leq V\big(\vartheta(0)\big)\,, \quad t \in [0, T]\,.
\end{gathered}
\end{equation*}
For our purposes, in view of its efficient evaluation, it is advantageous to reformulate the canonical potential as follows 
\begin{subequations}
\label{eq:PotentialOrderParameterConservation}
\begin{equation}
\label{eq:Potential}
\begin{split}
V: \RR^M &\longrightarrow \RR: \\
\vartheta
&\longmapsto 
V(\vartheta) = V(\vartheta_1, \dots, \vartheta_M) \\
&\qquad\qquad\,\,= - \sum_{m=1}^M \omega_m \, \vartheta_m + \tfrac{K}{2 M} \sum_{\ell,m=1}^{M} \big(1 - \cos(\vartheta_{\ell} - \vartheta_m)\big) \\
&\qquad\qquad\,\,= - \, \omega^T \vartheta + \tfrac{K M}{2} \, \Big(1 - \big(C_M(\vartheta)\big)^2 - \big(S_M(\vartheta)\big)^2\Big)\,.
\end{split}
\end{equation}
\colred{As mentioned below, the choice of the constant of integration is linked to the special case of synchronisation.}

\MyParagraphAbstand{Order parameter}
The modulus $r: \Sone^M \to \RR$ and the angle $\psi: \Sone^M \to \RR$ of the complex order parameter are given by
\begin{equation}
\label{eq:OrderParameter}
r(\vartheta) \, \ee^{\ii \, \psi(\vartheta)} = \tfrac{1}{M} \sum_{m=1}^{M} \ee^{\ii \, \vartheta_m} = C_M(\vartheta) + \ii \, S_M(\vartheta)\,, \quad \vartheta \in \Sone^M\,.
\end{equation}
Multiplying by~$\ee^{- \, \ii \, \vartheta_{\ell}}$ for $\ell \in \{1, 2, \dots, M\}$ and considering the imaginary part of the resulting relation 
\begin{equation*}
\begin{gathered}
r(\vartheta) \, \ee^{\ii \, (\psi(\vartheta) - \vartheta_{\ell})} = \tfrac{1}{M} \sum_{m=1}^{M} \ee^{\ii \, (\vartheta_m - \vartheta_{\ell})}\,, \quad
r(\vartheta) \, \MySin\big(\psi(\vartheta) - \vartheta_{\ell}\big) = \tfrac{1}{M} \, \sum_{m=1}^{M} \MySin\big(\vartheta_m - \vartheta_{\ell}\big)\,, \\
\vartheta \in \Sone^M\,, \quad \ell \in \{1, 2, \dots, M\}\,,
\end{gathered}
\end{equation*}
 the Kuramoto model~\eqref{eq:Kuramoto}--\eqref{eq:KuramotoReformulation} rewrites as 
\begin{equation*}
\begin{cases}
\vartheta'(t)
= \omega + K \, r\big(\vartheta(t)\big) \, \MySin\Big(\psi\big(\vartheta(t)\big) - \vartheta(t)\Big)\,, \\
\colblue{\vartheta(0)} \text{ given}\,, \quad t \in (0, T)\,.
\end{cases}
\end{equation*}
We point out that the order parameters contain the entire information about the interactions of the oscillators. 
Although this reformulation looks like a mean-field equation for a single oscillator, it completely represents the original system.  

\MyParagraphAbstand{Indicator for synchronisation}
\colblue{For configurations, where all cosine and sine values are close-by, the modulus of the complex order parameter has values nearly one
\begin{equation*}
\begin{gathered}
\cos(\vartheta_m) \approx \cos(\vartheta_1)\,, \quad \sin(\vartheta_m) \approx \sin(\vartheta_1)\,, \quad m \in \{2, 3, \dots, M\}\,, \\
C_M(\vartheta) \approx \cos(\vartheta_1)\,, \quad S_M(\vartheta) \approx \sin(\vartheta_1)\,, \\
r(\vartheta) = \sqrt{\big(C_M(\vartheta)\big)^2 + \big(S_M(\vartheta)\big)^2} \approx 1\,, \quad \vartheta \in \Sone^M\,.
\end{gathered}
\end{equation*}
Hence, this quantity indicates synchronisation. 
Furthermore, in this situation, the above stated choice of the constant of integration in the potential implies 
\begin{equation*}
V(\vartheta) \approx - \, \omega^T \vartheta\,, \quad \vartheta \in \Sone^M\,.
\end{equation*}
}

\MyParagraphAbstand{Conserved quantity}
\colred{A straightforward calculation shows that summation over all governing equations yields the identity
\begin{equation*}
\tfrac{1}{M} \sum_{m=1}^{M} \vartheta_m'(t) = \tfrac{1}{M} \sum_{m=1}^{M} \omega_m\,, \quad t \in [0, T]\,, 
\end{equation*}
see~\eqref{eq:Kuramoto}--\eqref{eq:KuramotoReformulation}.
Performing integration, this implies that the mean values of the intrinsic frequencies and the initial phases determine the mean values of the phases at later times
\begin{equation*}
\tfrac{1}{M} \sum_{m=1}^{M} \vartheta_m(t) = \tfrac{1}{M} \sum_{m=1}^{M} \vartheta_m(0) + t \, \tfrac{1}{M} \sum_{m=1}^{M} \omega_m\,, \quad t \in [0, T]\,.
\end{equation*}
In other words, the dynamics of the classical Kuramoto model is restricted to a time-dependent submanifold defined by the constraint 
\begin{equation}
\label{eq:5c}
\tfrac{1}{M} \sum_{m=1}^{M} \big(\vartheta_m(t) - \vartheta_m(0) - t \, \omega_m\big) = 0\,, \quad t \in [0, T]\,.
\end{equation}
We anticipate that this characteristic result extends to Kuramoto models on graphs under a certain symmetry condition, see Section~\ref{sec:SectionKuramotoGraph}.
However, in general, this conservation property does not hold and then the quantity   
\begin{equation}
\label{eq:Means}
\tfrac{1}{M} \sum_{m=1}^{M} \big(\vartheta_m(t) - \vartheta_m(0) - t \, \omega_m\big)
\end{equation}
reflects the deviation.

In the limit $M \to \infty$, the sums in \eqref{eq:5c}--\eqref{eq:Means} are replaced by integrals.
In case of randomly chosen problem data, e.g., they are interpreted as expected phases and intrinsic frequency, respectively.
}
\end{subequations}

\begin{table}[t]
\tbl{Classical Kuramoto model. Computational cost for the evaluation of the right-hand side based on the original formulation~\eqref{eq:Kuramoto} and the reformulation~\eqref{eq:KuramotoReformulation}, respectively. Number of function evaluations (sine, cosine) in dependence of the total number of oscillators. \\}
{\qquad\qquad\qquad
\begin{tabular}{c|c} 
Original formulation & $M \, (M-1)$ \\\midrule
Reformulation & $4 \, M$ 
\end{tabular}\qquad\qquad\qquad}
\label{tab:Table1}
\end{table}

\MyParagraphAbstand{Implementation and computational cost}
For systems involving a high number of oscillators $M >\!\!> 1$, the above stated approach leading to a reformulation
of the classical Kuramoto model is beneficial in several respects.
It permits to reduce the required memory capacity as well as the computational cost for the evaluation of the right-hand side significantly, see Table~\ref{tab:Table1}.
Moreover, parallelisation techniques can be used.

\MyParagraphAbstand{Exemplification (Euler)}
In order to exemplify our procedure for the efficient numerical integration of the Kuramoto model, we consider the simplest first-order one-step method, the explicit Euler method.
For a time grid with associated stepsizes 
\begin{equation*}
0 = t_0 < \dots < t_n < \dots < t_N = T\,, \quad \tau_n = t_{n+1} - t_n\,, \quad n \in \{0, 1, \dots, N-1\}\,, 
\end{equation*}
and a prescribed initial approximation $\vartheta^{(0)} \approx \vartheta(0)$, the explicit Euler solution is given by the recurrence 
\begin{equation*}
\vartheta^{(n+1)} = \vartheta^{(n)} + \tau_n \, F\big(\vartheta^{(n)}\big)\,, \quad n \in \{0, 1, \dots, N-1\}\,, 
\end{equation*}
see~\eqref{eq:KuramotoReformulation}.
In each time step, the evaluation of the defining function relies on the precomputation of the sums
\begin{subequations}
\label{eq:EulerCost}
\begin{equation}
S_M^{(n)} = S_M\big(\vartheta^{(n)}\big)\,, \quad C_M^{(n)} = C_M\big(\vartheta^{(n)}\big)\,, 
\end{equation}
and subsequently on the computation of 
\begin{equation}
F\big(\vartheta^{(n)}\big) = \omega + K \, \Big(S_M^{(n)} \MyCos\big(\vartheta^{(n)}\big) - C_M^{(n)} \MySin\big(\vartheta^{(n)}\big)\Big)\,.
\end{equation}
\end{subequations}
Altogether, this requires~$4 \, M$ evaluations of sine and cosine functions compared to~$M (M-1)$ evaluations of the sine function needed for the original formulation~\eqref{eq:Kuramoto}. 
The generalisation to higher-order explicit or implicit time integration methods is straightforward. 
Long-term simulations are ideally based on geometric integrators.
\colblue{Their benefits over standard methods are demonstrated in~\cite{BlanesCasas,IserlesQuispel,SanzSerna,SanzSernaCalvo}, e.g.} 

\MyParagraphAbstand{Remark}
In our study, the focus is on the numerical simulation of a high number of oscillators.
In this situation, the reduction from $M^2$ to $M^2 - M$ evaluations of the sine function by taking into account the evident identity for coinciding indices
\begin{equation*}
\ell = m: \quad \sin(\vartheta_{\ell} - \vartheta_m) = 0
\end{equation*}
is of minor relevance and thus will be neglected. 

\begin{figure}[t!]
\begin{center}
\includegraphics[width=4.4cm]{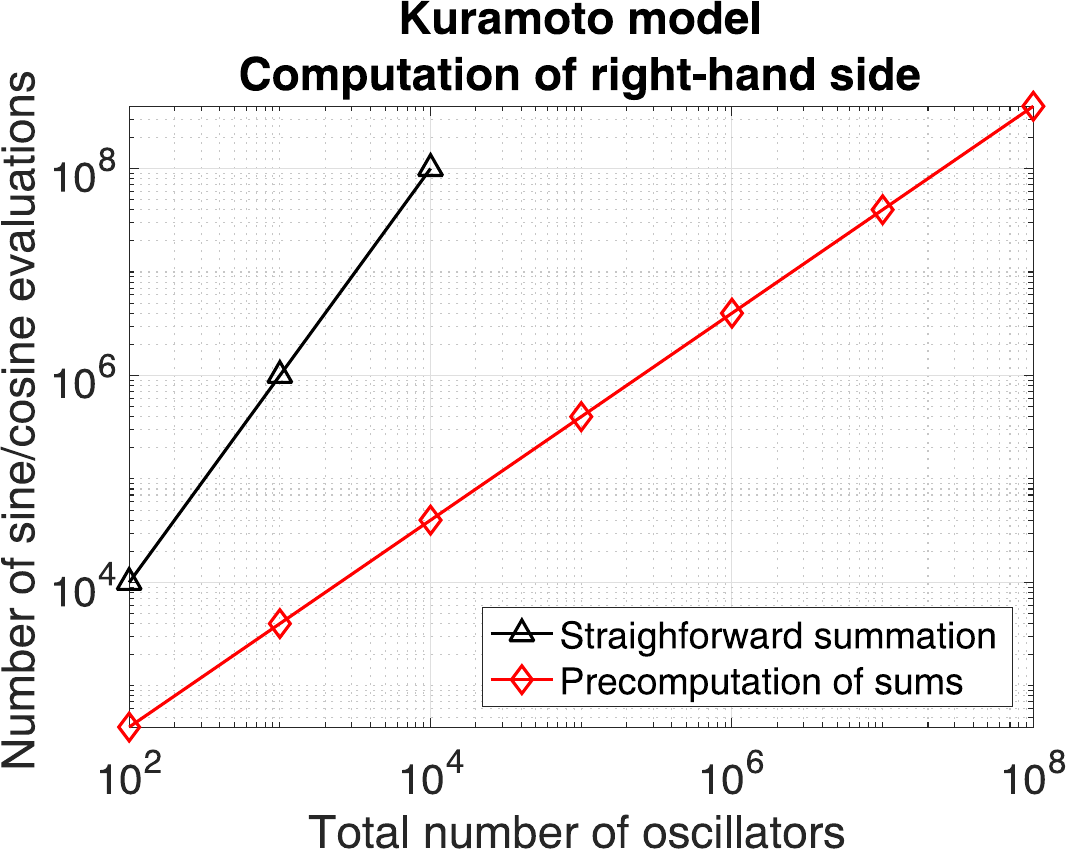} \quad
\includegraphics[width=4.4cm]{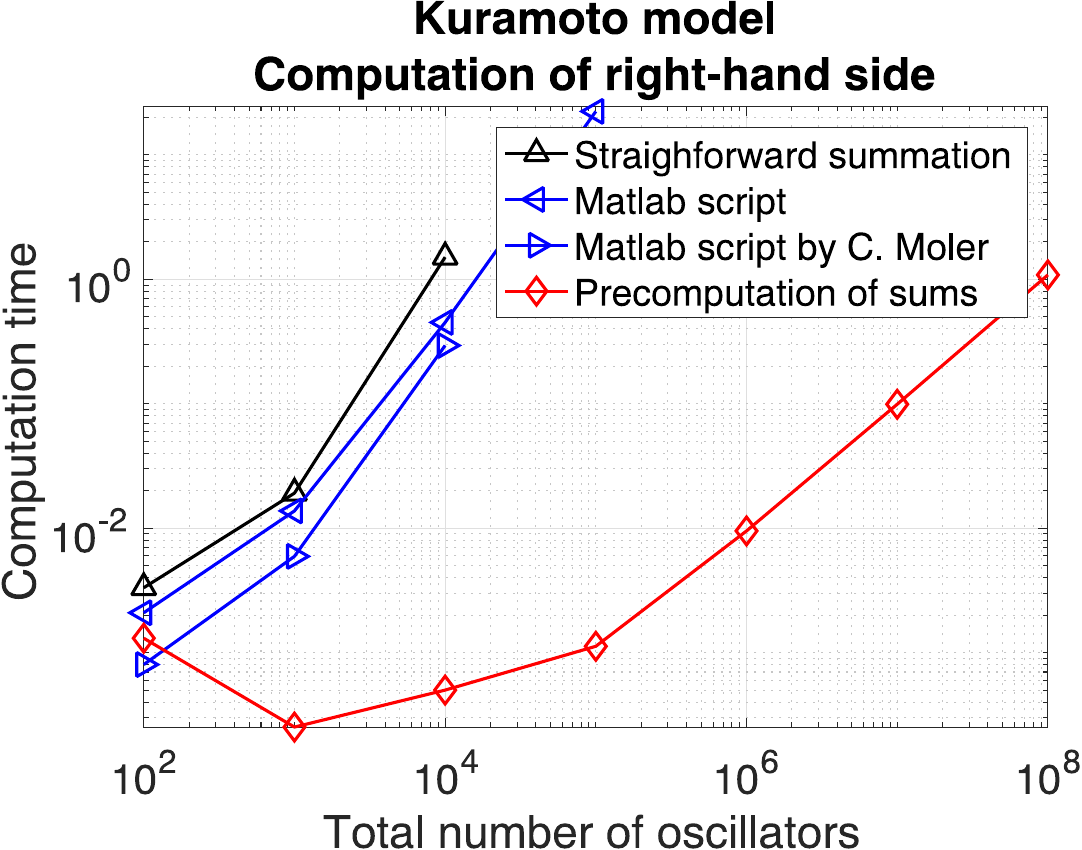} \quad
\includegraphics[width=4.4cm]{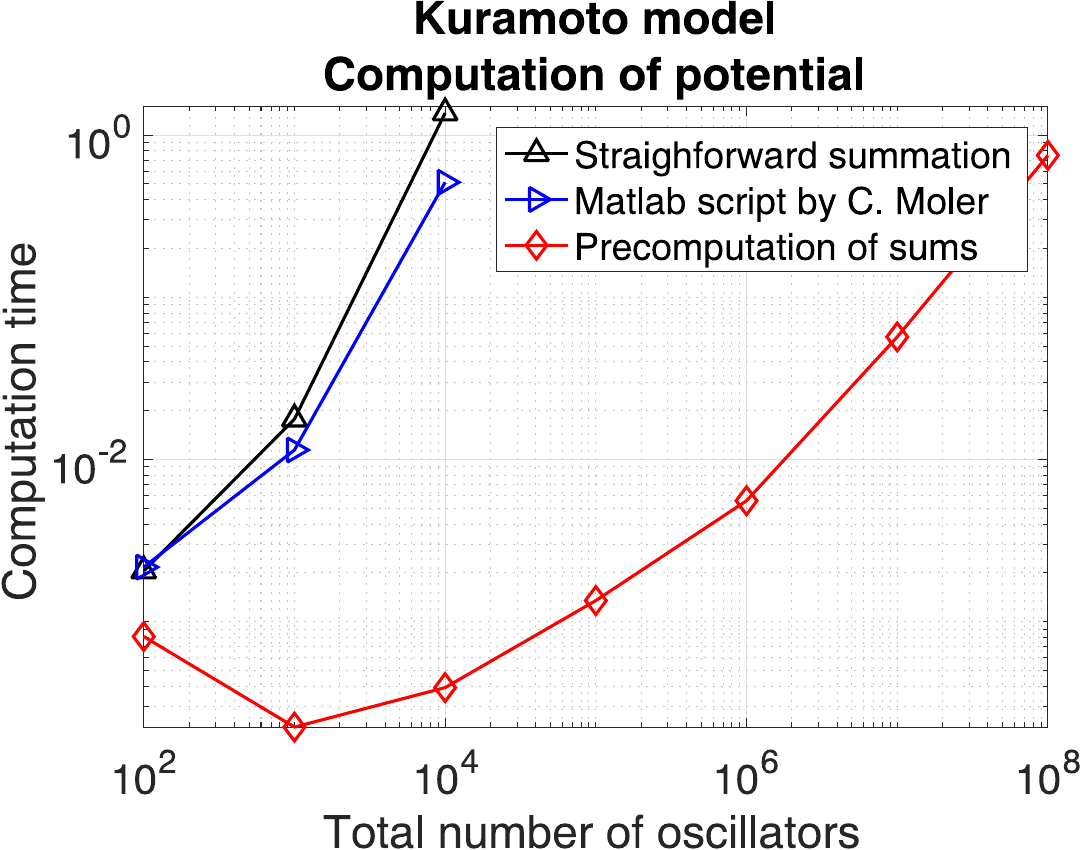}
\end{center}
\caption{Classical Kuramoto model. Computational cost for the evaluation of the right-hand side and the potential.
Left: Number of sine and cosine evaluations versus the total number of oscillators when using straighforward summation and the precomputation of sums, respectively. 
Middle: Numerical comparison of the computation time for different implementations in \textsc{Matlab}.
Right: Corresponding results for the potential.}
\label{fig:MyFigure1}
\end{figure}

\MyParagraphAbstand{Numerical comparisons}
Numerical comparisons of different approaches for the evaluation of the right-hand side of the classical Kuramoto model~\eqref{eq:Kuramoto}--\eqref{eq:KuramotoReformulation} and the associated potential~\eqref{eq:PotentialOrderParameterConservation} are displayed in Figure~\ref{fig:MyFigure1}.
We focus on an implementation in \textsc{Matlab} and expect that analogous conclusions hold for other software packages. 
We vary the total number of oscillators from~$10^2$ to~$10^8$, taking into account maximum array size preferences as mentioned in the introduction.
A randomly chosen number $\omega_0 \in (0, 1)$ defines the intrinsic frequencies.
The phases $\vartheta_1, \vartheta_2, \dots, \vartheta_M$ are uniformly distributed in $[0, 2 \pi]$. 
We compare 
\begin{enumerate}[(i)]
\item   
straightforward summation of $\sin(\vartheta_{\ell}-\vartheta_m)$ for $\ell, m \in \{1, 2, \dots, M\}$ realised by a double loop over all rows,
\item 
a simple script using $\text{sum}(\sin(\vartheta-\vartheta_m))$ for $m \in \{1, 2, \dots, M\}$ and thus involving a loop over all rows,
\item 
a script by Cleve Moler that generates the matrix $\sin(\vartheta-\vartheta')$ and then adds up each row,%
\footnote{The \textsc{Matlab} script kuramoto.m by Cleve Moler is available at \\ \url{https://de.mathworks.com/matlabcentral/fileexchange/72534-kuramoto-s-model-of-synchronizing-oscillators}.}
\item 
the above stated procedure using the precomputation of sums, see~\eqref{eq:KuramotoReformulation} and~\eqref{eq:EulerCost}.
\end{enumerate}
For a higher number of oscillators, it is probable that the evaluation of functions and the computation of sums are the most time consuming components. 
The numerical results confirm that our approach~(iv) is favourable in this regime and permits the efficient simulation of a high number of oscillators. 
Similar conclusions hold for the evaluation of the potential.

\begin{figure}[t!]
\begin{center}
\includegraphics[width=9cm]{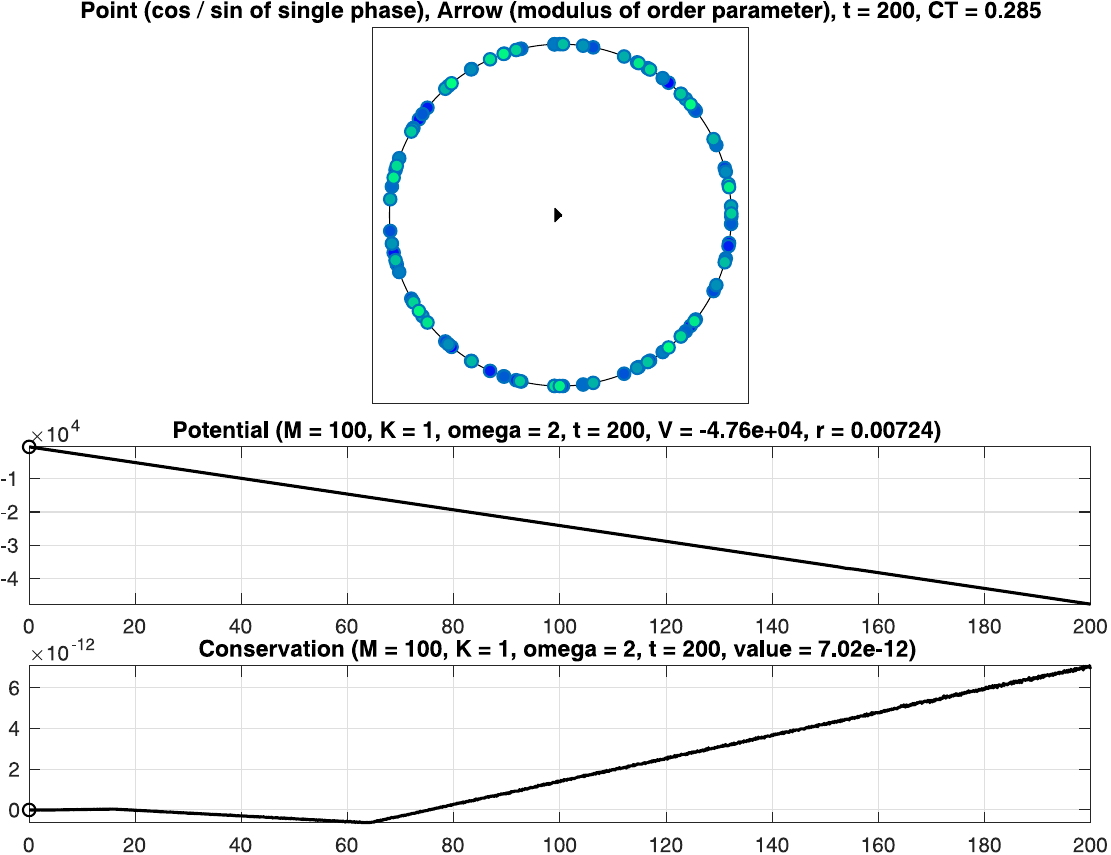}
\end{center}
\caption{\colblue{Numerical integration of the classical Kuramoto model involving $M = 10^2$ oscillators. Coupling constant $K = 1$ (no synchronisation).
Top: Visualisation of the phases at the final time.  
Middle: The time series confirms decreasing potential values. 
Bottom: A conserved quantity is numerically preserved with high accuracy.}}
\label{fig:MyFigure2}
\end{figure}

\begin{figure}[t!]
\begin{center}
\includegraphics[width=9cm]{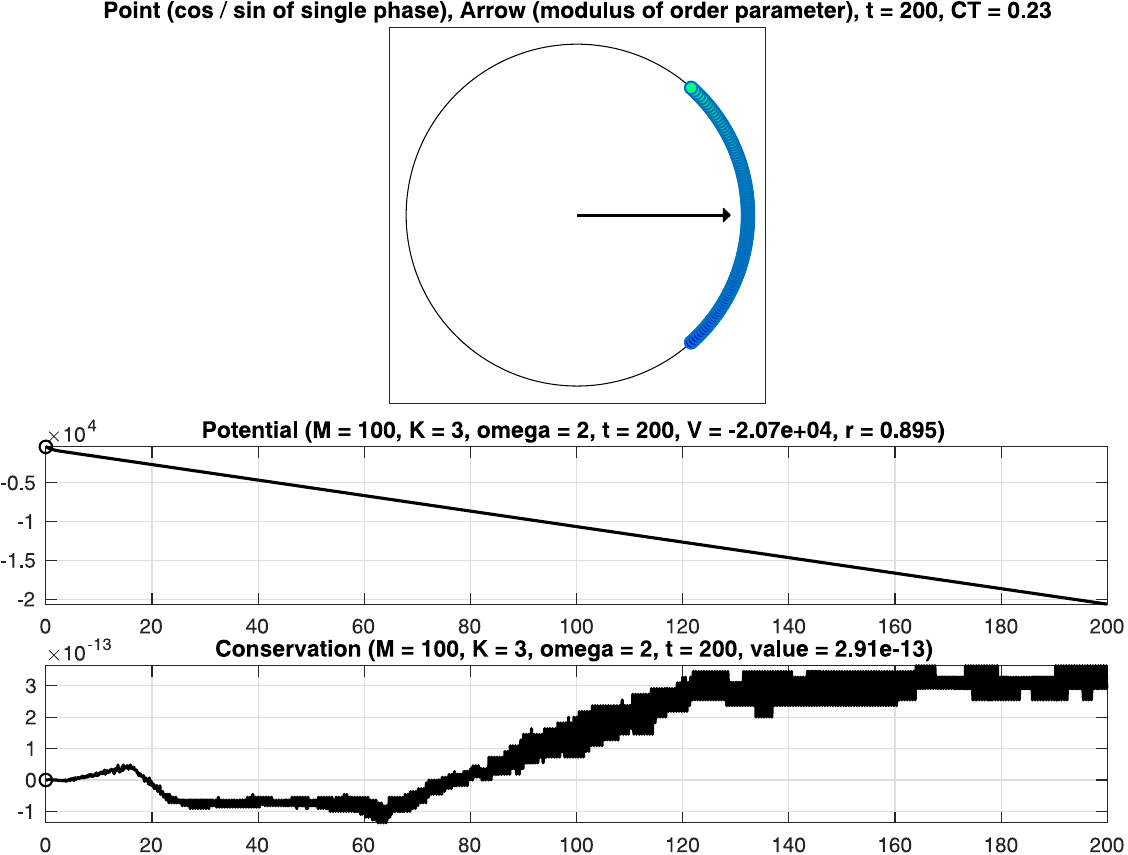} \\[4mm]
\includegraphics[width=9cm]{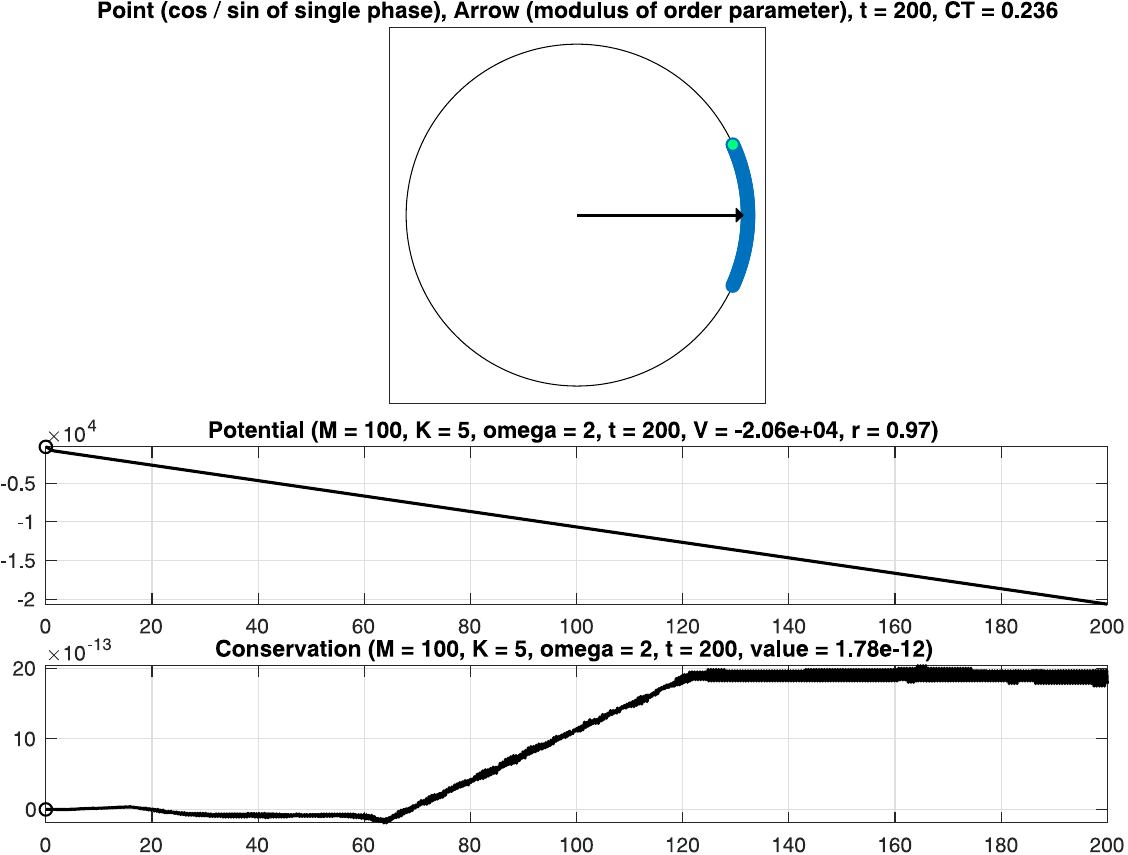} 
\end{center}
\caption{\colblue{Numerical integration of the classical Kuramoto model involving $M = 10^2$ oscillators. Coupling constant $K \in \{3, 5\}$ (gradual synchronisation).}}
\label{fig:MyFigure3}
\end{figure}

\begin{figure}[t!]
\begin{center}
\includegraphics[width=9cm]{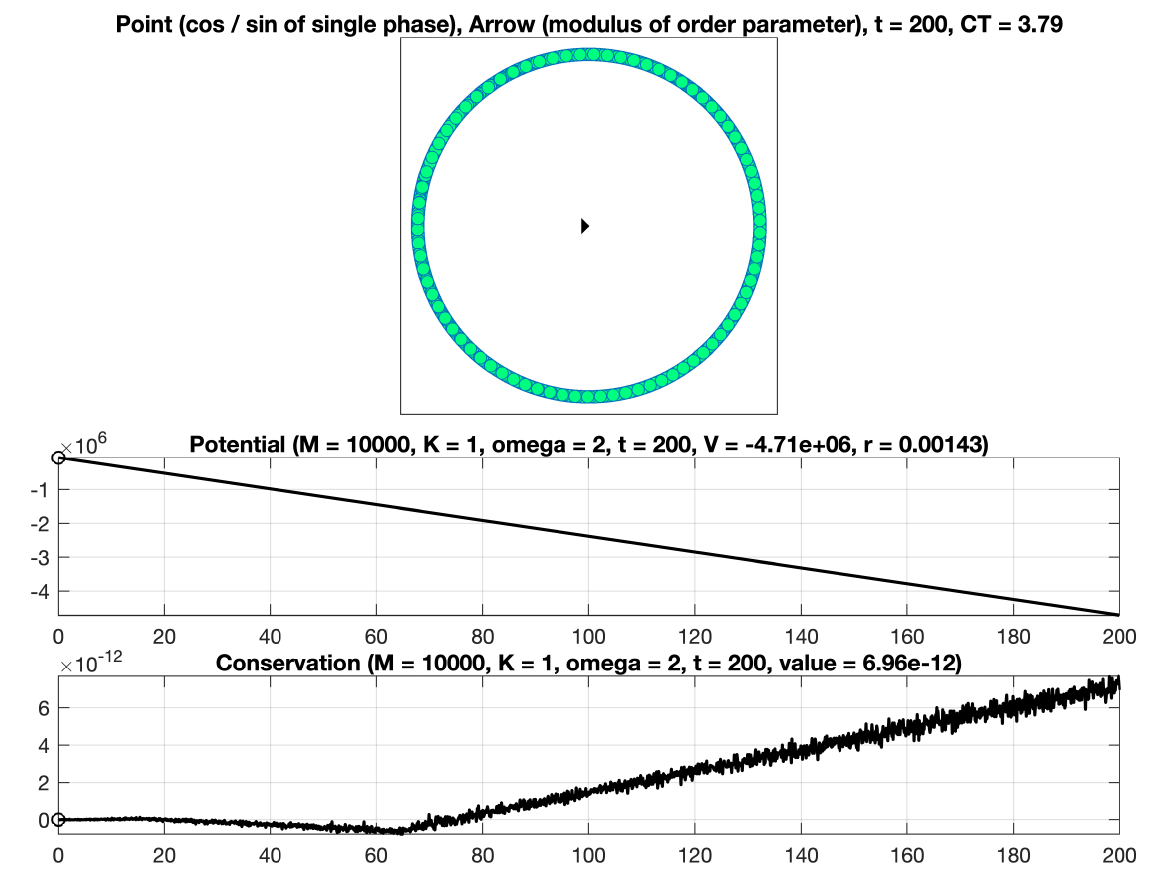}
\end{center}
\caption{\colblue{Numerical integration of the classical Kuramoto model involving $M = 10^4$ oscillators. Coupling constant $K = 1$ (no synchronisation).}}
\label{fig:MyFigure4}
\end{figure}

\begin{figure}[t!]
\begin{center}
\includegraphics[width=9cm]{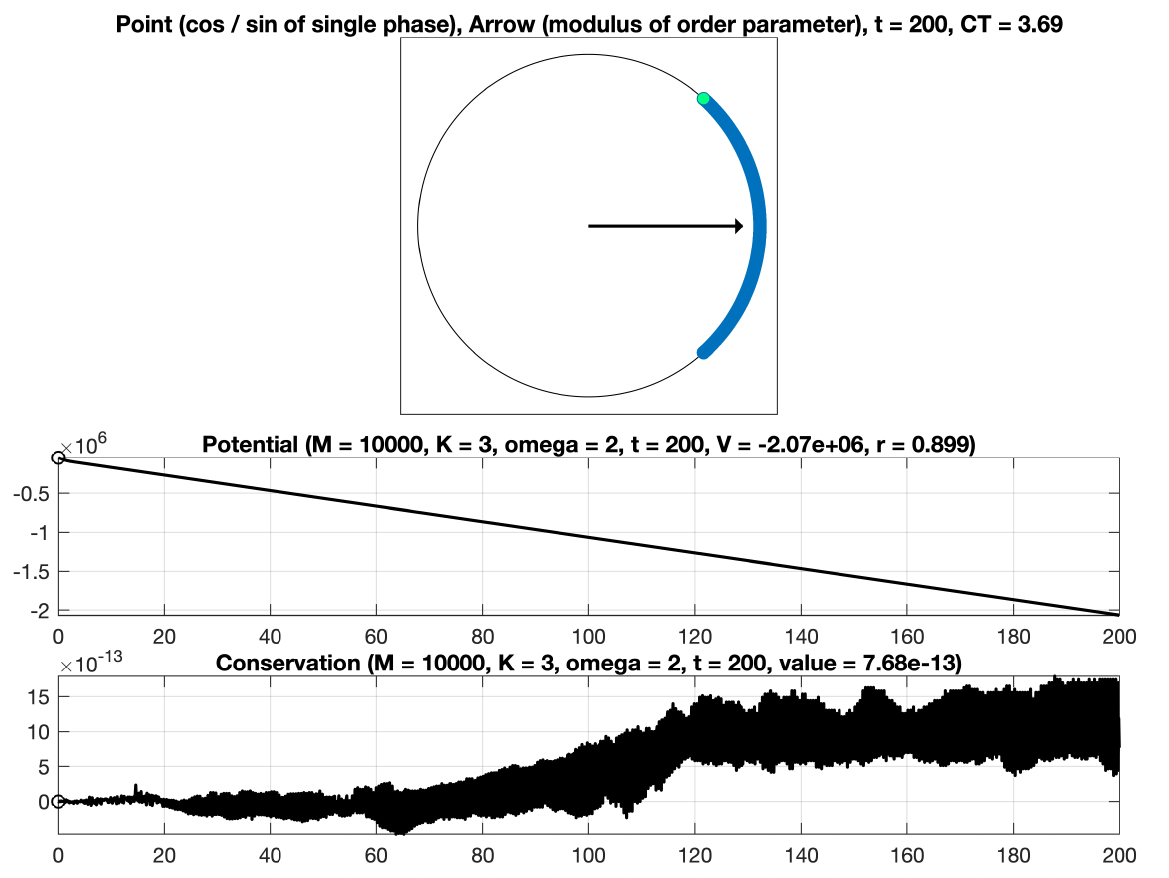} \\[4mm]
\includegraphics[width=9cm]{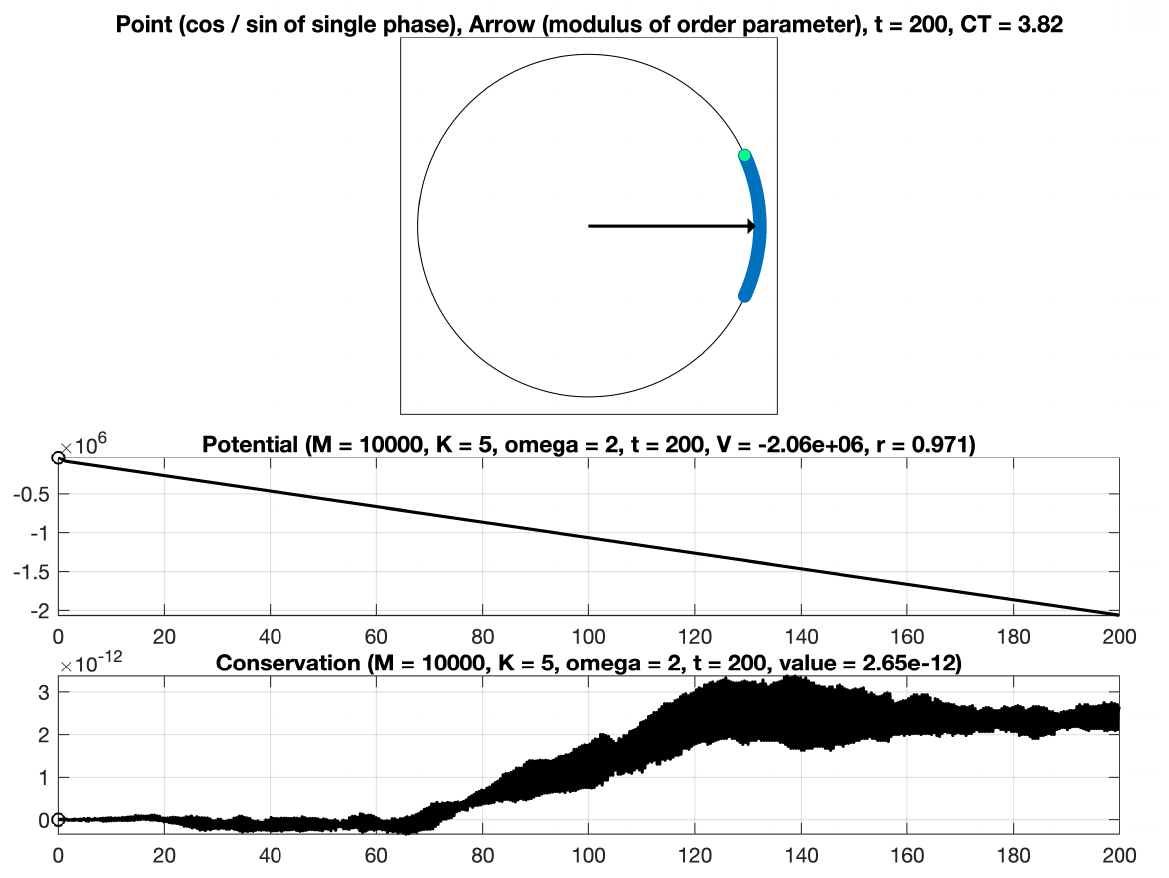} 
\end{center}
\caption{\colblue{Numerical integration of the classical Kuramoto model involving $M = 10^4$ oscillators. Coupling constant $K \in \{3, 5\}$ (gradual synchronisation).}}
\label{fig:MyFigure5}
\end{figure}

\begin{figure}[t!]
\begin{center}
\includegraphics[width=9cm]{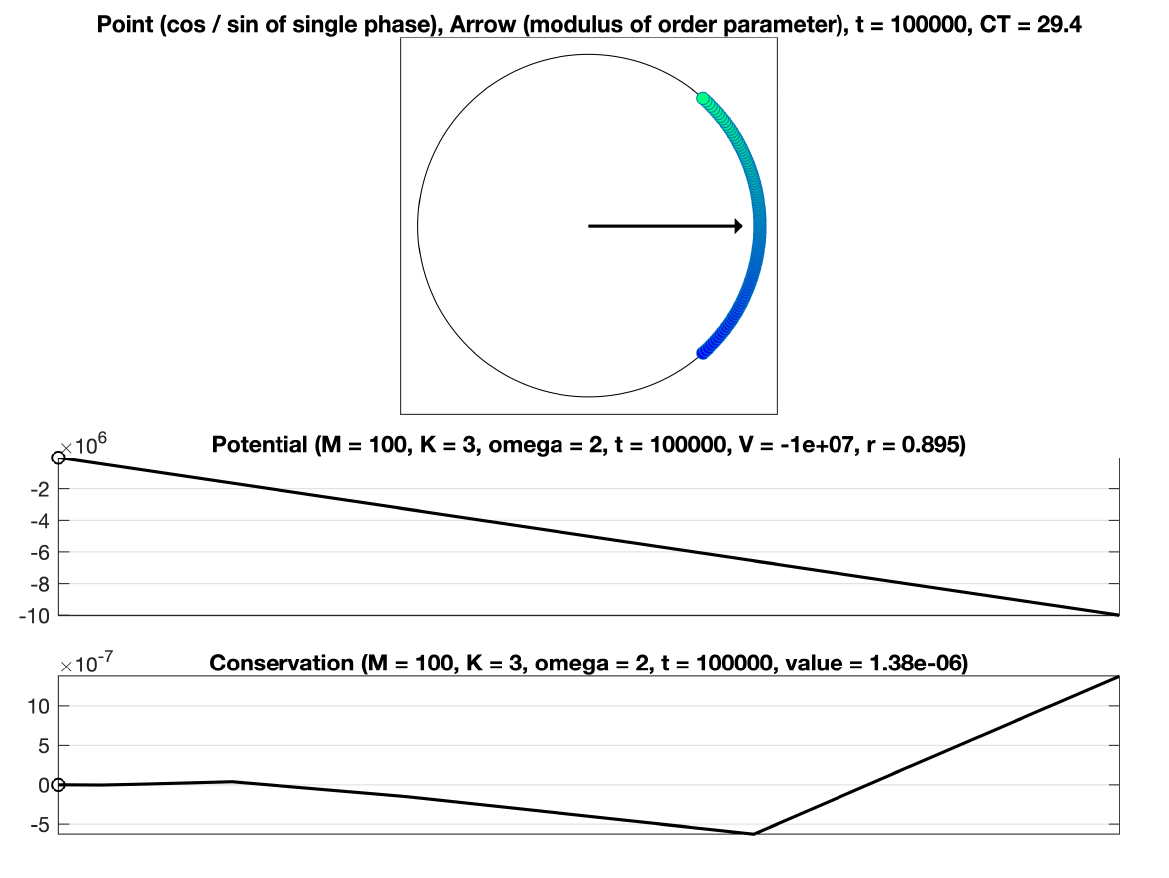} \\[4mm]
\includegraphics[width=9cm]{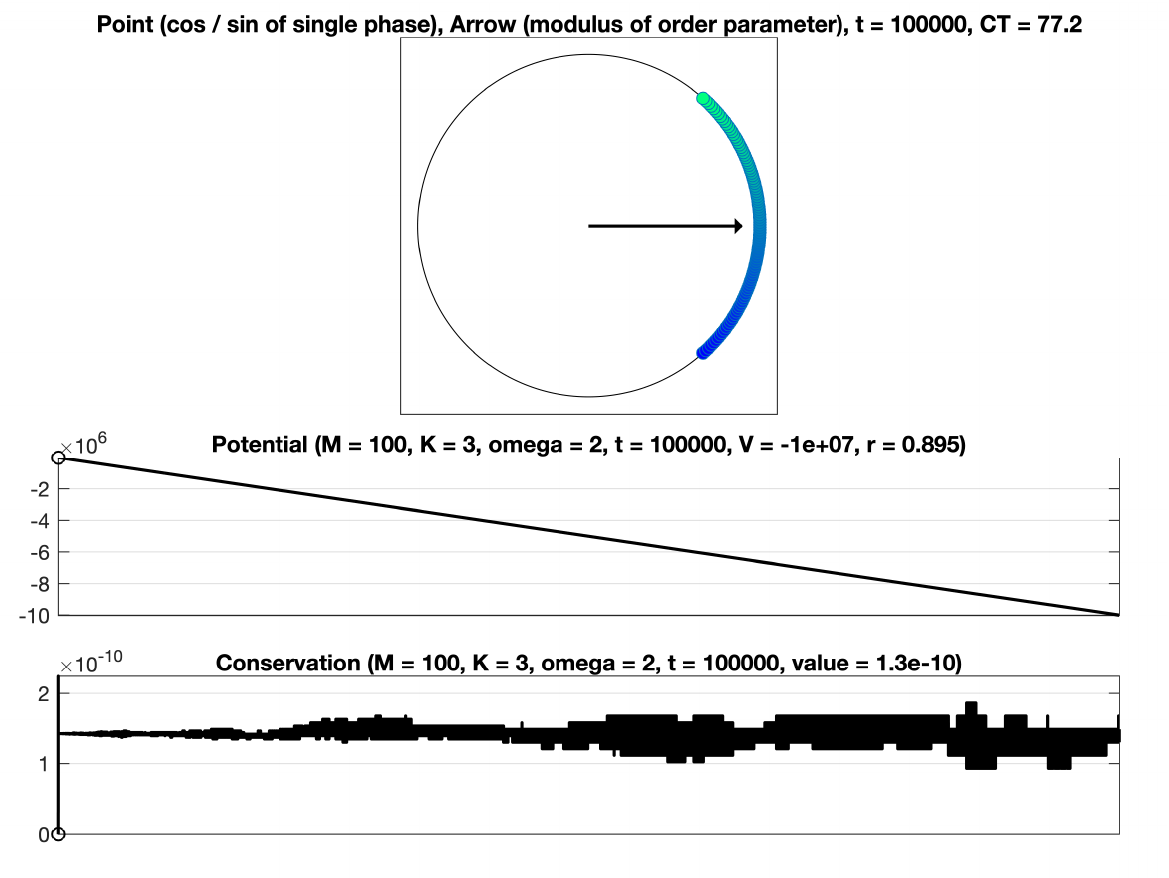} 
\end{center}
\caption{\colblue{Long-term integration of the classical Kuramoto model based on a second-order explicit Runge--Kutta method (top) and a second-order implicit Runge--Kutta method with improved numerical preservation of a conserved quantity (bottom). 
To avoid a further diminishment of the relevant vertical axis, ticks along the horizontal axis are omitted.}}
\label{fig:MyFigure6}
\end{figure}

\MyParagraphAbstand{Numerical integration}
The numerical integration of the classical Kuramoto model~\eqref{eq:Kuramoto}--\eqref{eq:KuramotoReformulation} based on the precomputation of sums is illustrated in Figures~\ref{fig:MyFigure2}--\ref{fig:MyFigure6}.
Movies showing the time evolution are found at 
\begin{center}
\small{\url{http://techmath.uibk.ac.at/mecht/MyHomepage/Research/MovieKuramotoClassicalK1.m4v}} \\
\small{\url{http://techmath.uibk.ac.at/mecht/MyHomepage/Research/MovieKuramotoClassicalK3.m4v}} \\
\small{\url{http://techmath.uibk.ac.at/mecht/MyHomepage/Research/MovieKuramotoClassicalK5.m4v}} \\
\end{center}
We consider $M = 10^2$ as well as $M = 10^4$ oscillators and set $T = 200$. 
The interplay between the coupling constant $K \in \{1, 3, 5\}$ and the constant $\omega_0 = 2$ defining the intrinsic frequencies determines the strength of synchronisation. 
\colred{
The points and the arrow reflect the values of the phases and the complex order parameter at the final time
\begin{equation*}
\Big(\MyCos\big(\vartheta_m(t) - \psi(t)\big), \MySin\big(\vartheta_m(t) - \psi(t)\big)\Big)\,, \quad m \in \{1, 2, \dots, M\}\,, \quad r\big(\vartheta(t)\big)\,, \quad t = T\,, 
\end{equation*}
see~\eqref{eq:OrderParameter}. 
Moreover, the graphs of the associated quantities
\begin{equation*}
V\big(\vartheta(t)\big)\,, \quad \tfrac{1}{M} \sum_{m=1}^{M} \big(\vartheta_m(t) - \vartheta_m(0) - t \, \omega_m\big)\,, \quad t \in [0, T]\,,
\end{equation*}
are shown, see~\eqref{eq:Potential} and \eqref{eq:Means}.}
In addition, the computation time CT, measured in seconds, is displayed.
For the considered time interval, a variable stepsize fourth-order explicit Runge--Kutta method leads to a reliable result.
In particular, it is seen that the values of the potential decrease when time evolves and that the modulus of the complex order parameter is close to one for $K = 5$.
When applying instead a second-order geometric integrator, we observed an improved behaviour regarding the conserved quantity over long times, see Figure~\ref{fig:MyFigure6}.
For comprehensive information about the benefits of geometric numerical integrators in comparison with standard integration methods, we refer to~\cite{BlanesCasas,IserlesQuispel,SanzSerna,SanzSernaCalvo}. 
\section{Kuramoto models on graphs}
\label{sec:SectionKuramotoGraph}
In this section, we investigate extensions of the classical Kuramoto model on graph topologies and propose suitable modifications of our approach for their efficient numerical integration.

\MyParagraphAbstand{General formulation}
Henceforth, we study the following system of ordinary differential equations 
\begin{subequations}
\label{eq:KuramotoGraph}
\begin{equation}
\begin{cases}
\displaystyle \vartheta_m'(t) = \omega_m + \tfrac{K}{\nM_m} \sum_{\ell=1}^{M} A_{m \ell} \, \MySin\big(\vartheta_{\ell}(t) - \vartheta_m(t)\big)\,, \\
\vartheta_m(0) \text{ given}\,, \quad t \in (0, T)\,, \quad m \in \{1, 2, \dots, M\}\,.
\end{cases}
\end{equation}
In contrast to the classical Kuramoto model~\eqref{eq:Kuramoto}, the all-to-all coupling is replaced by the interactions of certain communities of oscillators, which are described by the associated adjacency matrix
\begin{equation}
A = \big(A_{m \ell}\big)_{\ell,m \in \{1, \dots, M\}}\,, \quad A_{m \ell} \in \{0, 1\}\,, \quad \ell, m \in \{1, 2, \dots, M\}\,.
\end{equation}
We distinguish two kinds of scalings affecting in particular the strength of synchronisation.
In case the $m$th row of~$A$ is zero, the corresponding equation reduces to
\begin{equation*}
\vartheta_m'(t) = \omega_m\,, \quad t \in (0, T)\,, 
\end{equation*}
and hence the scaling is dispensible.
\begin{enumerate}[(i)]
\item
\emph{Uniform scaling.} \, 
With regard to the literature, a common choice is 
\begin{equation}
\nM_m = M\,, \quad m \in \{1, 2, \dots, M\}\,.
\end{equation}
\item
\emph{Non-uniform scaling.} \, 
As exemplified below, an alternative is 
\begin{equation}
\nM_m = \sum_{\ell=1}^{M} A_{m \ell}\,, \quad m \in \{1, 2, \dots, M\}\,.
\end{equation}
\end{enumerate}
For our considerations, it is convenient to introduce the scaled adjacency matrix
\begin{equation}
\nA = \big(\tfrac{1}{\nM_m} \, A_{m \ell}\big)_{\ell,m \in \{1, \dots, M\}}\,.
\end{equation}
\end{subequations}

\MyParagraphAbstand{Auxiliary identities}
The decisive terms defining the right-hand side of the extended Kuramoto model~\eqref{eq:KuramotoGraph} rewrite as follows 
\begin{equation*}
\begin{split}
&\sum_{\ell=1}^{M} \nA_{m \colblue{\ell}} \, \MySin\big(\vartheta_{\ell}(t) - \vartheta_m(t)\big) \\
&\quad= \sum_{\ell=1}^{M} \nA_{m \colblue{\ell}} \, \MySin\big(\vartheta_{\ell}(t)\big) \MyCos\big(\vartheta_m(t)\big) \\
&\quad\qquad - \sum_{\ell=1}^{M} \nA_{m \colblue{\ell}} \, \MyCos\big(\vartheta_{\ell}(t)\big) \MySin\big(\vartheta_m(t)\big)\,, \quad
t \in [0, T]\,, \quad m \in \{1, 2, \dots, M\}\,.
\end{split}
\end{equation*}
Denoting the componentwise product of two columns by 
\begin{equation*}
v \CompProduct w
= \begin{pmatrix} v_1 \\ v_2 \\ \vdots \\ v_M \end{pmatrix} \CompProduct \begin{pmatrix} w_1 \\ w_2 \\ \vdots \\ w_M \end{pmatrix}
= \begin{pmatrix} v_1 \, w_1 \\ v_2 \, w_2 \\ \vdots \\ v_M \, w_M \end{pmatrix}\,, \quad v, w \in \RR^M\,,
\end{equation*}
and employing the compact vector notation 
\begin{subequations}
\label{eq:KuramotoGraphReformulation}
\begin{equation}
\begin{gathered}
G\big(\vartheta(t)\big) = \omega + K \, \bigg(\Big(\nA \MySin\big(\vartheta(t)\big)\Big) \CompProduct \MyCos\big(\vartheta(t)\big)
- \Big(\nA \MyCos\big(\vartheta(t)\big)\Big) \CompProduct \MySin\big(\vartheta(t)\big)\,, \\
t \in [0, T]\,, 
\end{gathered}
\end{equation}
we obtain the following reformulation of the governing equations 
\begin{equation}
\begin{cases}
\vartheta'(t) = G\big(\vartheta(t)\big)\,, \quad t \in (0, T)\,, \\
\vartheta(0) \text{ given}\,,
\end{cases}
\end{equation}
\end{subequations}
see also~\eqref{eq:KuramotoReformulation}.
Similar considerations yield the identity 
\begin{equation*}
\begin{split}
&\sum_{\ell,m=1}^{M} \nA_{m \colblue{\ell}} \, \MyCos\big(\vartheta_{\ell}(t) - \vartheta_m(t)\big) \\
&\quad= \sum_{\ell,m=1}^{M} \nA_{m \colblue{\ell}} \, \MyCos\big(\vartheta_{\ell}(t)\big) \MyCos\big(\vartheta_m(t)\big)
+ \sum_{\ell,m=1}^{M} \nA_{m \colblue{\ell}} \, \MySin\big(\vartheta_{\ell}(t)\big) \MySin\big(\vartheta_m(t)\big) \\ 
&\quad= \Big(\MyCos\big(\vartheta(t)\big)\Big)^T \Big(\nA \MyCos\big(\vartheta(t)\big)\Big)
+ \Big(\MySin\big(\vartheta(t)\big)\Big)^T \Big(\nA \MySin\big(\vartheta(t)\big)\Big)\,, \quad t \in [0, T]\,.
\end{split}
\end{equation*}

\MyParagraphAbstand{Uniform scaling}
In situations, where the adjacency matrix is symmetric and the scaling is uniform
\begin{equation*}
A = A^T\,, \quad \nM_m = M\,, \quad m \in \{1, 2, \dots, M\}\,, \quad \nA = \tfrac{1}{M} A\,, \quad \nA^T = \nA\,,
\end{equation*}
the existence of a potential function is ensured and the sum over all phases leads to a conserved quantity, see also~\eqref{eq:PotentialOrderParameterConservation}.
With regard to the classical Kuramoto model and the above stated auxiliary identity, we define the potential through 
\begin{subequations}
\label{eq:ExtendedPotentialConservation}
\begin{equation}
\begin{split}
V: \RR^M &\longrightarrow \RR: \\
\vartheta
&\longmapsto 
V(\vartheta) = V(\vartheta_1, \dots, \vartheta_M) \\
&\qquad\qquad \,\,= - \sum_{m=1}^M \omega_m \, \vartheta_m + \tfrac{K}{2} \sum_{\ell,m=1}^{M} \nA_{m \colblue{\ell}} \, \big(1 - \cos(\vartheta_{\ell} - \vartheta_m)\big) \\
&\qquad\qquad \,\,= - \, \omega^T \vartheta + \tfrac{K}{2} \, \bigg(\sum_{\ell,m=1}^{M} \nA_{m \colblue{\ell}} - \big(\cos(\vartheta)\big)^T \big(\nA \cos(\vartheta)\big) \\
&\qquad\qquad\qquad\qquad\qquad\qquad\qquad\quad \;\;\;\, - \big(\sin(\vartheta)\big)^T \big(\nA \sin(\vartheta)\big)\bigg)\,.
\end{split}
\end{equation}
Due to the fact that the relation 
\begin{equation*}
\begin{gathered}
\nA_{m \ell} \, \MySin\big(\vartheta_{\ell}(t) - \vartheta_m(t)\big) + \nA_{\ell m} \, \MySin\big(\vartheta_m(t) - \vartheta_{\ell}(t)\big) = 0\,, \\
t \in [0, T]\,, \quad \ell, m \in \{1, 2, \dots, M\}\,,   
\end{gathered}
\end{equation*}
holds, summation and integration with respect to time implies 
\colred{
\begin{equation}
\label{eq:KuramotoGraphConservation}
\tfrac{1}{M} \sum_{m=1}^{M} \big(\vartheta_m(t) - \vartheta_m(0) - t \, \omega_m\big) = 0\,, \quad t \in [0, T]\,.
\end{equation}
}
\end{subequations}

\MyParagraphAbstand{Non-uniform scaling}
In order to illustrate our motive for the non-uniform scaling, we consider the special case, where the adjacency matrix comprises a square submatrix with coefficients equal to one and is zero otherwise 
\begin{equation*}
\begin{gathered}
A = \begin{pmatrix} B^{(11)} & B^{(12)} \\ B^{(21)} & B^{(22)} \end{pmatrix} \in \RR^{M \times M}\,, \\
B^{(11)} = \big(1\big)_{\ell,m \in \{1, \dots, M_0\}}\,, \quad B^{(12)} = B^{(21)} = B^{(22)} = 0 \in \RR^{M_0 \times M_0}\,, \\
M_0 \in \{2, \dots, M-1\}\,. 
\end{gathered}
\end{equation*}
This would correspond to a configuration with pairwise interaction of the first part of the oscillators and a decoupling of the second part 
\begin{equation*}
\begin{cases}
\displaystyle \vartheta_m'(t) = \omega_m + \tfrac{K}{\nM_m} \sum_{\ell=1}^{M_0} \MySin\big(\vartheta_{\ell}(t) - \vartheta_m(t)\big)\,,
\quad t \in (0, T)\,, \quad m \in \{1, 2, \dots, M_0\}\,, \\
\displaystyle \vartheta_m'(t) = \omega_m\,, \quad t \in (0, T)\,, \quad m \in \{M_0 + 1, M_0 + 2, \dots, M\}\,.
\end{cases}
\end{equation*}
Here, the non-uniform scaling seems to be more natural
\begin{equation*}
\nM_m = \sum_{\ell=1}^{M} A_{m \ell} = M_0 \neq M\,, \quad m \in \{1, 2, \dots, M_0\}\,.
\end{equation*}
We point out that a potential function of the form~\eqref{eq:ExtendedPotentialConservation} remains valid for non-symmetric scaled adjacency matrices, whereas summation over all governing equations \colblue{does not lead} to a conserved quantity, in general. 

\MyParagraphAbstand{Extension of our approach}
The extension of our approach for the classical Kuramoto model based on the precomputation of sums to Kuramoto models on graphs~\eqref{eq:KuramotoGraph} requires a careful incorporation of the structure of the associated adjacency matrix.
The following considerations prove to be particularly expedient for situations, where the adjacency matrix has a block structure and can be divided into relatively sparse and relatively dense submatrices, respectively. 
In view of an efficient evaluation of the decisive terms in~\eqref{eq:KuramotoGraph}, we consider a submatrix of the form 
\begin{equation*}
\big(A_{m \ell}\big)_{m \in \{M_1, \dots, M_2\}, \; \ell \in \{M_3, \dots, M_4\}}\,, 
\end{equation*}
defined by positive integers $M_1, M_2, M_3, M_4 \in \NN$ such that $1 \leq M_1 < M_2 \leq M$ as well as $1 \leq M_3 < M_4 \leq M$.
Our basic concept is to optimise the required memory capacity and the number of functions evaluations. 
More precisely, in order to avoid the storage of each submatrix and to accelerate the computation of the sums
\begin{equation*}
\sum_{\ell=M_3}^{M_4} A_{m \ell} \, \MySin\big(\vartheta_{\ell}(t) - \vartheta_m(t)\big)\,, \quad m \in \{M_1, \dots, M_2\}\,, 
\end{equation*}
we distinguish two complementary cases.
Whenever the number of coefficients equal to zero is relatively high, we store all indices corresponding to non-zero coefficients and sum over these coefficients
\begin{subequations}
\label{eq:KuramotoGraphApproach}
\begin{equation}
\underset{A_{m \colblue{\ell}} = 1}{\sum_{\ell \in \{M_3, \dots, M_4\}}} A_{m \ell} \, \MySin\big(\vartheta_{\ell}(t) - \vartheta_m(t)\big)\,, \quad m \in \{M_1, \dots, M_2\}\,.
\end{equation}
Whenever the number of non-zero coefficients is relatively high, we instead store all indices corresponding to coefficients equal to zero and make use of the precomputation of sums.
That is, where applicable, we employ the reformulation  
\begin{equation*}
\begin{split}
&\sum_{\ell=M_3}^{M_4} A_{m \ell} \, \MySin\big(\vartheta_{\ell}(t) - \vartheta_m(t)\big) \\
&\quad = \sum_{\ell=M_3}^{M_4} A_{m \ell} \, \MySin\big(\vartheta_{\ell}(t)\big) \MyCos\big(\vartheta_m(t)\big)
- \sum_{\ell=M_3}^{M_4} A_{m \ell} \, \MyCos\big(\vartheta_{\ell}(t)\big) \MySin\big(\vartheta_m(t)\big) \\
&\quad = \Bigg(\sum_{\ell=M_3}^{M_4} \MySin\big(\vartheta_{\ell}(t)\big)
- \sum_{\ell=M_3}^{M_4} \big(1 - A_{m \ell}\big) \MySin\big(\vartheta_{\ell}(t)\big)\Bigg) \MyCos\big(\vartheta_m(t)\big) \\
&\quad\qquad - \Bigg(\sum_{\ell=M_3}^{M_4} \MyCos\big(\vartheta_{\ell}(t)\big)
- \sum_{\ell=M_3}^{M_4} \big(1 - A_{m \ell}\big) \MyCos\big(\vartheta_{\ell}(t)\big)\Bigg) \MySin\big(\vartheta_m(t)\big)\,, 
\end{split}
\end{equation*}
precompute the sums 
\begin{equation}
S_{M_3,M_4}\big(\vartheta(t)\big) = \sum_{\ell=M_3}^{M_4} \MySin\big(\vartheta_{\ell}(t)\big)\,, \quad 
C_{M_3,M_4}\big(\vartheta(t)\big) = \sum_{\ell=M_3}^{M_4} \MyCos\big(\vartheta_{\ell}(t)\big)\,, 
\end{equation}
and then substract the terms that correspond to coefficients equal to zero
\colred{
\begin{equation}
\begin{split}
&\sum_{\ell=M_3}^{M_4} A_{m \ell} \, \MySin\big(\vartheta_{\ell}(t) - \vartheta_m(t)\big) \\
&\quad= \Bigg(S_{M_3,M_4}\big(\vartheta(t)\big) - \underset{A_{m \colblue{\ell}} = 0}{\sum_{\ell \in \{M_3, \dots, M_4\}}} \MySin\big(\vartheta_{\ell}(t)\big)\Bigg) \MyCos\big(\vartheta_m(t)\big) \\
&\quad\qquad - \Bigg(C_{M_3,M_4}\big(\vartheta(t)\big) \\
&\quad\qquad - \underset{A_{m \colblue{\ell}} = 0}{\sum_{\ell \in \{M_3, \dots, M_4\}}} \MyCos\big(\vartheta_{\ell}(t)\big)\Bigg) \MySin\big(\vartheta_m(t)\big)\,, \quad m \in \{M_1, \dots, M_2\}\,.
\end{split}
\end{equation}
}
\end{subequations}

\MyParagraphAbstand{Remark}
In situations, where it is desirable to evaluate the right-hand side and the potential in parallel, it is advantageous to modify the above stated approach in view of an efficient computation of
\begin{equation*}
\begin{gathered}
\sum_{\ell=M_3}^{M_4} A_{m \ell} \, \MySin\big(\vartheta_{\ell}(t)\big) = \sum_{\ell=M_3}^{M_4} \MySin\big(\vartheta_{\ell}(t)\big)
- \sum_{\ell=M_3}^{M_4} \big(1 - A_{m \ell}\big) \MySin\big(\vartheta_{\ell}(t)\big) \\
\sum_{\ell=M_3}^{M_4} A_{m \ell} \, \MyCos\big(\vartheta_{\ell}(t)\big) = \sum_{\ell=M_3}^{M_4} \MyCos\big(\vartheta_{\ell}(t)\big)
- \sum_{\ell=M_3}^{M_4} \big(1 - A_{m \ell}\big) \MyCos\big(\vartheta_{\ell}(t)\big)\,, \\
m \in \{M_1, \dots, M_2\}\,,
\end{gathered}
\end{equation*}
either by straightforward summation or based on the precomputation of sums, see also~\eqref{eq:KuramotoGraphReformulation}--\eqref{eq:ExtendedPotentialConservation}.

\begin{figure}[t!]
\begin{center}
\includegraphics[width=3.4cm]{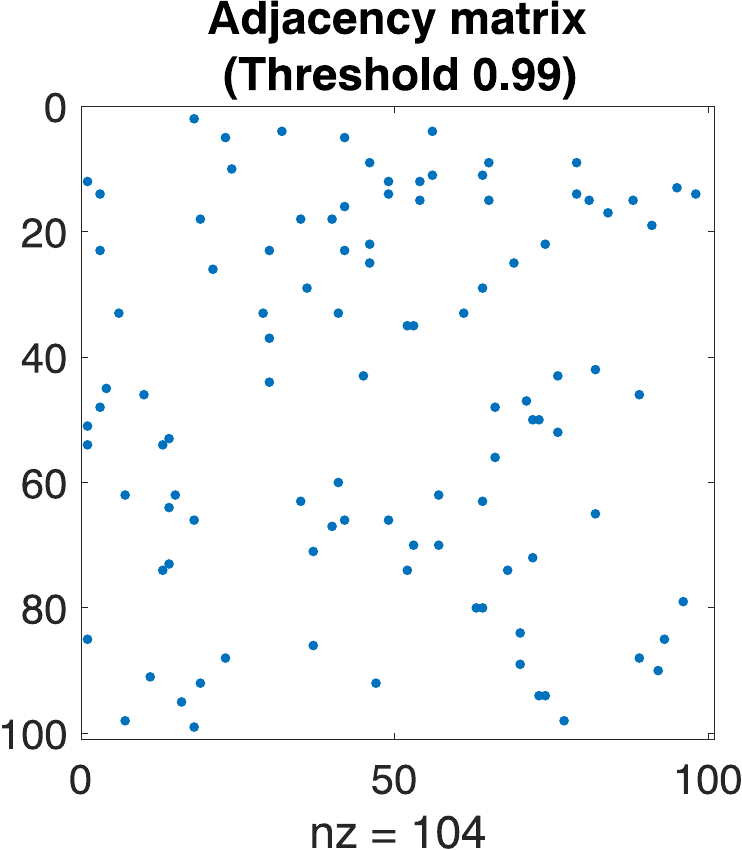} \quad
\includegraphics[width=4.9cm]{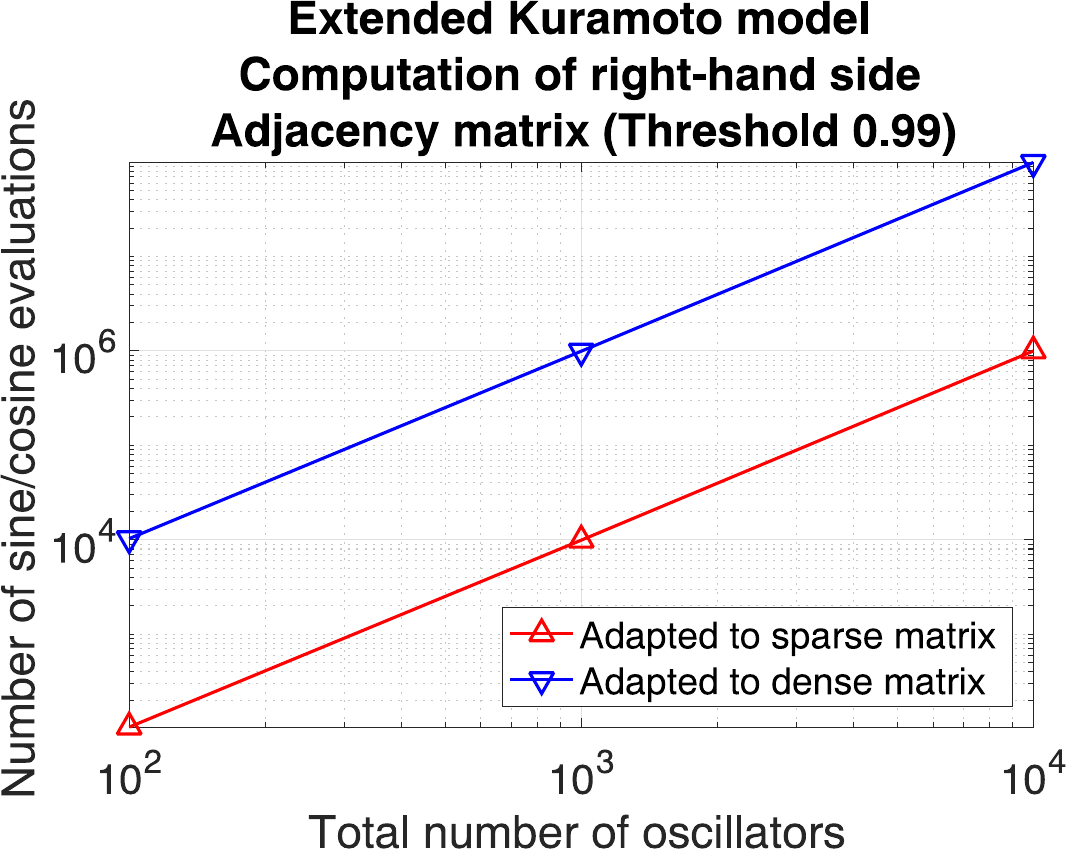} \quad
\includegraphics[width=4.9cm]{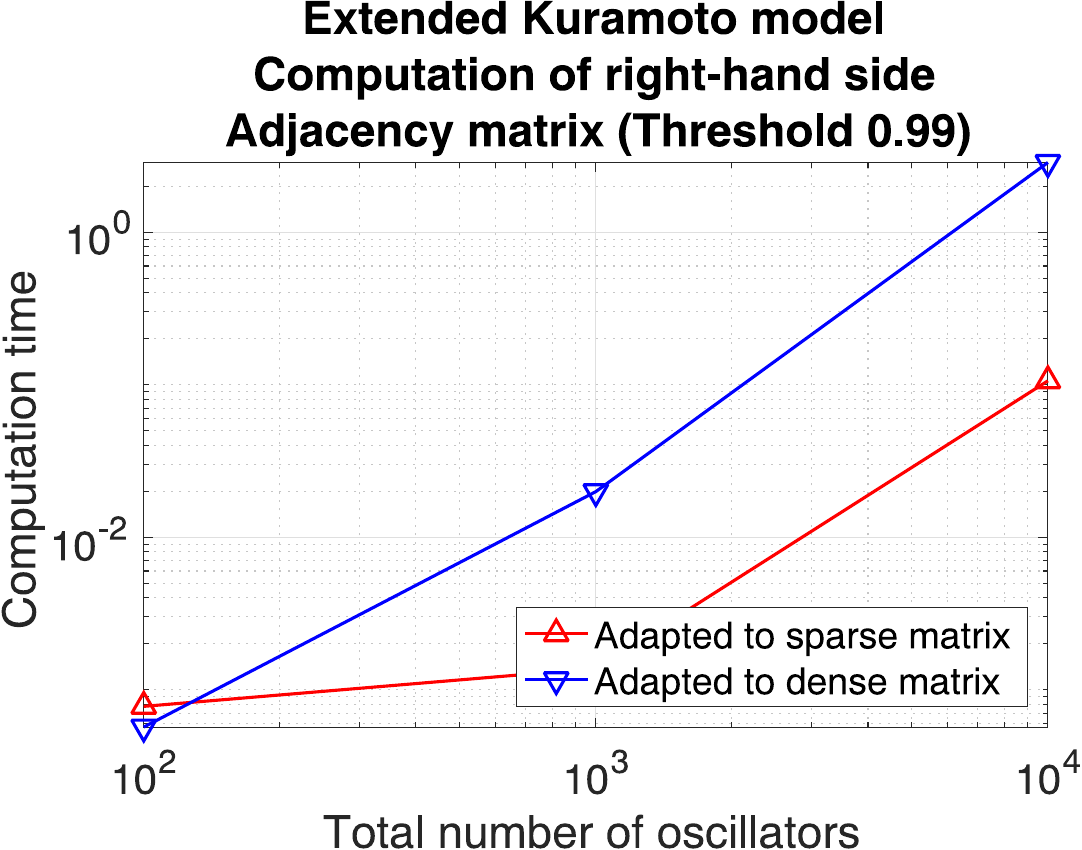} \\[4mm]
\includegraphics[width=3.4cm]{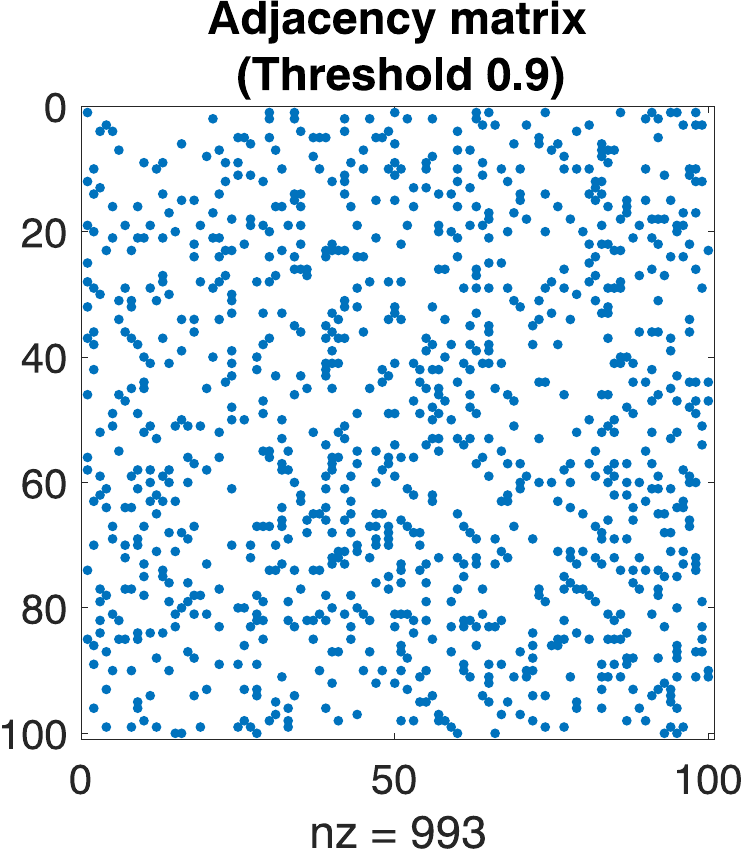} \quad
\includegraphics[width=4.9cm]{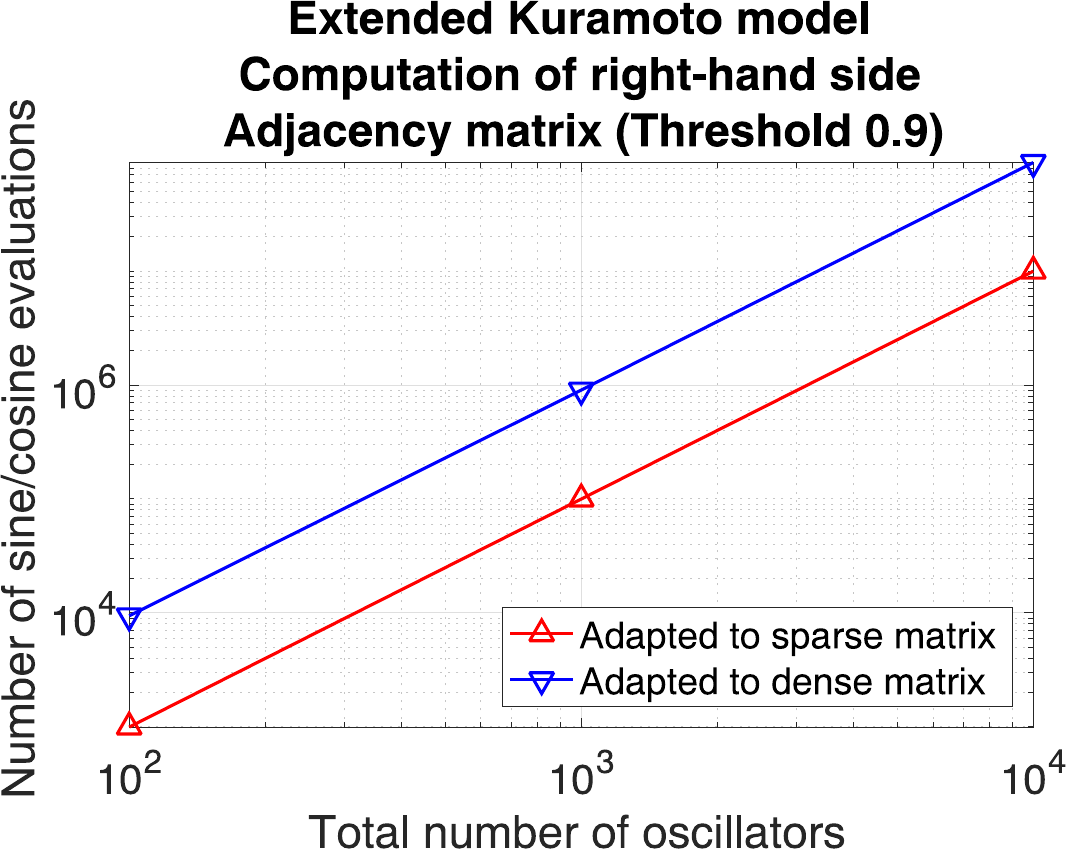} \quad
\includegraphics[width=4.9cm]{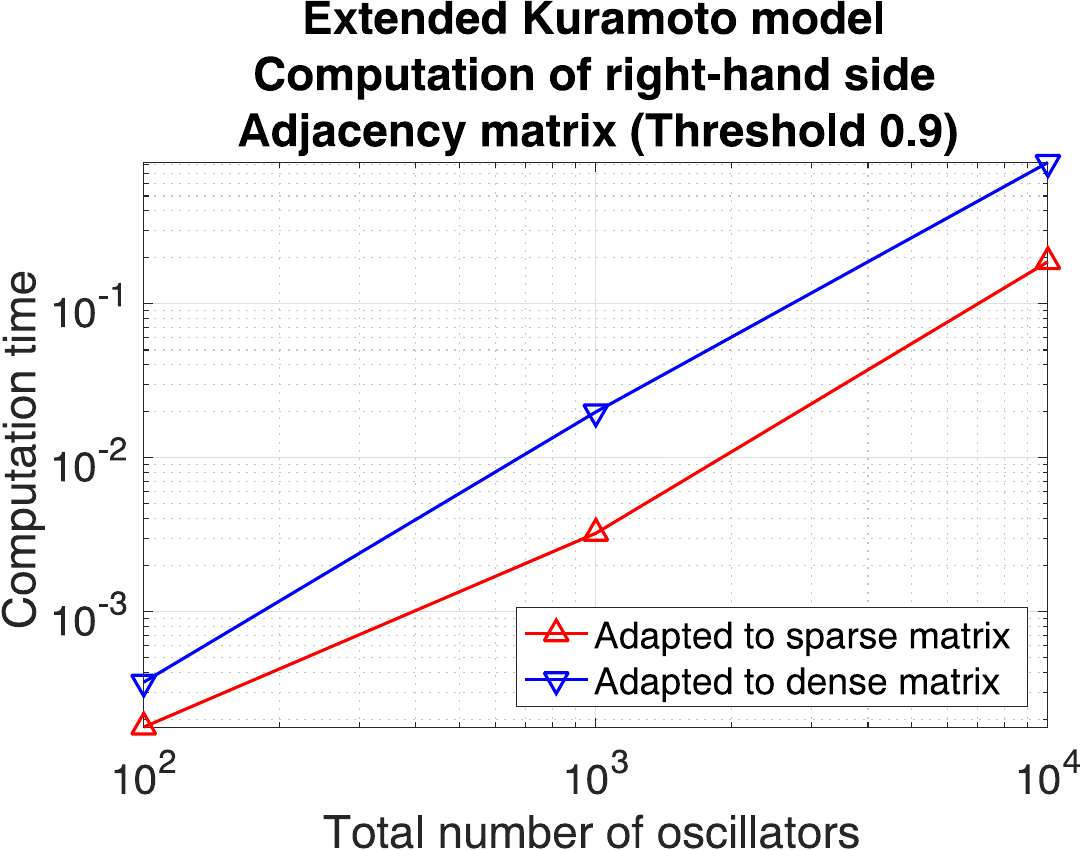}
\end{center}
\caption{\colblue{Extended Kuramoto models involving randomly generated adjacency matrices defined through the thresholds $0.99, 0.9$. 
Evaluation of the right-hand side by means of approaches adapted to sparse matrices (straightforward summation) and dense matrices (precomputation of sums), respectively.
Left: Illustration of the adjacency matrix for 100 oscillators. 
Middle: Number of sine and cosine evaluations versus the total number of oscillators.
Right: Comparison of the computation time.}}
\label{fig:MyFigure7}
\end{figure}

\begin{figure}[t!]
\begin{center}
\includegraphics[width=3.4cm]{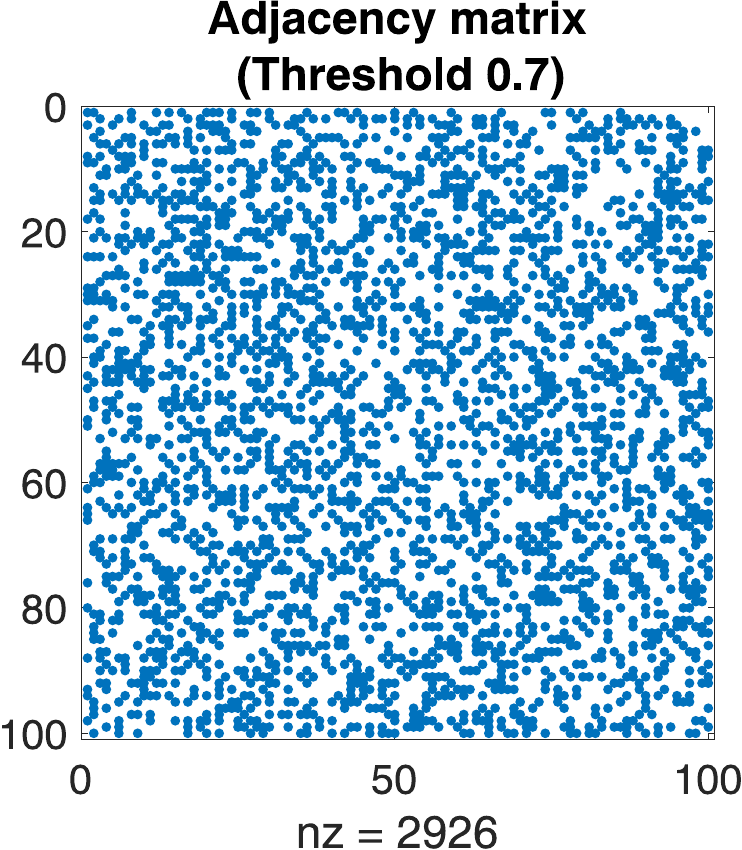} \quad
\includegraphics[width=4.9cm]{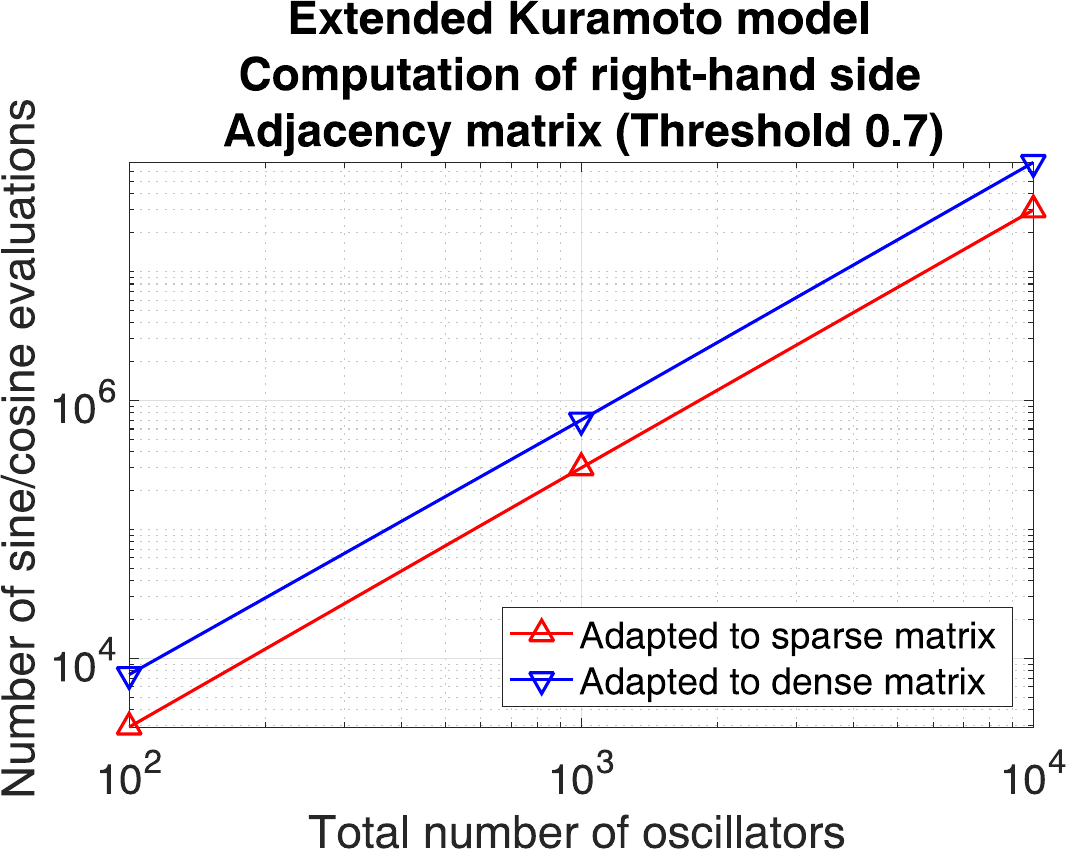} \quad
\includegraphics[width=4.9cm]{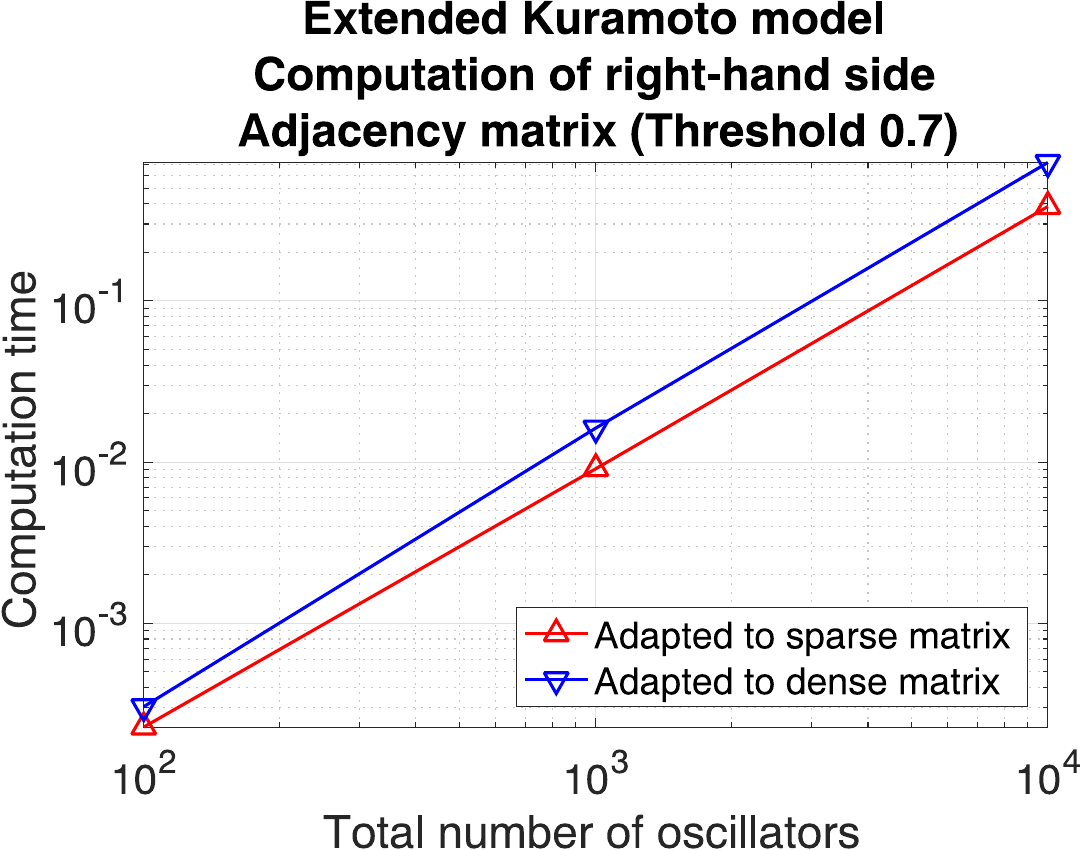} \\[4mm]
\includegraphics[width=3.4cm]{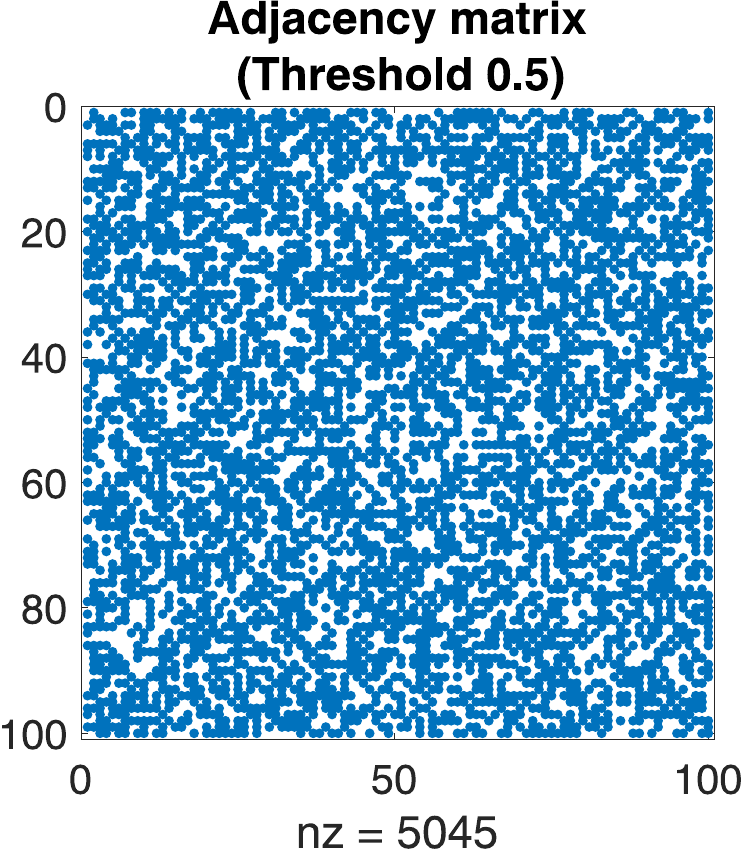} \quad
\includegraphics[width=4.9cm]{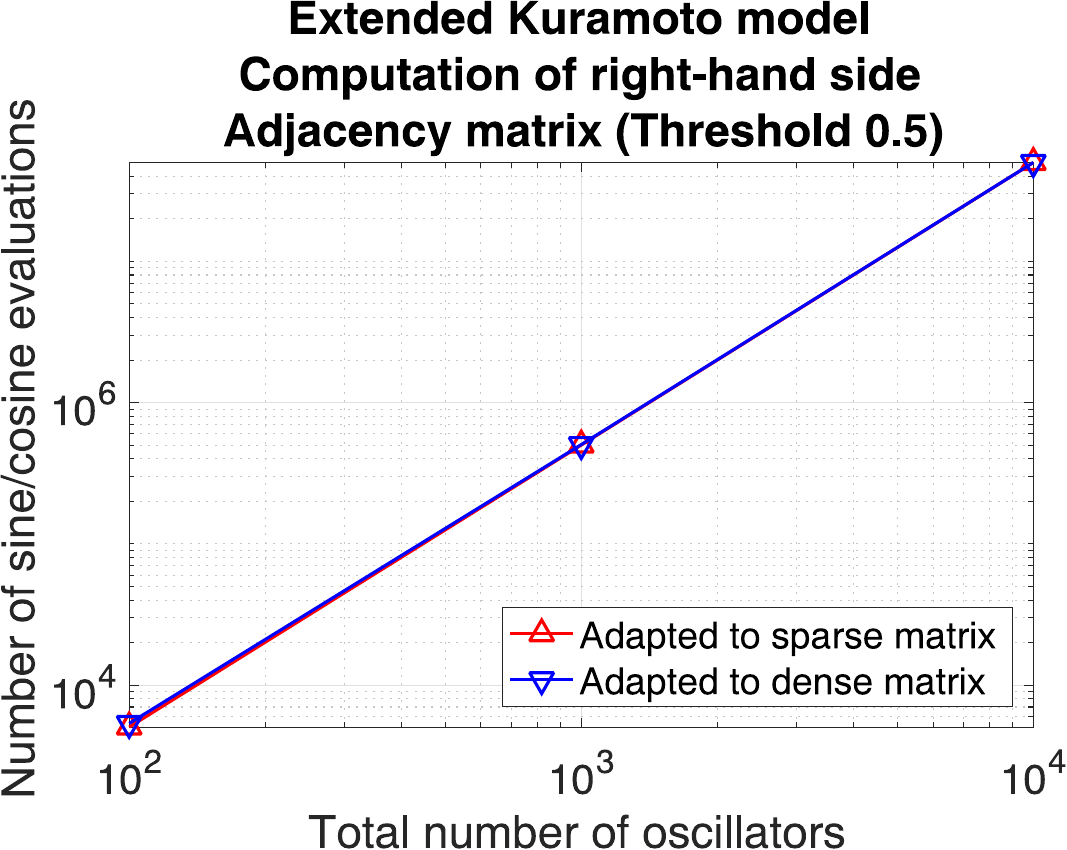} \quad
\includegraphics[width=4.9cm]{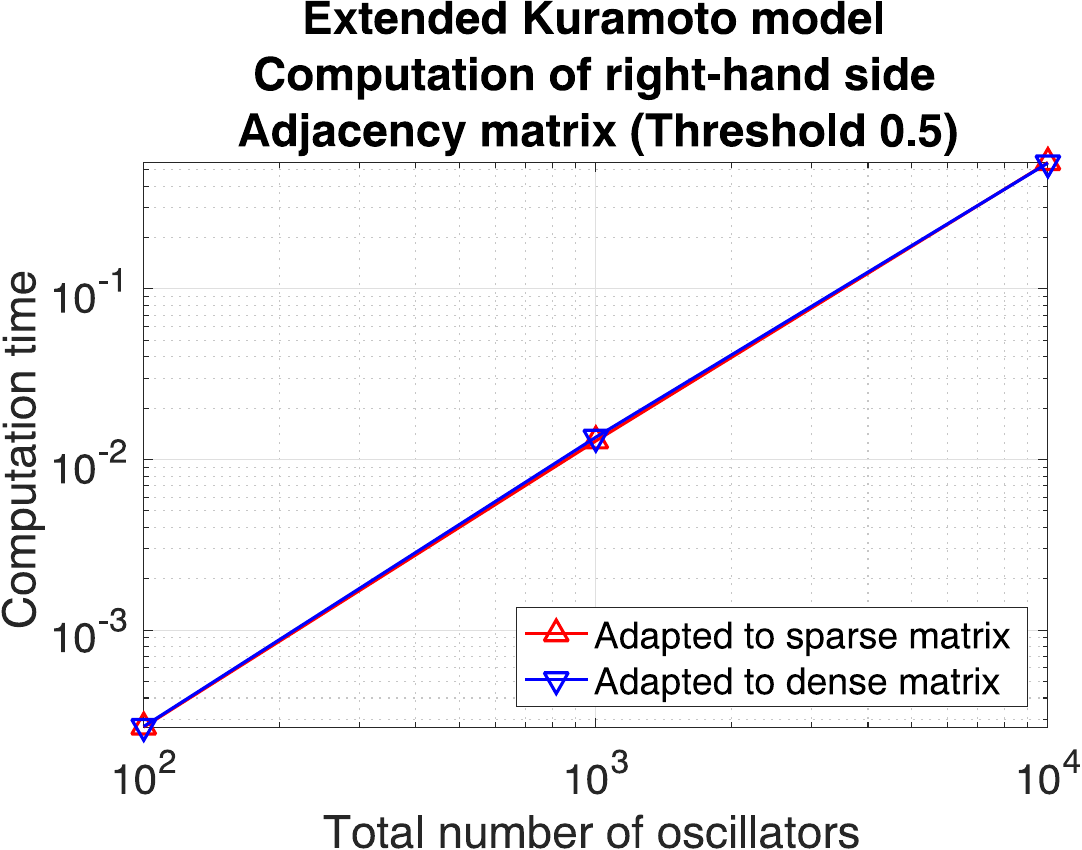} \\[4mm]
\includegraphics[width=3.4cm]{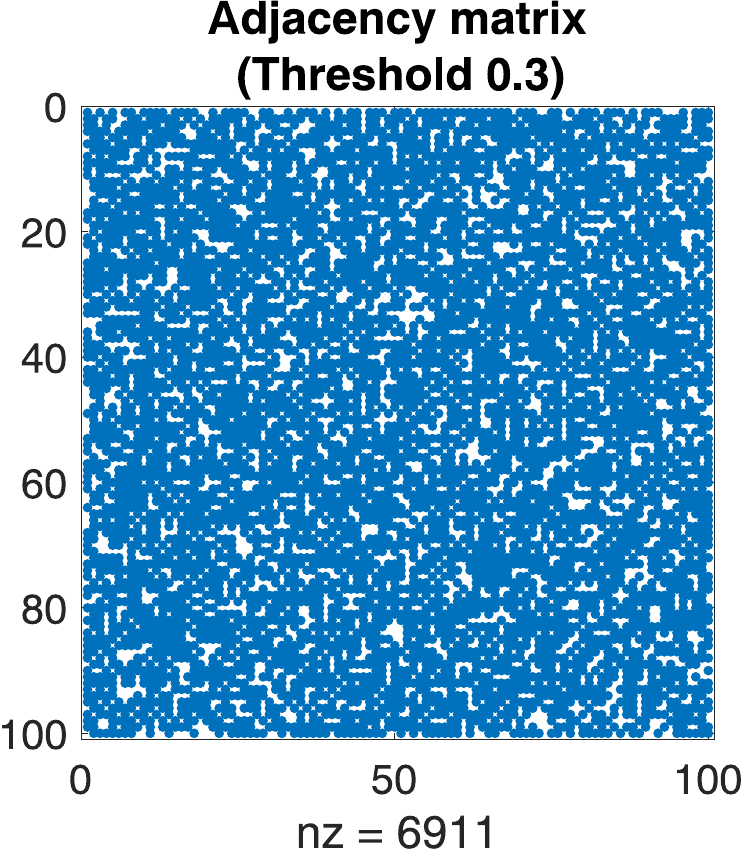} \quad
\includegraphics[width=4.9cm]{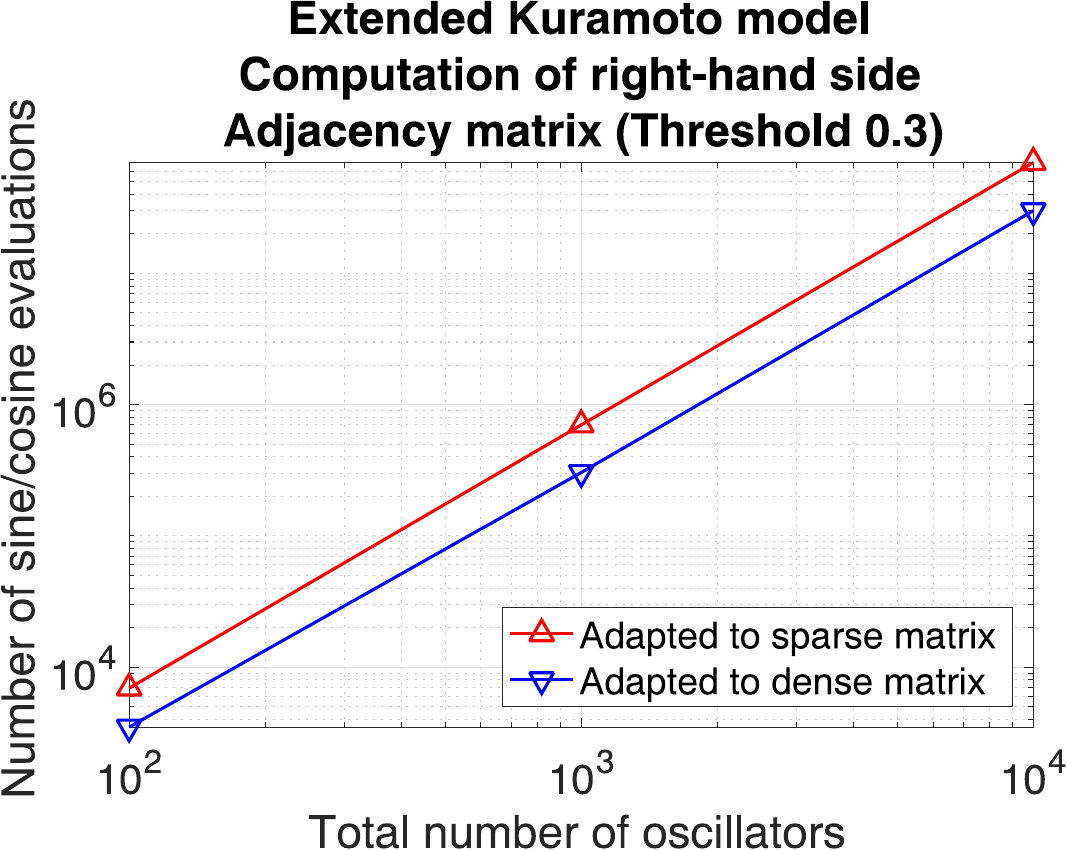} \quad
\includegraphics[width=4.9cm]{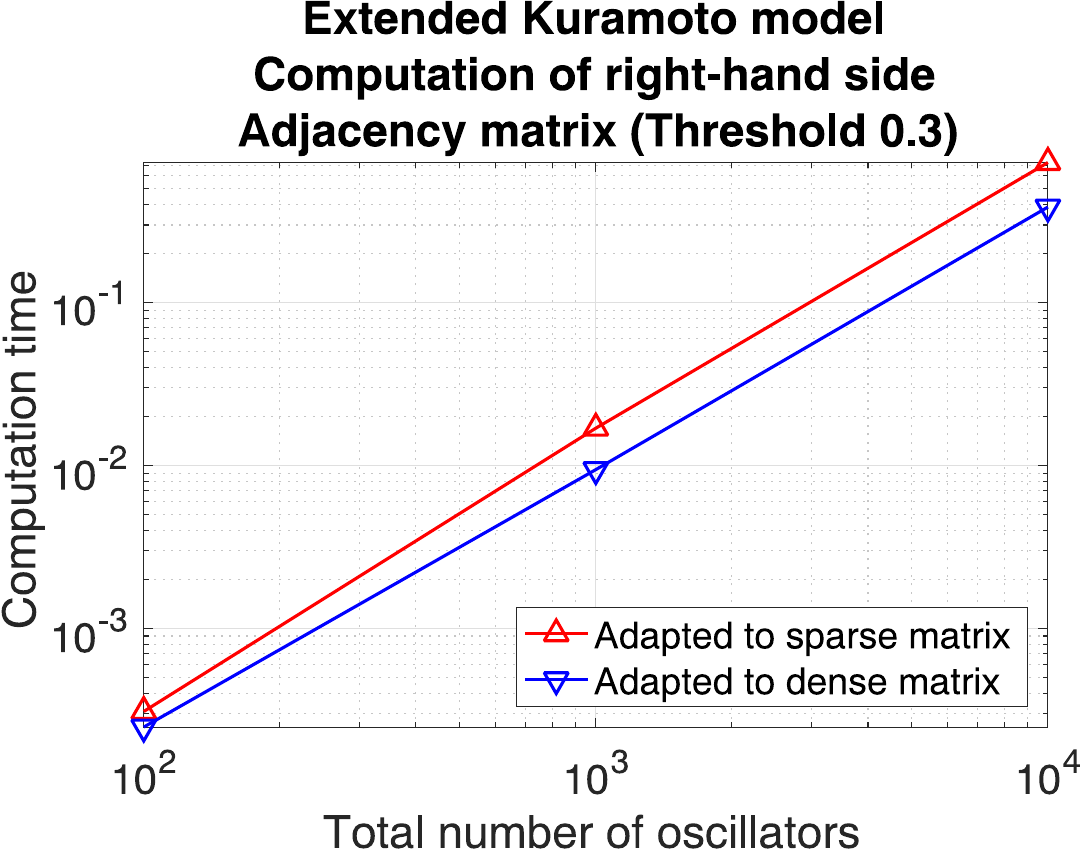} 
\end{center}
\caption{\colblue{Corresponding results for the thresholds $0.7, 0.5, 0.3$.}}
\label{fig:MyFigure8}
\end{figure}

\begin{figure}[t!]
\begin{center}
\includegraphics[width=3.4cm]{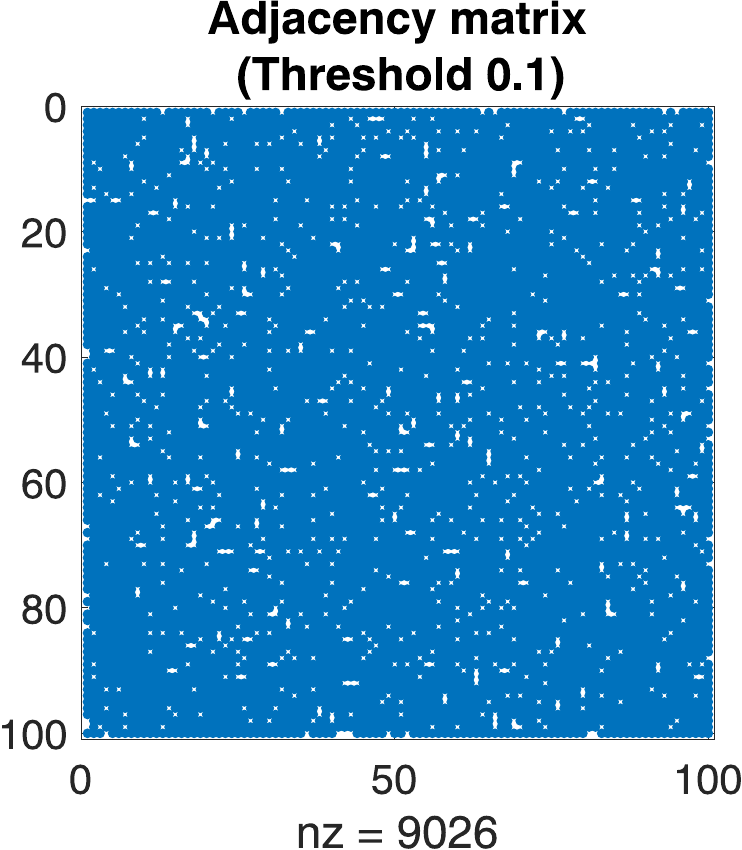} \quad
\includegraphics[width=4.9cm]{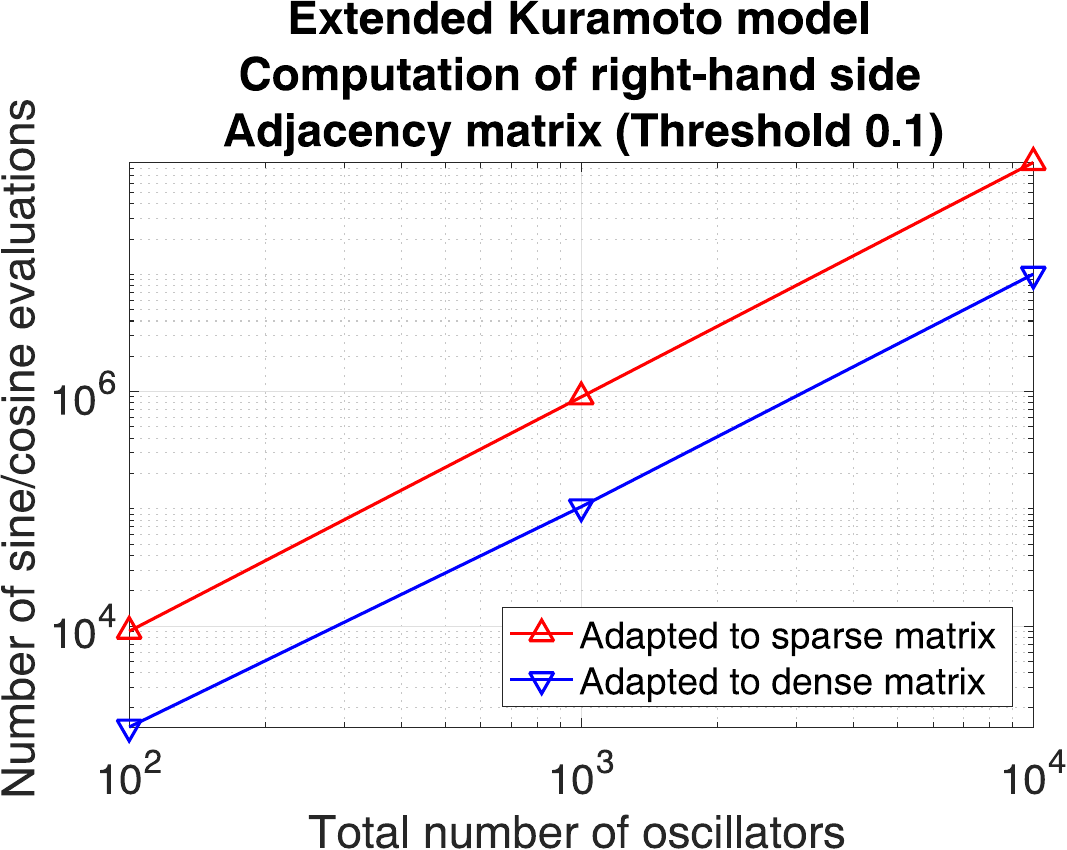} \quad
\includegraphics[width=4.9cm]{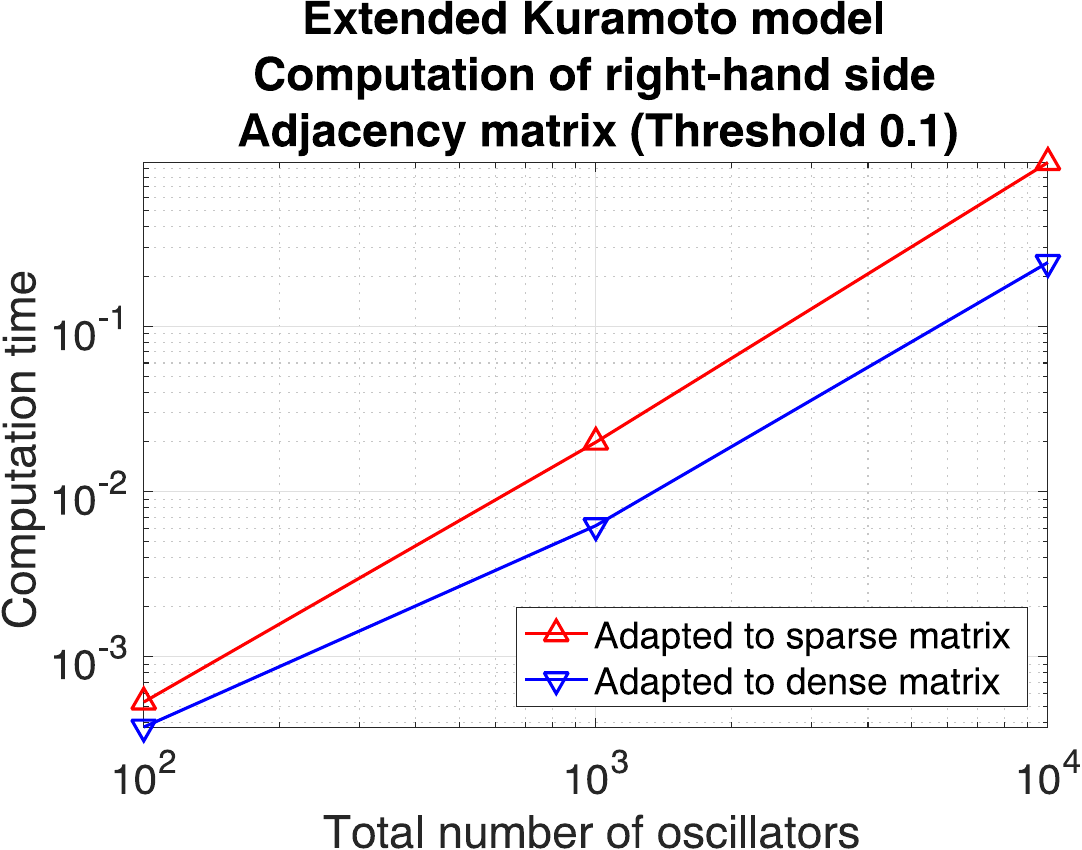} \\[4mm]
\includegraphics[width=3.4cm]{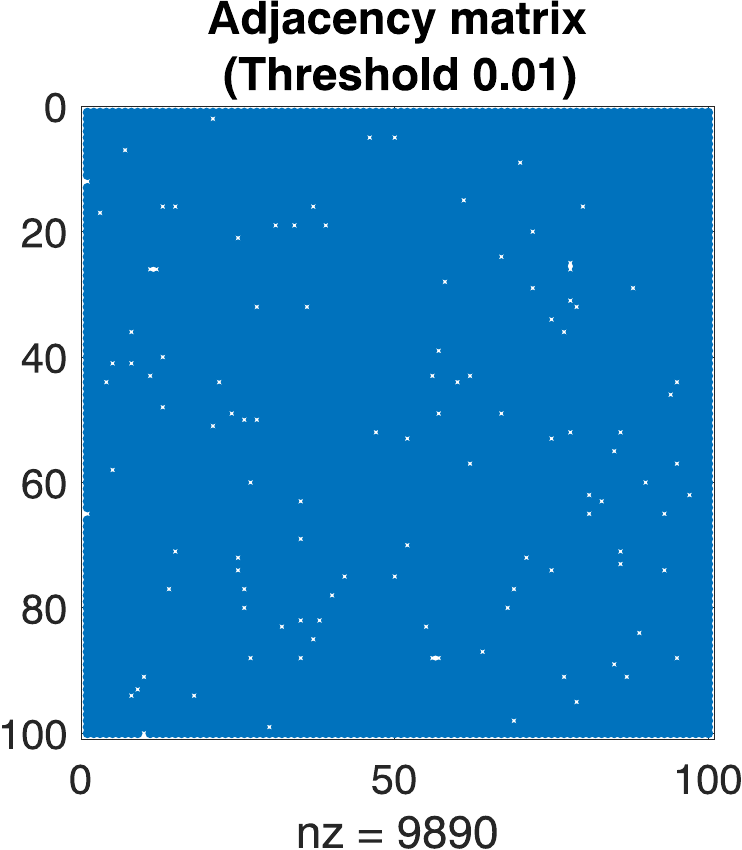} \quad
\includegraphics[width=4.9cm]{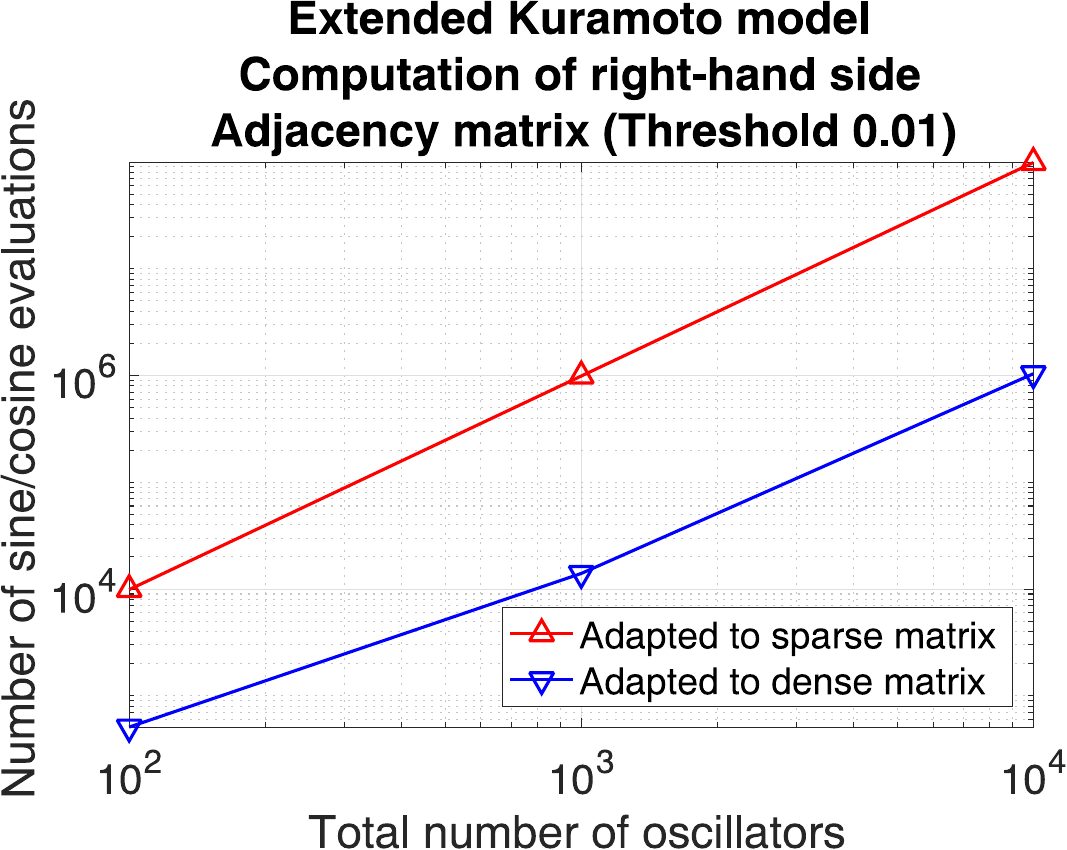} \quad
\includegraphics[width=4.9cm]{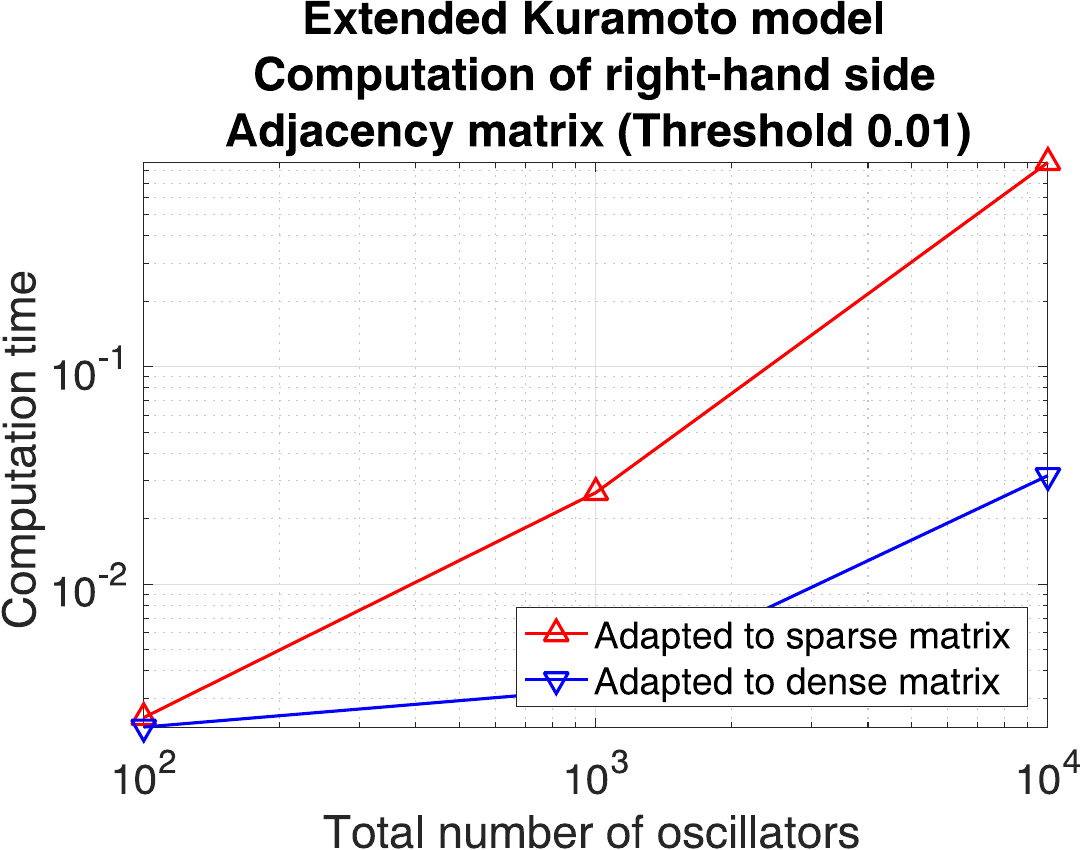} 
\end{center}
\caption{\colblue{Corresponding results for the thresholds $0.1, 0.01$.}}
\label{fig:MyFigure9}
\end{figure}

\begin{figure}[t!]
\begin{center}
\includegraphics[width=3.2cm]{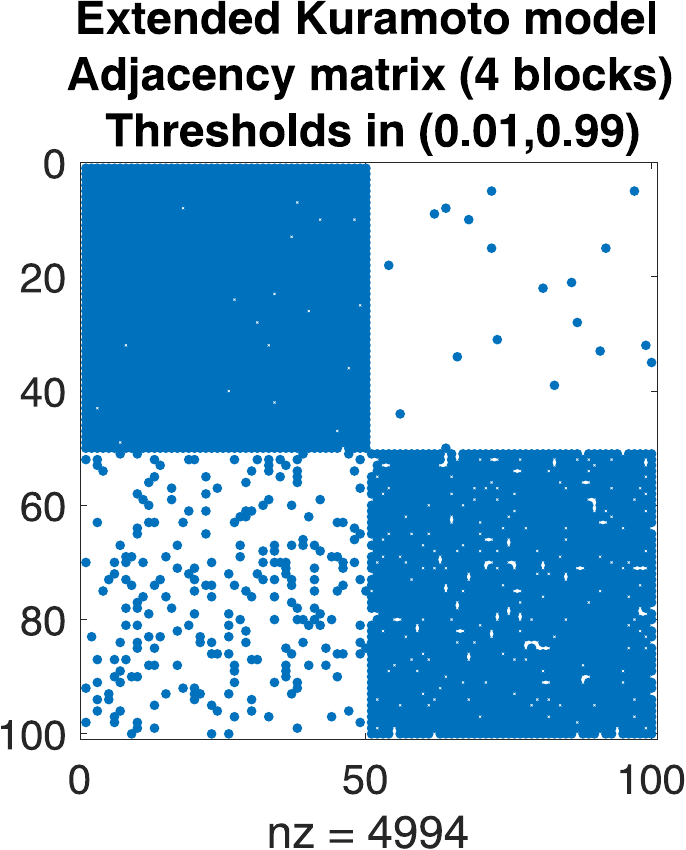} \quad
\includegraphics[width=5cm]{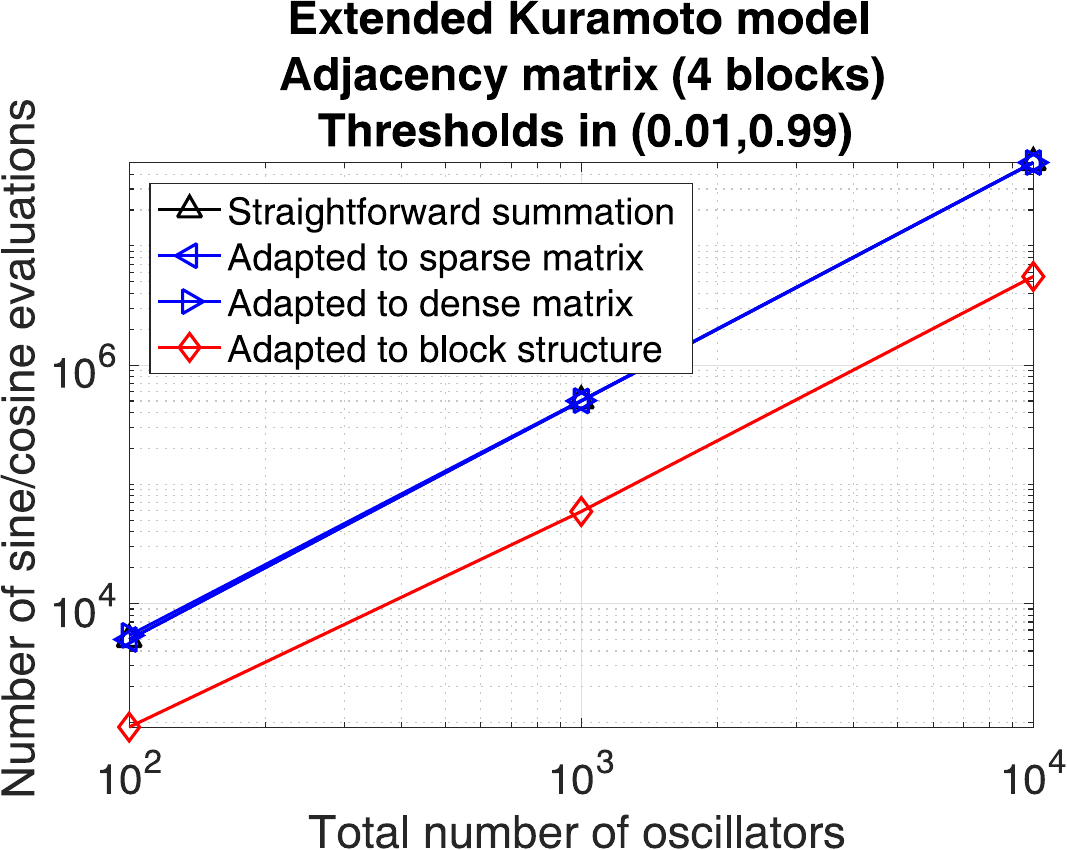} \quad
\includegraphics[width=5cm]{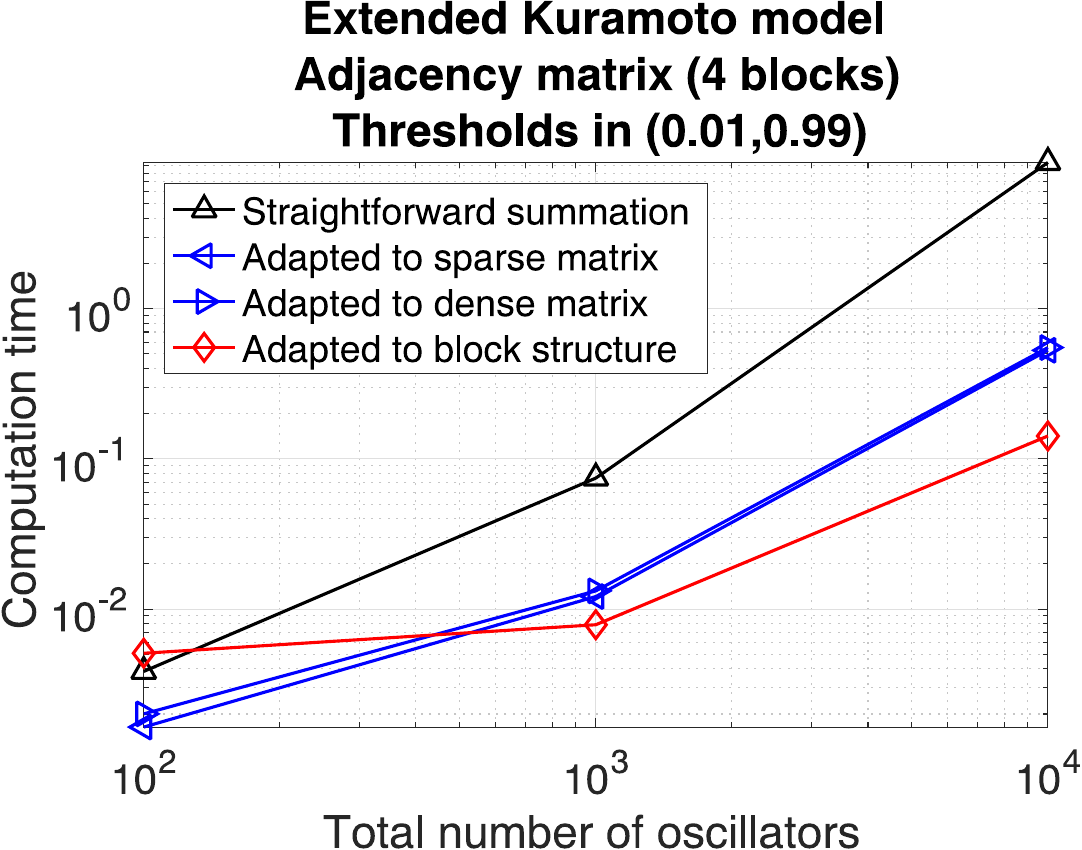} \\[4mm]
\includegraphics[width=3.2cm]{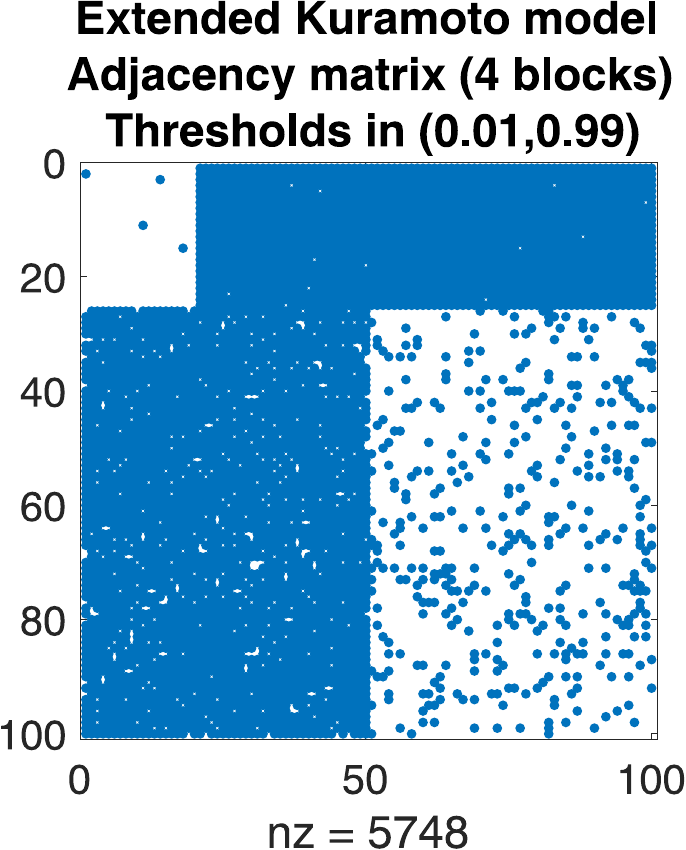} \quad
\includegraphics[width=5cm]{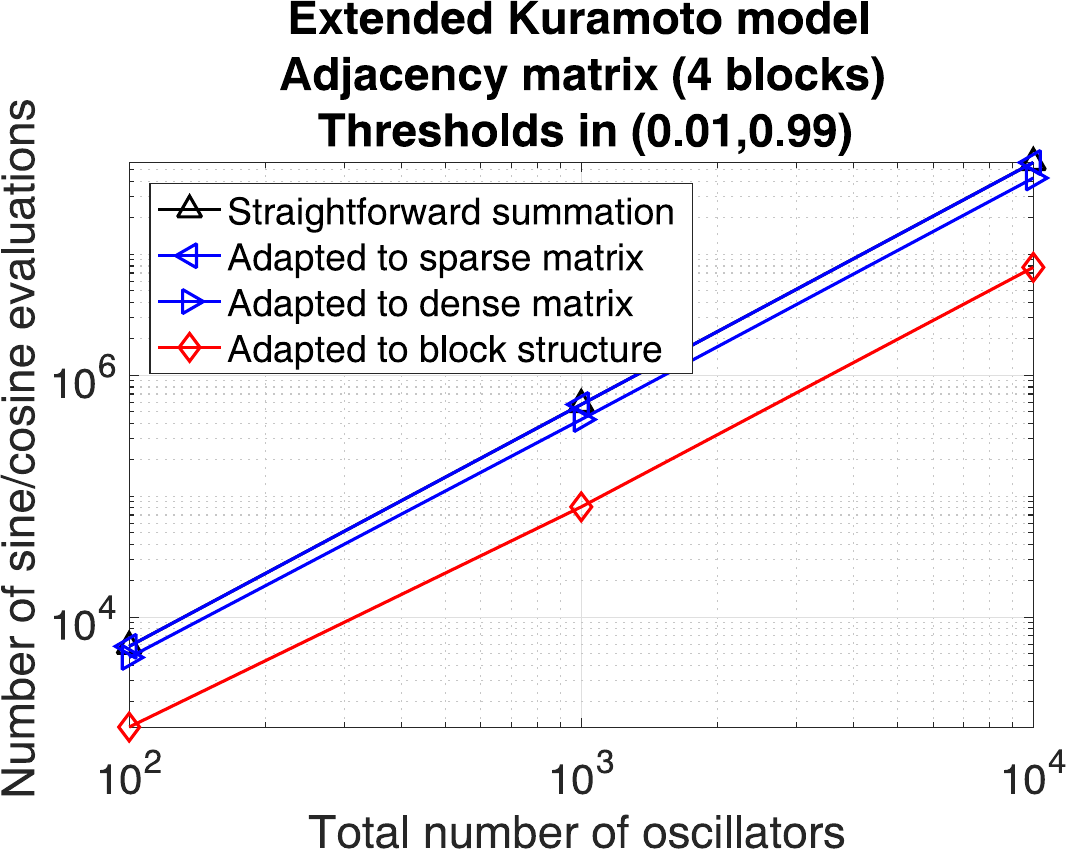} \quad
\includegraphics[width=5cm]{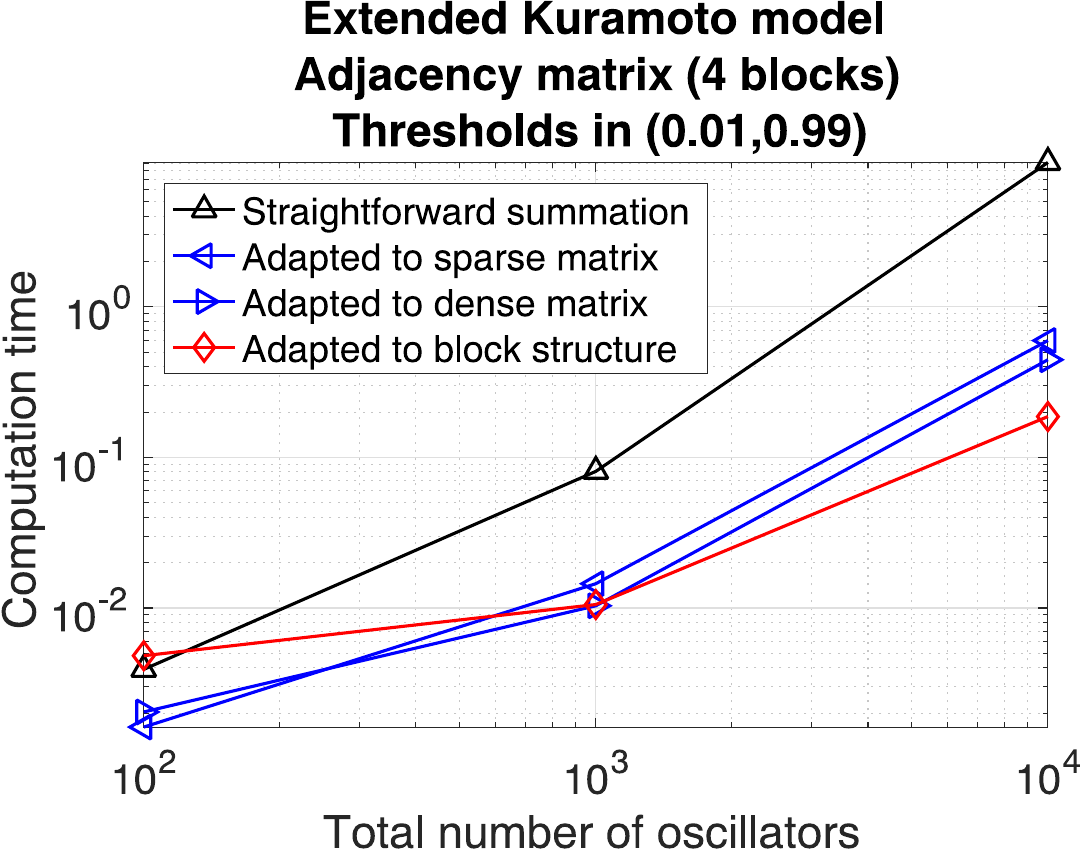}
\end{center}
\caption{Extended Kuramoto models involving randomly generated adjacency matrices with block structure defined through a threshold per block.
Evaluation of the right-hand side by means of different approaches adapted to sparse matrices, dense matrices, and block matrices, respectively.
Left: Illustration of the adjacency matrix for 100 oscillators. 
Middle: Number of sine and cosine evaluations versus the total number of oscillators.
Right: Comparison of the computation time.}
\label{fig:MyFigure10}
\end{figure}

\begin{figure}[t!]
\begin{center}
\includegraphics[width=3.2cm]{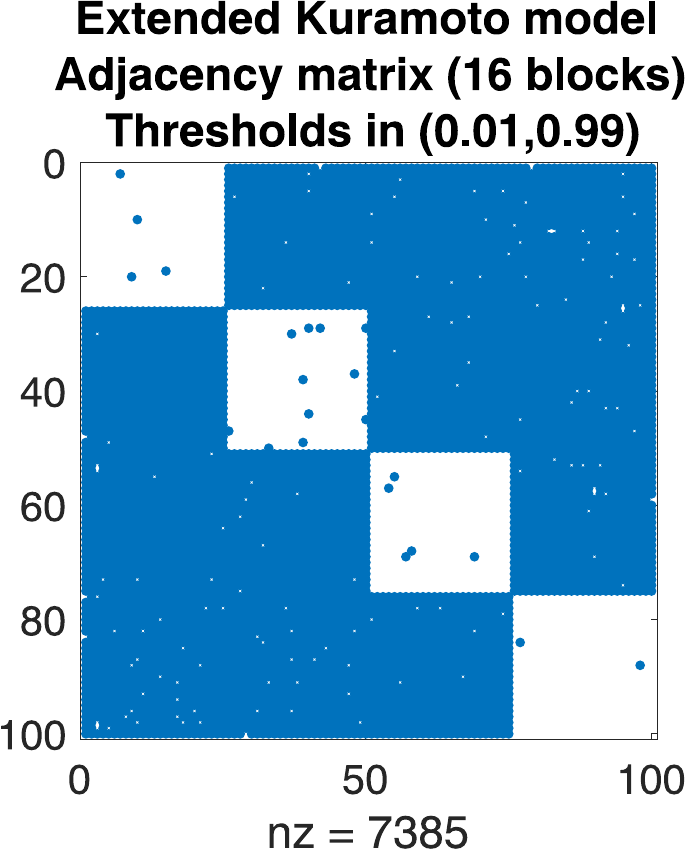} \quad
\includegraphics[width=5cm]{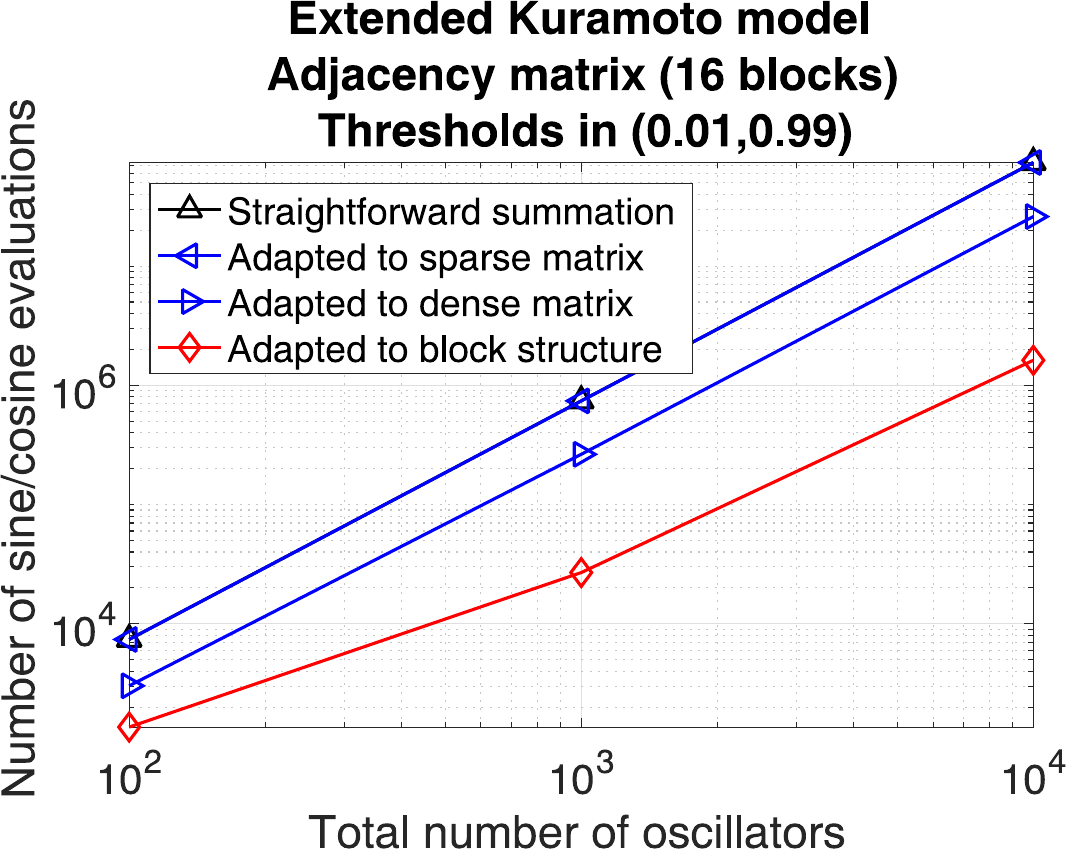} \quad
\includegraphics[width=5cm]{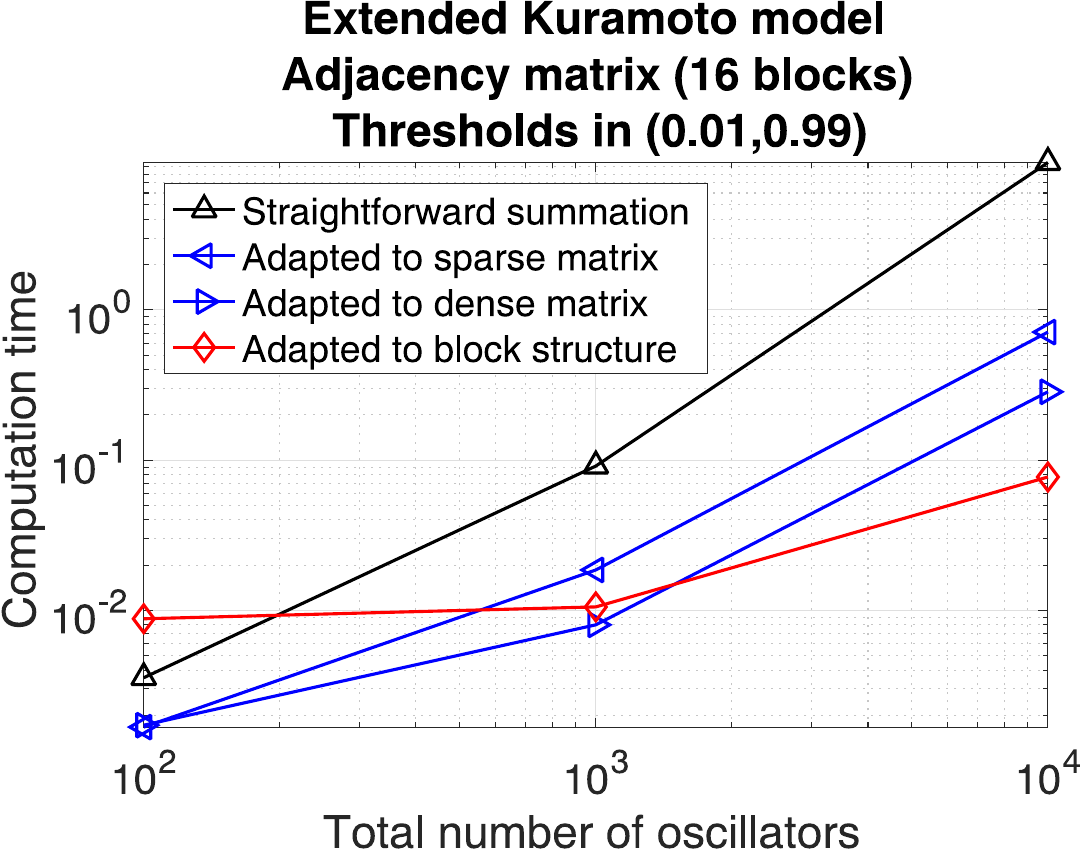} \\[4mm]
\includegraphics[width=3.2cm]{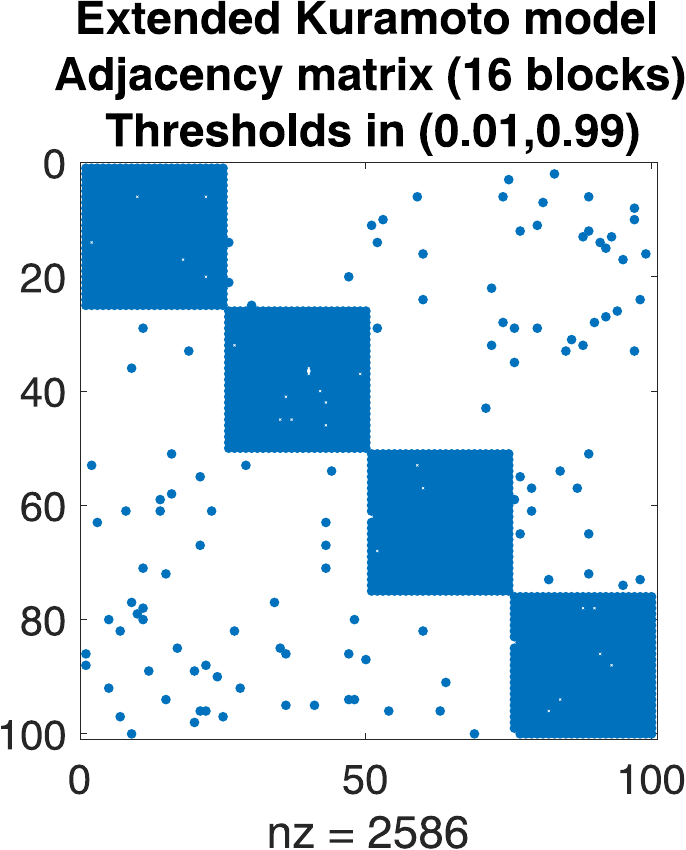} \quad
\includegraphics[width=5cm]{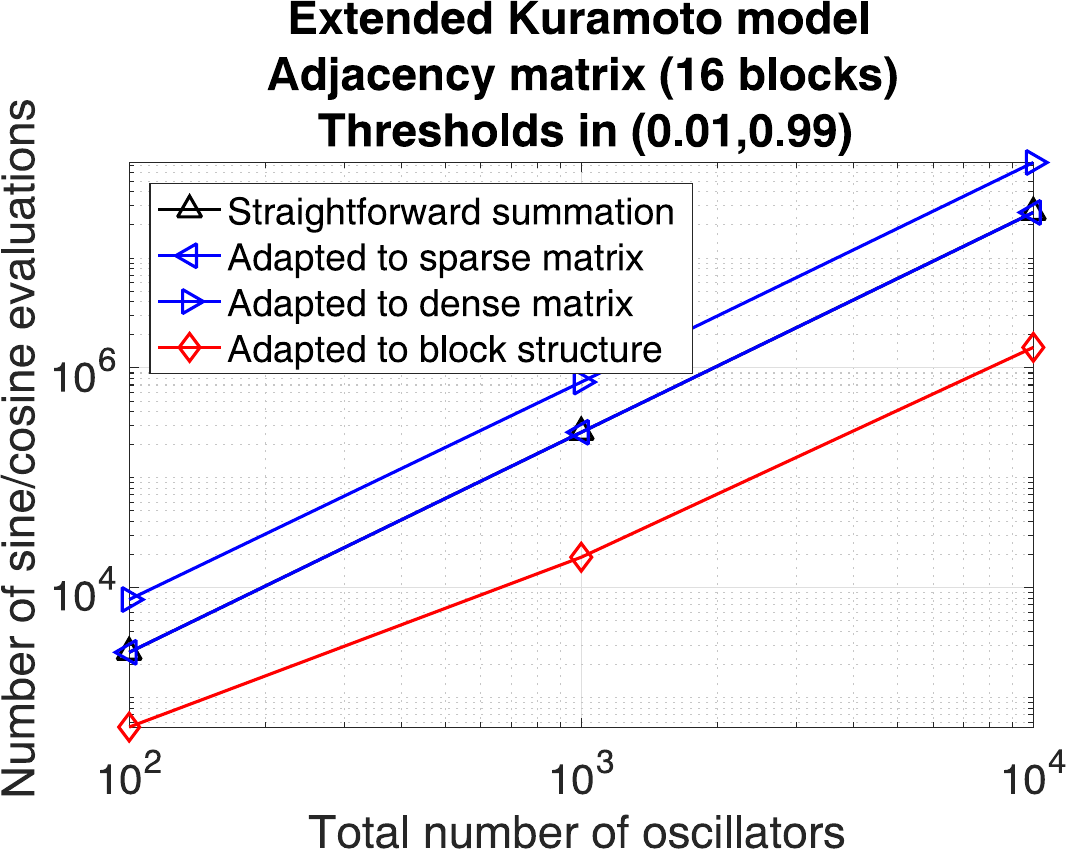} \quad
\includegraphics[width=5cm]{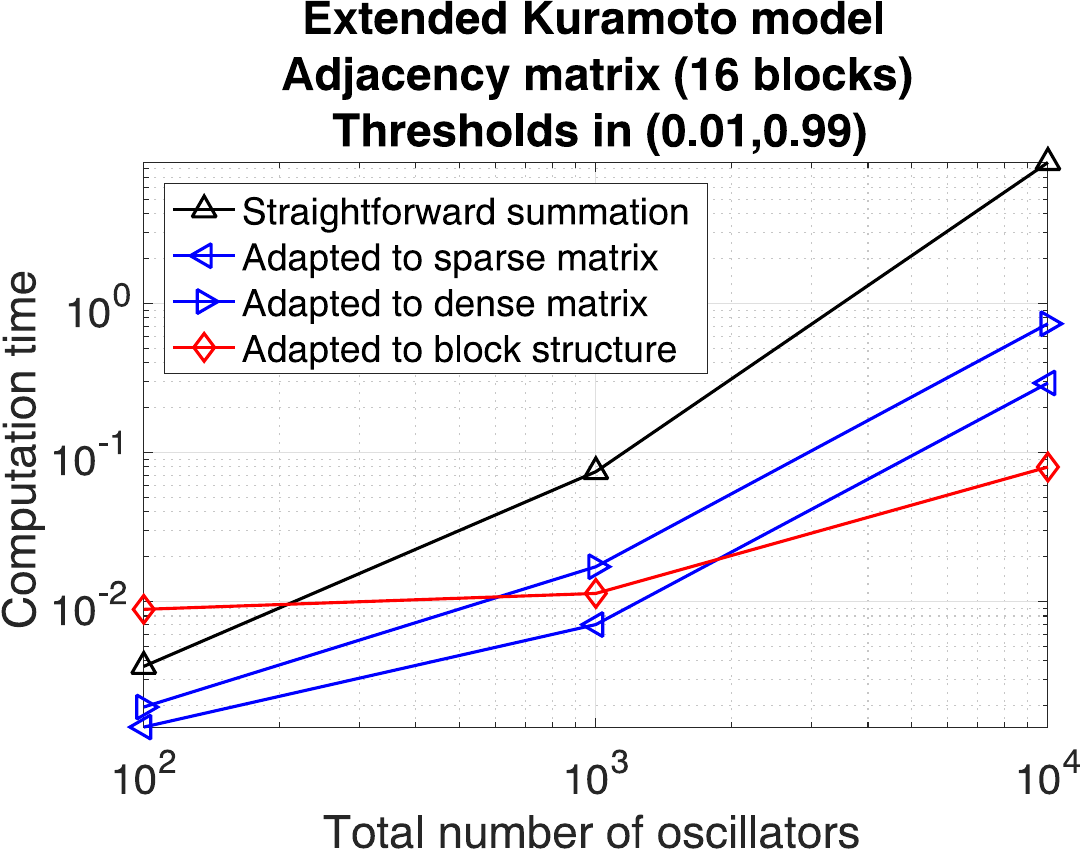} 
\end{center}
\caption{Corresponding results.}
\label{fig:MyFigure11}
\end{figure}

\begin{figure}[t!]
\begin{center}
\includegraphics[width=3.3cm]{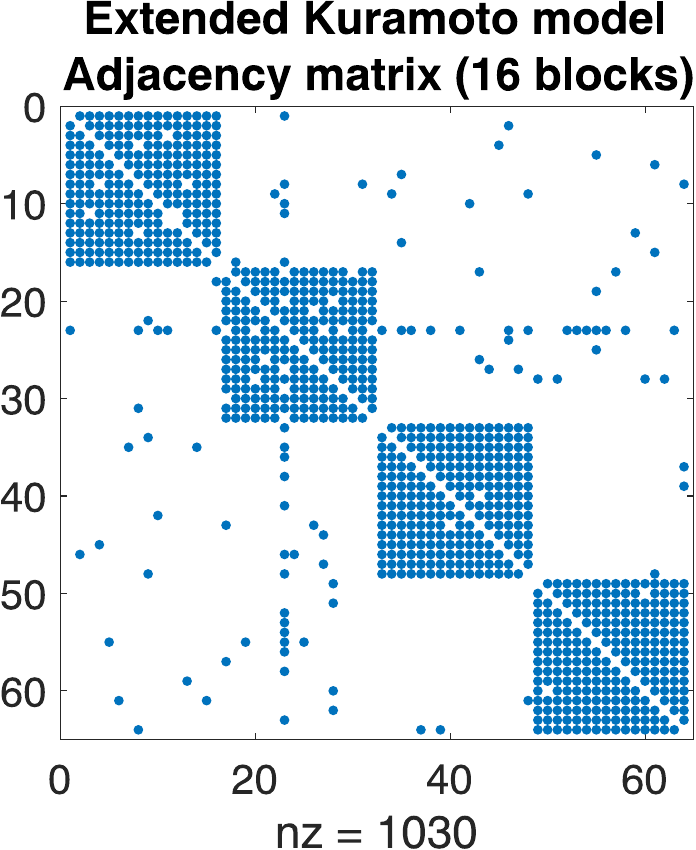} \quad
\includegraphics[width=5cm]{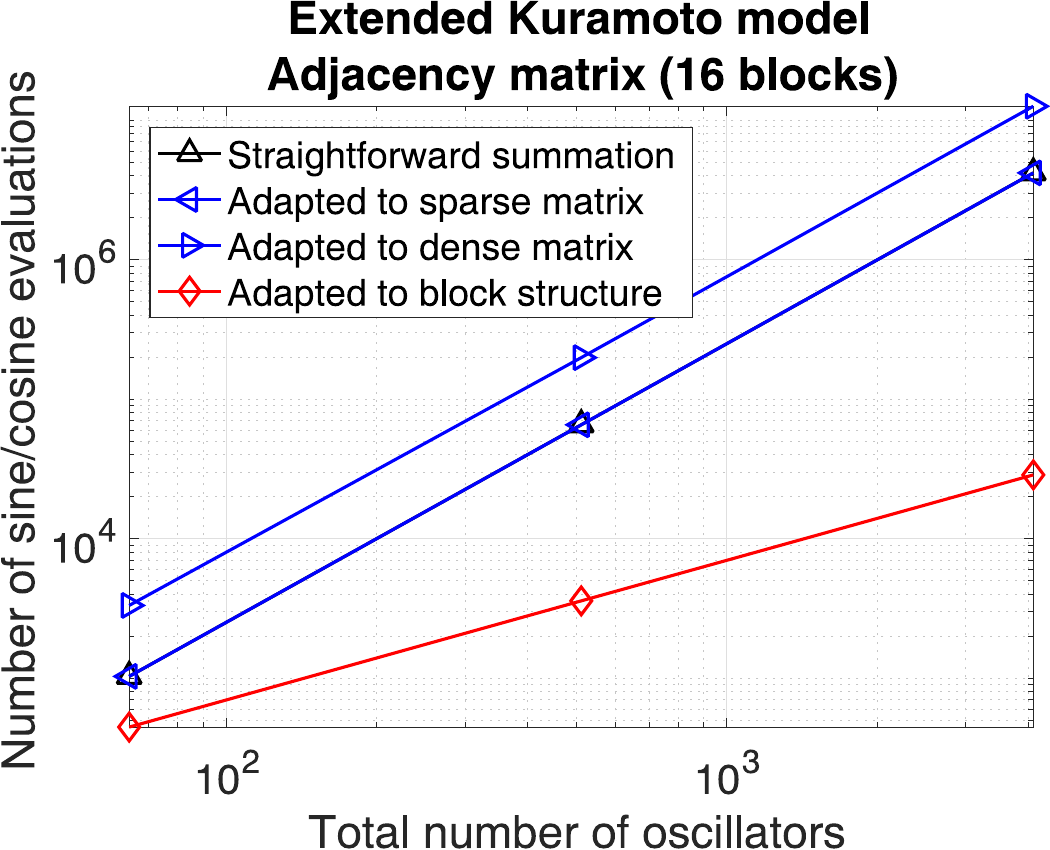} \quad
\includegraphics[width=5cm]{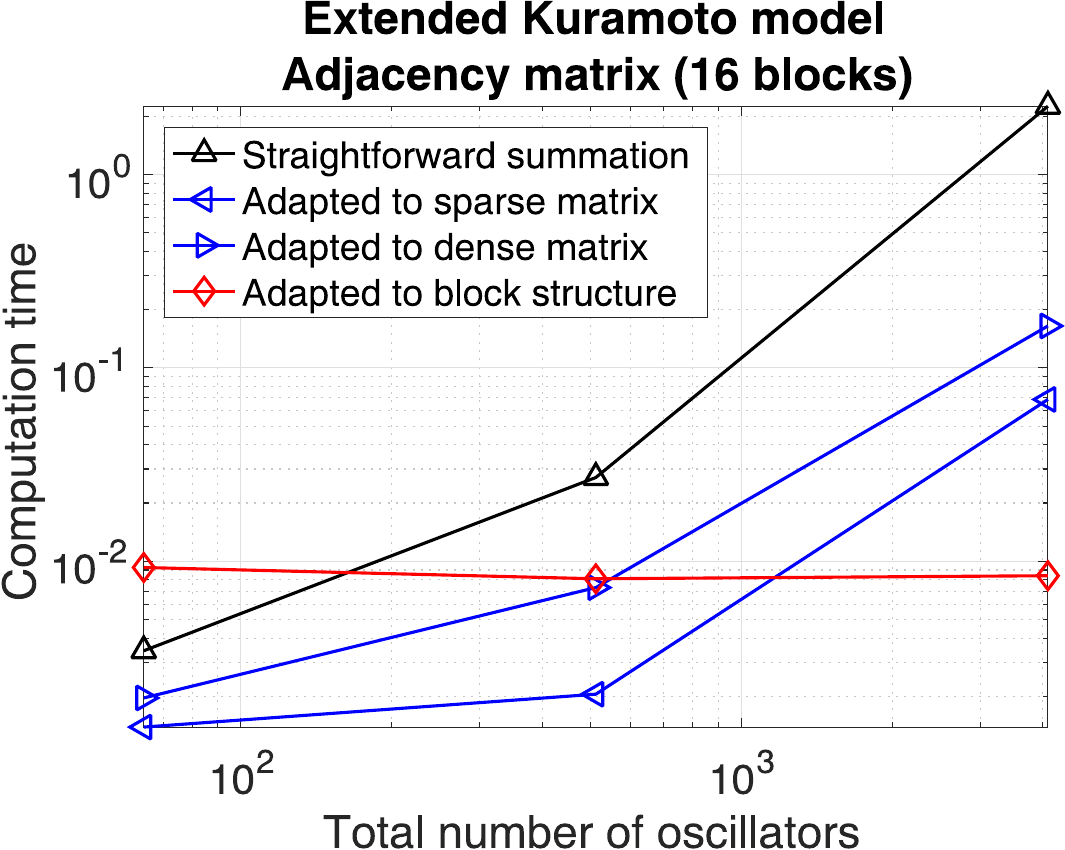} 
\end{center}
\caption{Corresponding results for a more realistic adjacency matrix describing the interactions of four communities of oscillators.}
\label{fig:MyFigure12}
\end{figure}

\MyParagraphAbstand{Numerical comparisons}
In order to test the above described strategy~\eqref{eq:KuramotoGraphApproach} for the efficient computation of the sums
\begin{equation*}
\sum_{\ell=1}^{M} A_{m \ell} \, \sin(\vartheta_{\ell} - \vartheta_m)\,, \quad \vartheta_m = \tfrac{2 \pi m}{M}\,, \quad m \in \{1, 2, \dots, M\}\,,
\end{equation*}
we study two situations, which we consider to be of practical relevance in view of the numerical integration of Kuramoto models on graphs~\eqref{eq:Kuramoto}, see also~\eqref{eq:KuramotoGraph}.
\begin{enumerate}[(i)]
\item
\emph{Matrices without block structure.} \,
On the one hand, we define well-balanced adjacency matrices and compare the number of function evaluations as well as the computation time of an implementation in \textsc{Matlab}.  
For this purpose, we prescribe a threshold $p \in [0,1]$ and generate row-by-row a sequence of uniformly distributed random numbers $z \in [0,1]$. 
Whenever $z > p$, the corresponding coefficient is set to one, otherwise, it is set to zero. 
For $p = 0.99$, e.g., the resulting adjacency matrix is sparse and straightforward summation is advantageous, see Figure~\ref{fig:MyFigure7}.
For the complementary case $p = 0.01$, e.g., the adjacency matrix comprises few non-zero coefficients and the precomputation of sums is favourable, see Figure~\ref{fig:MyFigure9}.
For the threshold $p = 0.5$ both approaches lead to essentially the same counts, see  Figure~\ref{fig:MyFigure8}.
\item
\emph{Matrices with block structure.} \,
On the other hand, we generate adjacency matrices that are composed of well-balanced submatrices.
Each block is connected with a different threshold.
We use 
\begin{enumerate}[(a)]
\item   
a double loop over all rows and straightforward summation of $\sin(\vartheta_{\ell} - \vartheta_m)$ provided that $A_{m \ell} \neq 0$ for $\ell, m \in \{1, 2, \dots, M\}$,
\item
a single loop to sum over all indices that correspond to non-zero coefficients, realised by a script of the form $\text{sum}(\sin(\vartheta(\text{NonZero})-\vartheta_m))$ for $m \in \{1, 2, \dots, M\}$,
\item
the precomputation of sums and a single loop to subtract terms that correspond to coefficients equal to zero, applying the script 
$\text{sum}(\sin(\vartheta(\text{Zero})-\vartheta_m))$ for $m \in \{1, 2, \dots, M\}$,
\item 
the procedure in (b) and (c) adapted to each submatrix. 
\end{enumerate}
The obtained results, displayed in Figures~\ref{fig:MyFigure10}--\ref{fig:MyFigure12}, confirm that the latter approach is beneficial for a higher number of oscillators, where the evaluation of functions and the computation of sums are expected to be the most time consuming components. 
\end{enumerate}

\begin{figure}[t!]
\begin{center}
\includegraphics[width=4cm]{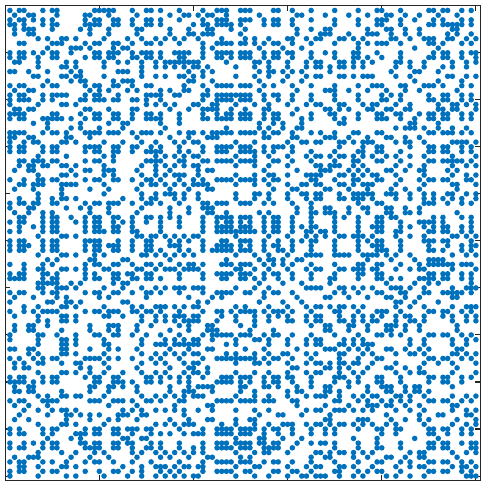} \quad
\includegraphics[width=4cm]{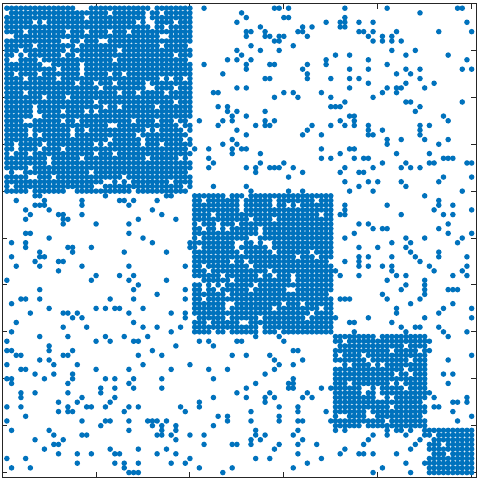}
\end{center}
\caption{Adjacency matrices~$A$ and~$P A \, P^T$ with suitably chosen permutation matrix~$P$.
Left: The structure of the underlying graph is not evident.   
Right: A separation into communities of oscillators is recognisable.
Matrices of this form are used for numerical tests of community detection algorithms.}
\label{fig:MyFigure13}
\end{figure}

\begin{figure}[t!]
\begin{center}
\includegraphics[width=4cm]{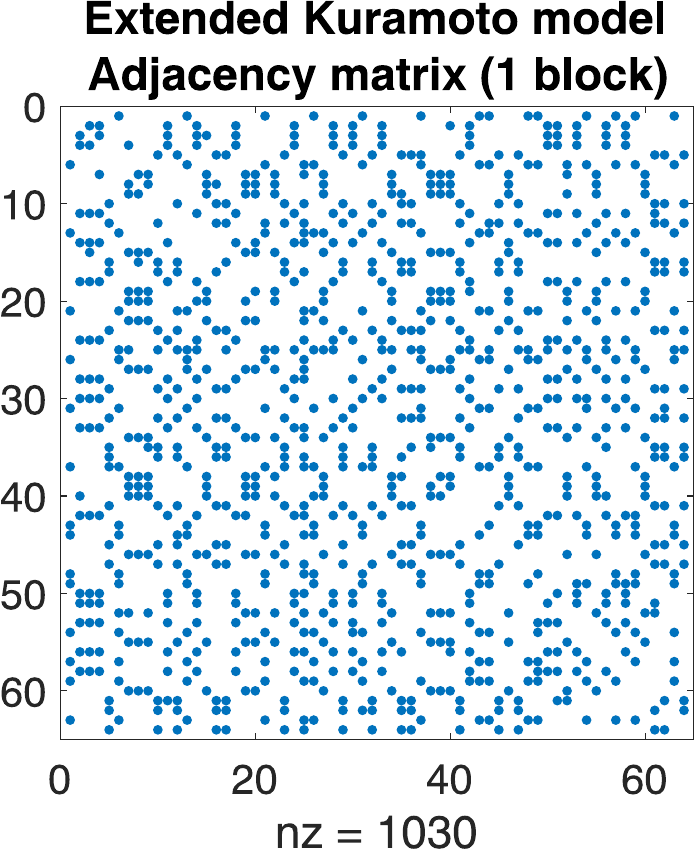} \quad
\includegraphics[width=4cm]{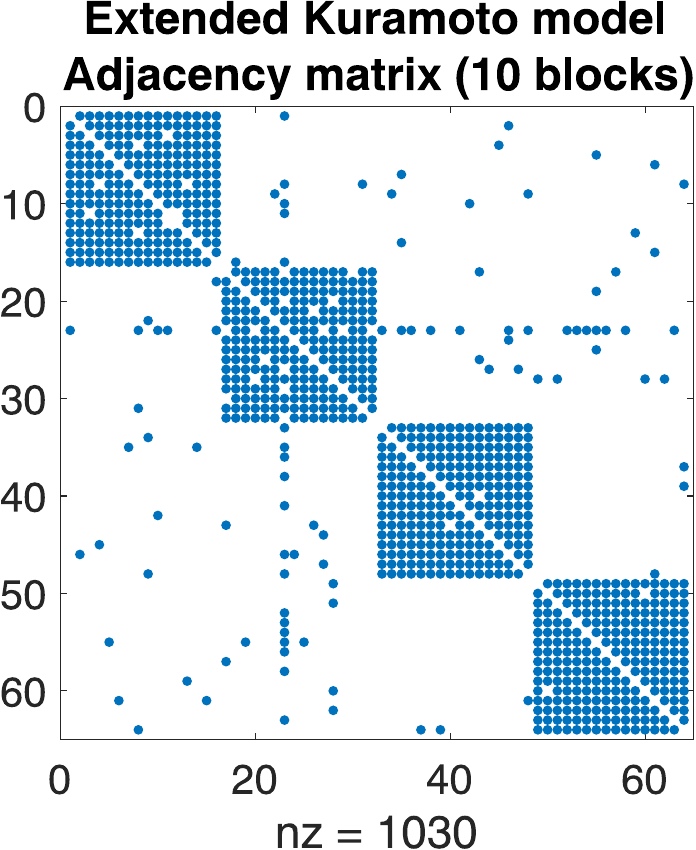} 
\end{center}
\caption{Adjacency matrices~$A$ and~$P A \, P^T$ with suitably chosen permutation matrix~$P$.
Left: The structure of the underlying graph is not evident.   
Right: A separation into communities of oscillators is recognisable.
Symmetric adjacency matrices of this form are considered in connection with the numerical integration of Kuramoto models on graphs.
See also Figure~\ref{fig:MyFigure12}.}
\label{fig:MyFigure14}
\end{figure}

\begin{figure}[t!]
\begin{center}
\includegraphics[width=5.6cm]{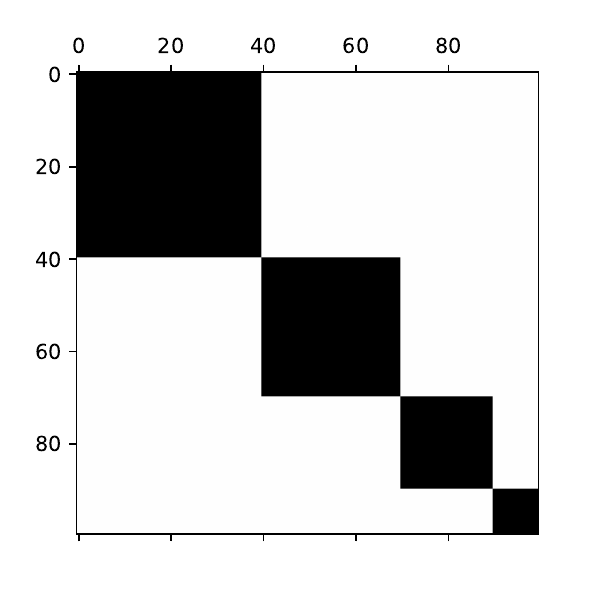} 
\includegraphics[width=8.5cm]{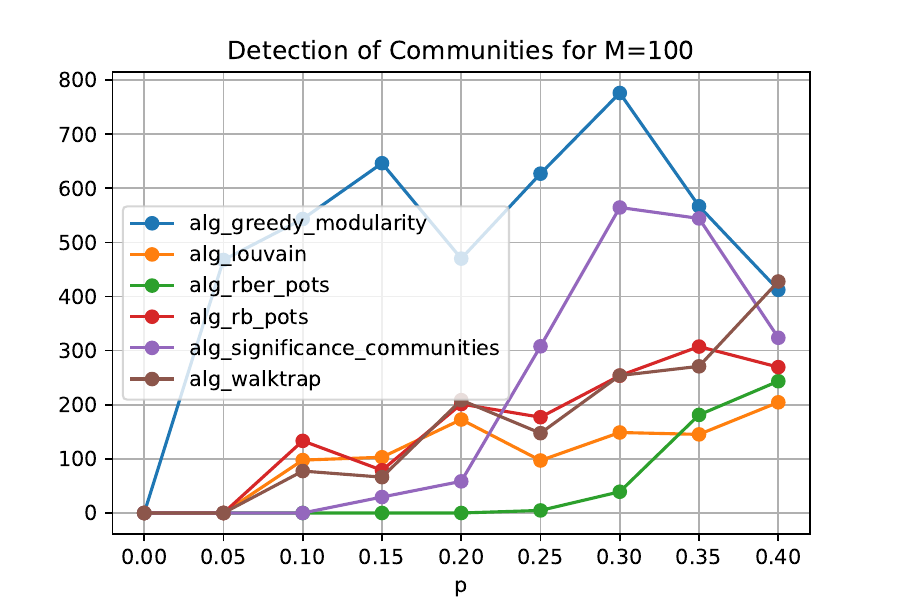} \\
\includegraphics[width=5.6cm]{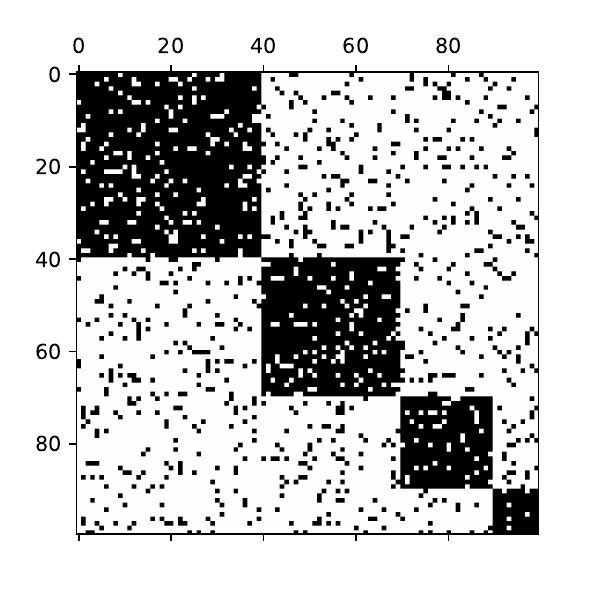} 
\includegraphics[width=8.5cm]{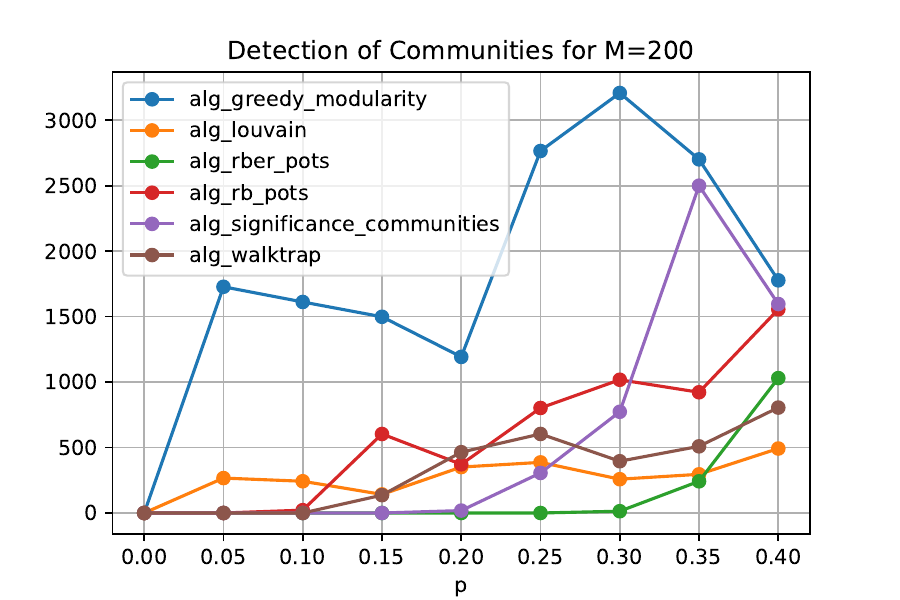} \\
\includegraphics[width=5.6cm]{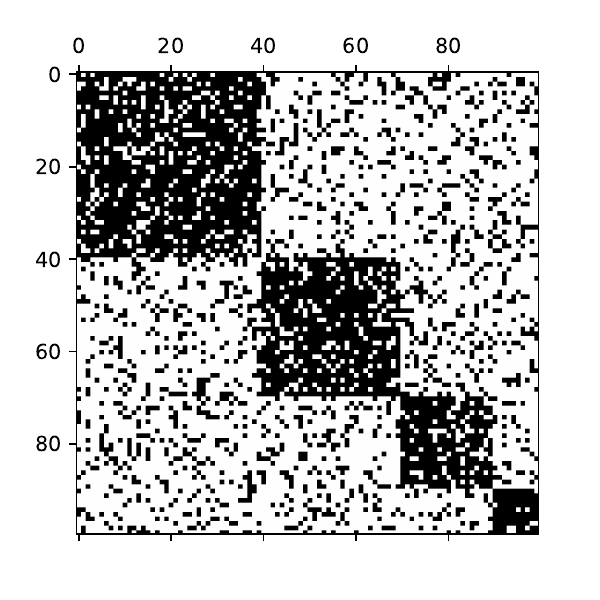} 
\includegraphics[width=8.5cm]{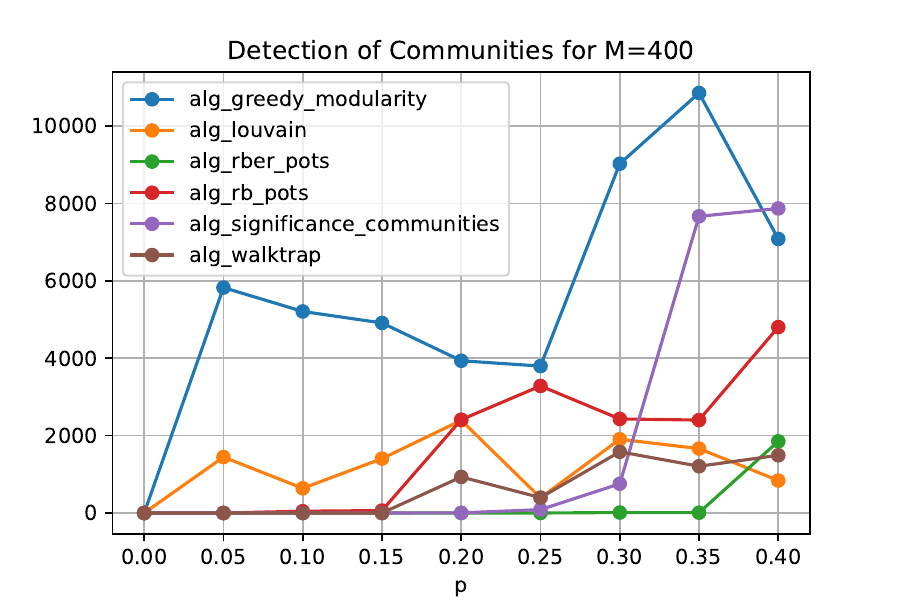} 
\end{center}
\caption{\colred{Left: Matrix comprising four blocks ($p = 0$, top) and related adjacency matrices for increasing thresholds $p \in \{0.1, 0.2\}$ (middle to bottom).
Right: Average of the quantity~\eqref{eq:QuantityPerformance} over eight runs, which is chosen as a measure for the performance of the considered community detection algorithms, see also Table~\ref{tab:Table2}. Obtained results for $M = 100$ (top), $M = 200$ (middle), $M = 400$ (bottom).}}
\label{fig:MyFigure15}
\end{figure}

\begin{figure}[t!]
\begin{center}
\includegraphics[width=5.6cm]{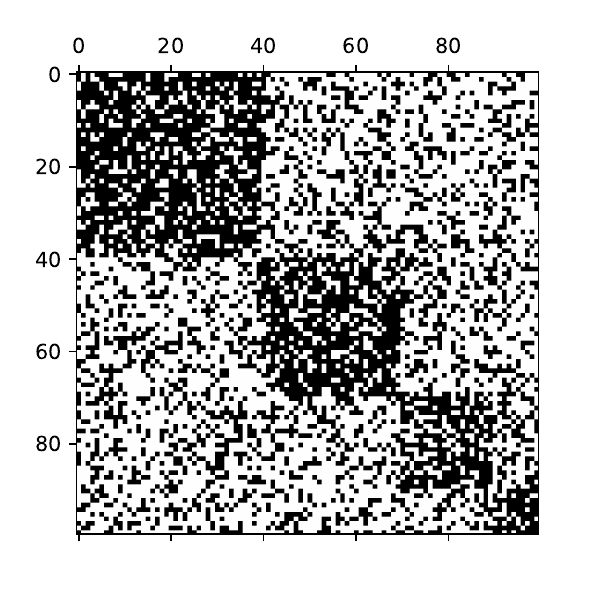} 
\includegraphics[width=8.5cm]{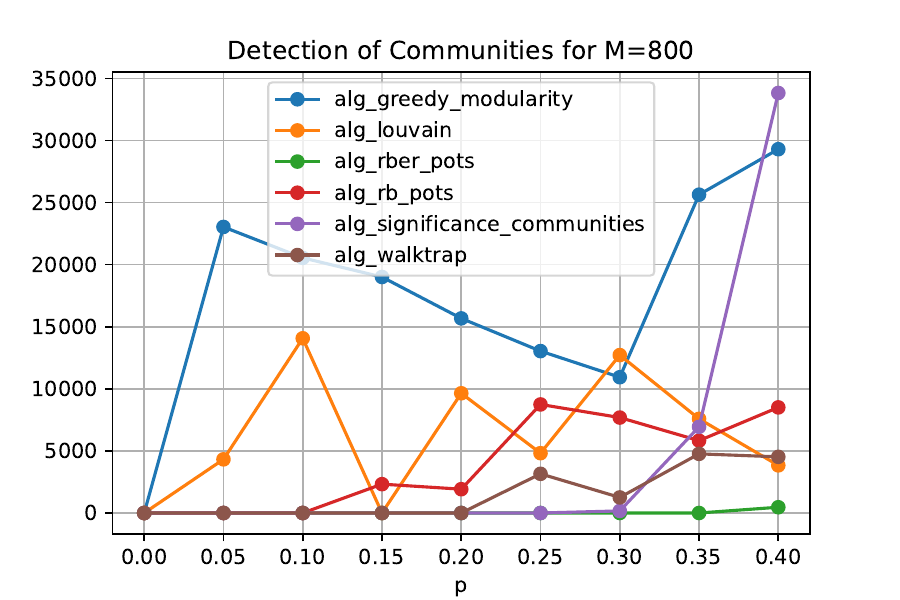} \\
\includegraphics[width=5.6cm]{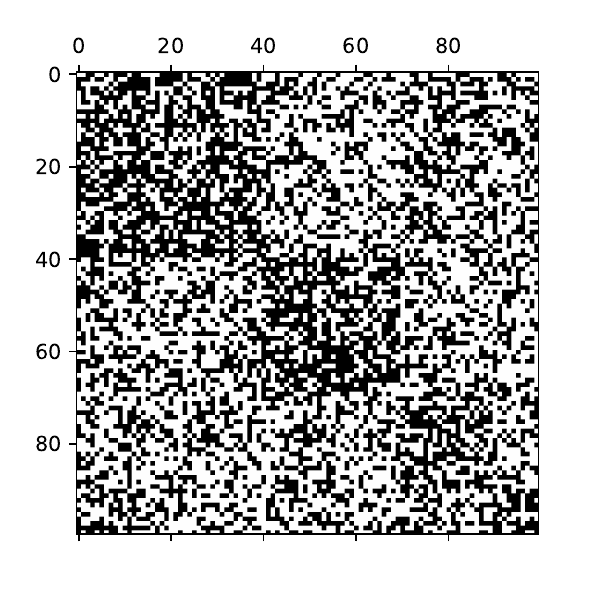} 
\includegraphics[width=8.5cm]{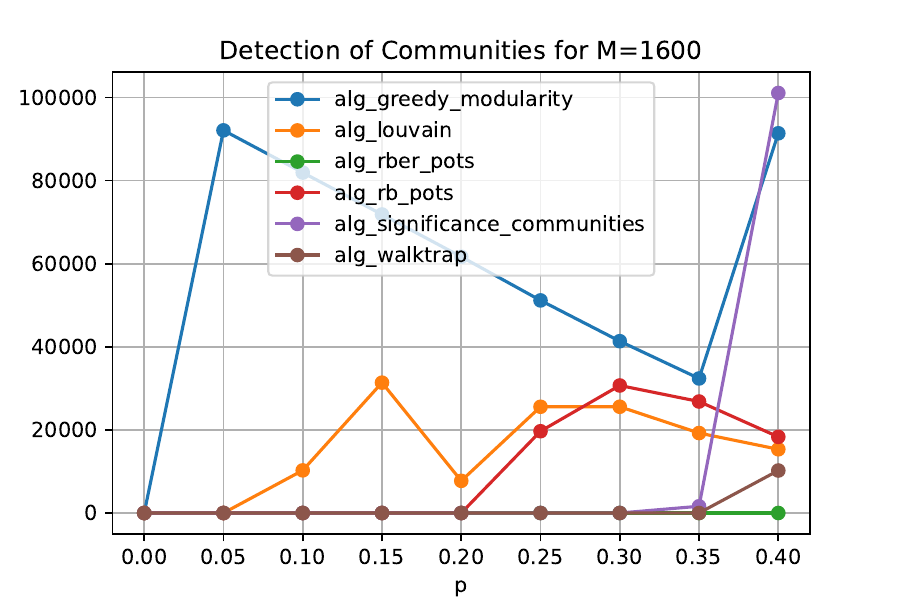} 
\end{center}
\caption{\colred{Left: Adjacency matrices for increasing thresholds $p \in \{0.3, 0.4\}$ (top to bottom). 
Right: Corresponding results for $M = 800$ (top), $M = 1600$ (bottom).}}
\label{fig:MyFigure16}
\end{figure}

\begin{figure}[t!]
\begin{center}
\includegraphics[width=7.1cm]{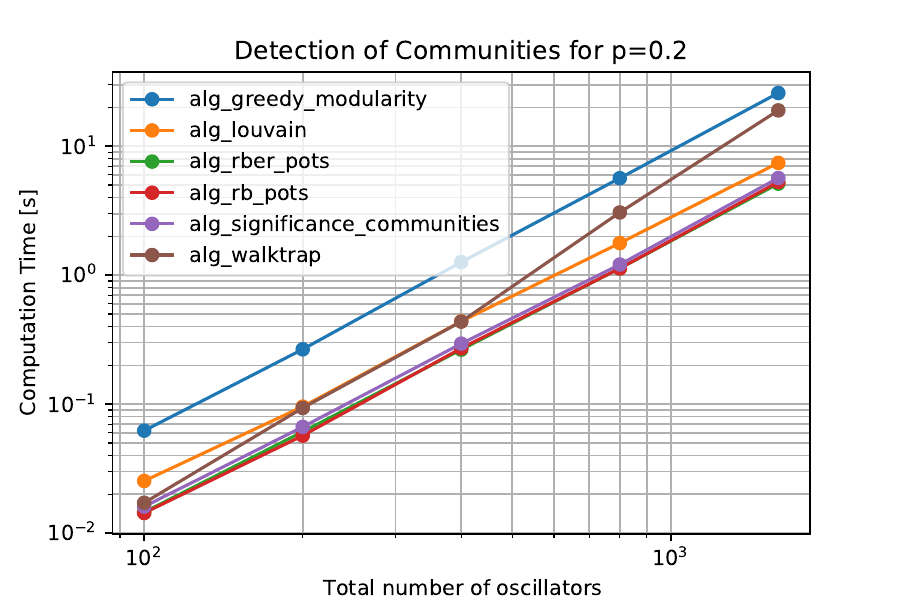} 
\includegraphics[width=7.1cm]{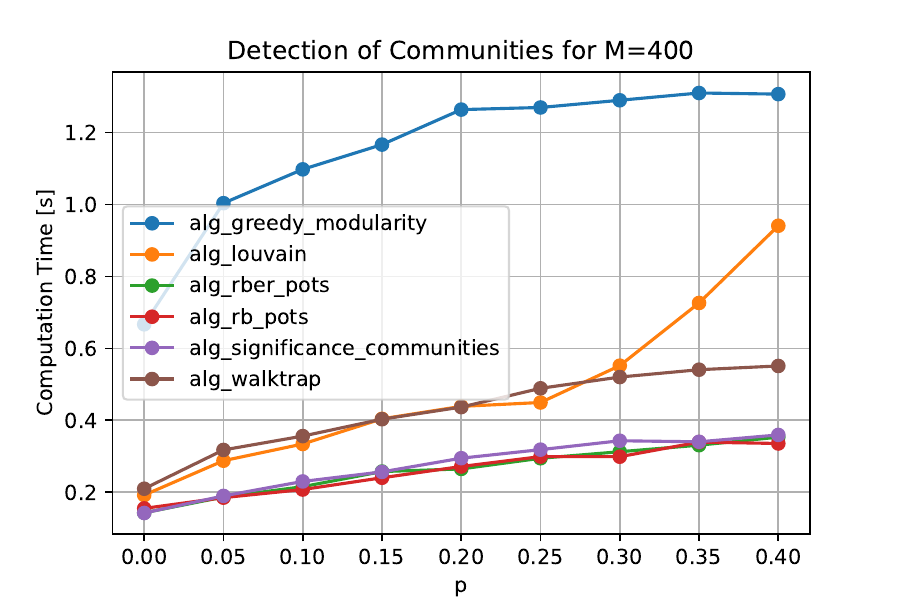} 
\end{center}
\caption{\colred{Average computation time of the considered community detection algorithms over eight runs. 
Left: The displayed numerical results, obtained for the threshold $p = 0.2$, reflect a quadratic scaling with respect to the total number of oscillators.
Right: For larger deviations of the adjacency matrices from the underlying block diagonal matrix, related to increasing values of the threshold~$p$, the computation time increases, in general.}}
\label{fig:MyFigure17}
\end{figure}

\MyParagraphAbstand{Community detection in graphs}
Our numerical comparisons show that it is to the best advantage to take the underlying structure of the considered Kuramoto model on a graph~\eqref{eq:KuramotoGraph} into account.
A suitable reordering of the governing equations accordingly to the separation of the oscillators into communities corresponds to the transformation of the associated adjacency matrix to a well-balanced block matrix and is a fundamental means in view of their efficient numerical integration, see Figures~\ref{fig:MyFigure13}--\ref{fig:MyFigure14}. 
In the following, we study the performance of various algorithms for the detection of communities in graphs. 
We restrict ourselves to algorithms from the \textsc{Python} packages \textsc{networkx} and \textsc{cdlib}, which were suitable for our purposes, did not require additional parameters, and took a reasonable computation time, see Table~\ref{tab:Table2}.
We focus on the relevant case, where interactions primarily take place within certain communities of oscillators.
\colred{As illustrated in Figure~\ref{fig:MyFigure15}, our starting point is a matrix (of suitably chosen dimension) with fully occupied blocks along the diagonal
\begin{equation*}
\begin{gathered}  
B = (B_{m \ell})_{m, \ell \in \{1, 2, \dots, M\}}
= \begin{pmatrix} B^{(11)} &&& 0 \\ & B^{(22)} && \\ && B^{(33)} & \\ 0 &&& B^{(44)} \end{pmatrix} \in \RR^{M \times M}\,, \\
B^{(kk)} = (1)_{i,j \in \{1, 2, \dots, (5-k) M/10\}} \in \RR^{(5-k) M/10 \times (5-k) M/10}\,, \quad k \in \{1, 2, 3, 4\}\,.
\end{gathered}  
\end{equation*}}%
Similarly to before, we prescribe a threshold $p \in [0,1]$ that determines the probability that a coefficient of the matrix~$B$ is changed from one to zero or from zero to one, respectively.
That is, we generate uniformly distributed random numbers $z_{m \ell} \in [0,1]$ for $m, \ell \in \{1, 2, \dots, M\}$ and define
\begin{equation*}
\begin{gathered}  
A = (A_{m \ell})_{m, \ell \in \{1, 2, \dots, M\}} \in \RR^{M \times M}\,, \\
A_{m \ell} = \begin{cases} 1 - B_{m \ell}\,, &z_{m \ell} < p\,, \\ B_{m \ell}\,, &z_{m \ell} \geq p\,, \end{cases} \quad m, \ell \in \{1, 2, \dots, M\}\,. 
\end{gathered}  
\end{equation*}
The higher the threshold, the stronger the deviation of the adjacency matrix from the related block matrix, see Figures~\ref{fig:MyFigure15}--\ref{fig:MyFigure16}.
The community detection algorithms were applied to the randomly permuted adjacency matrix and yield as outputs partitions of the nodes into communities.
We associate them with permutations of the sequence $(1, 2, \dots, M)$ such that nodes in one community are arranged one after the other.
These permutations correspond to transformations of the original matrices that can be cast into the form 
\begin{equation*}
\widetilde{A} = P A \, P^T
\end{equation*}
with a permutation matrix $P \in \RR^{M \times M}$ such that in particular the identity $P^{-1} = P^T$ is valid.
For each of these matrices, we identify a matrix with fully occupied blocks along the diagonal reflecting the detected communities
\begin{equation*}
\widetilde{B} = \begin{pmatrix} \widetilde{B}^{(11)} && 0 \\ & \widetilde{B}^{(22)} & \\ 0 && \ddots \end{pmatrix} \in \RR^{M \times M}\,.
\end{equation*}
We recall that the computation of pointwise products such as 
\begin{equation*}
\big(B \sin(\vartheta)\big) \CompProduct \cos(\vartheta)\,, \quad
\big(\widetilde{B} \sin(\vartheta)\big) \CompProduct \cos(\vartheta)\,, \quad \vartheta \in \RR^M\,, 
\end{equation*}
based on the precomputation of sums requires in total~$2 M$ evaluations of cosine and sine functions and that the number of non-zero coefficients of the matrices 
\begin{equation*}
C = A - B\,, \quad \widetilde{C} = \widetilde{A} - \widetilde{B}\,, 
\end{equation*}
match the additional computational costs. 
For this reason, we determine the quantity 
\begin{equation}
\label{eq:QuantityPerformance}  
\sum_{\ell,m=1}^{M} \absbig{\widetilde{C}_{m \ell}} - \sum_{\ell,m=1}^{M} \absbig{C_{m \ell}}\,,  
\end{equation}
\colred{in order to assess the performance of the algorithms, see \cite{PorterOnnelaMucha} for detailed explanations.}
The numerical results obtained for $M \in \{100, 200, 400, 800, 1600\}$ and thresholds in the range~$[0, 0.4]$ are displayed in  Figures~\ref{fig:MyFigure15}--\ref{fig:MyFigure17}.
In case $p = 0$, i.e.~for a random permutation of the fully occupied block diagonal matrix, each of the tested algorithms detected the four communities. 
For larger deviations of the adjacency matrices from the underlying block diagonal matrix, however, some algorithms failed. 
Overall, the community detection algorithm \textsc{rber\_pots}, available through the \textsc{Python} package \textsc{cdlib}, provided the most reliable results for a higher number of oscillators. 

\MyParagraphAbstand{Favourable community detection algorithm}
It is notable that there exist several variants of the algorithm \textsc{rber\_pots}~\cite{Reichardt2004,Reichardt2006} with the common objective to detect communities in a graph with associated adjacency matrix~$A$.
The standard implementation in \textsc{cdlib}~\cite{CDlib} minimises the quantity 
\begin{align*}
Q = - \sum_{\ell,m=1}^M (A_{m \ell}- p) \, \delta(\sigma_m, \sigma_{\ell})\,, 
\end{align*}
where $p \in (0, 1)$ represents the mean density of edges in a graph, that is, the ratio between the numbers of actually existing and potential edges.
Whenever two nodes $m, \ell \in \{1, 2, \dots, M\}$ belong to the same community, the quantity $\delta(\sigma_m, \sigma_{\ell})$ takes the value one,
and it is zero otherwise.
Provided that the mean edge density within a community is higher than the mean edge density of the complete network, a node is included in a community.
The favourable performance of \textsc{rber\_pots} in the considered example is explained by the fact that the four communities are well recognisable by the mean edge densities, even for $p = 0.4$.

\begin{table}[t]
\tbl{Algorithms from \textsc{Python} packages applied for community detection in graphs. \\}
{\qquad\qquad\qquad
\begin{tabular}{c|c} 
\textsc{networkx}~\cite{NetworkX} & \textsc{greedy\_modularity}~\cite{Clauset2004} \\\midrule
\textsc{cdlib}~\cite{CDlib} & \textsc{louvain}~\cite{Blondel2008} \\
& \textsc{rber\_pots}~\cite{Reichardt2004, Reichardt2006} \\
& \textsc{rb\_pots}~\cite{Leicht2008, Reichardt2006} \\
& \textsc{significance\_communities}~\cite{Traag2013} \\
& \textsc{walktrap}~\cite{Pons2005} 
\end{tabular}\qquad\qquad\qquad}
\label{tab:Table2}
\end{table}

\begin{figure}[t!]
\begin{center}
\includegraphics[width=7.1cm]{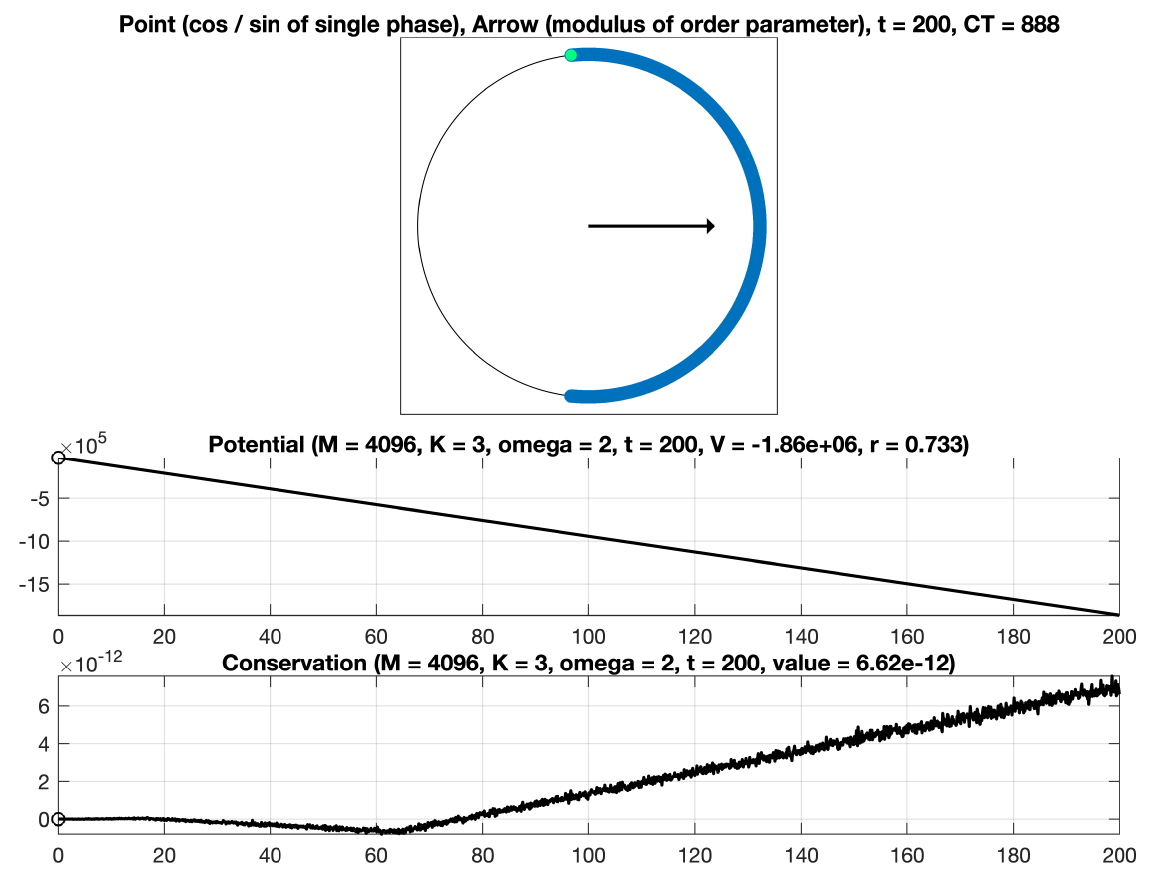} 
\includegraphics[width=7.1cm]{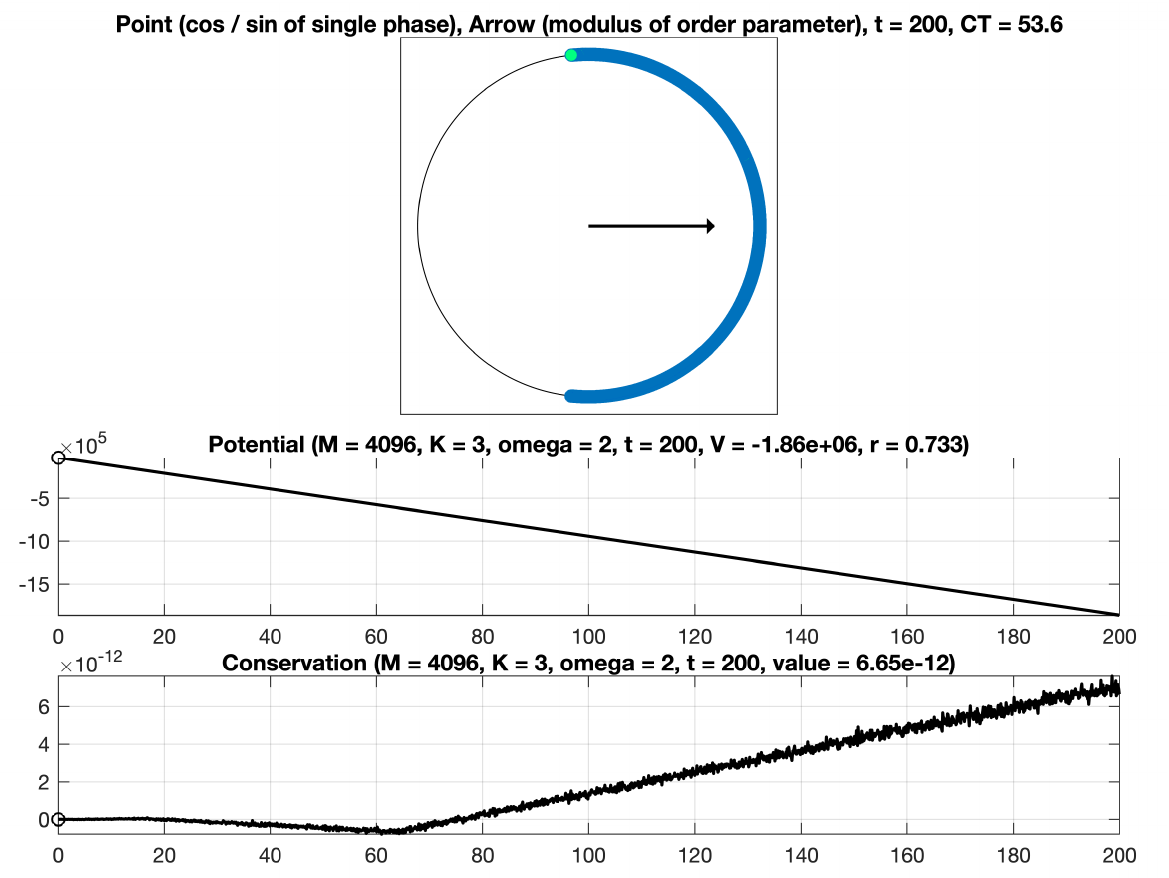}
\end{center}
\caption{\colred{Numerical integration of a Kuramoto model on a graph with the common uniform scaling. 
The associated adjacency matrix has a structure as illustrated in Figure~\ref{fig:MyFigure14} (left: matrix without recognisable block structure, right: permuted matrix with recognisable block structure). 
Left: Evaluation of the right-hand side by straightforward summation.
Right: Employing the underlying block structure and using the precomputation of sums permits a significant reduction of the computation time CT.}}
\label{fig:MyFigure18}
\end{figure}

\begin{figure}[t!]
\begin{center}
\includegraphics[width=7.1cm]{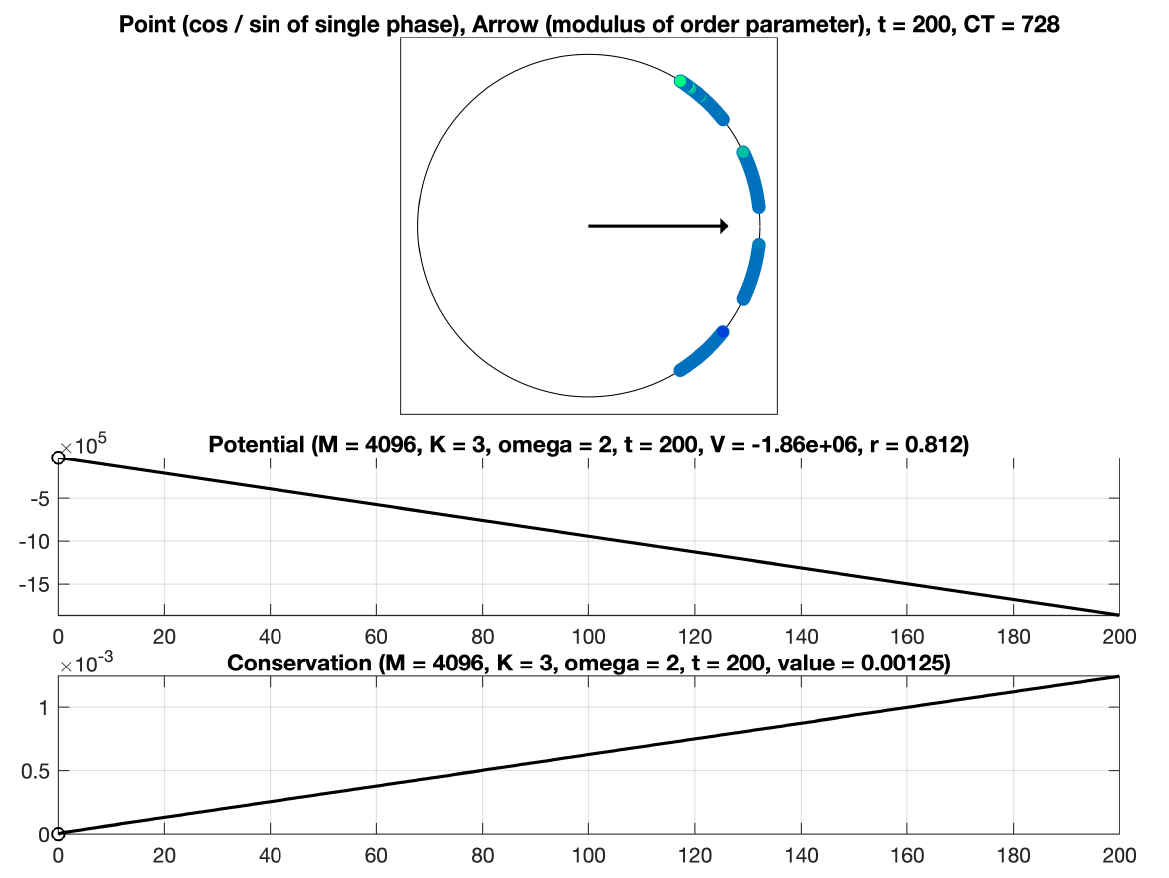} 
\includegraphics[width=7.1cm]{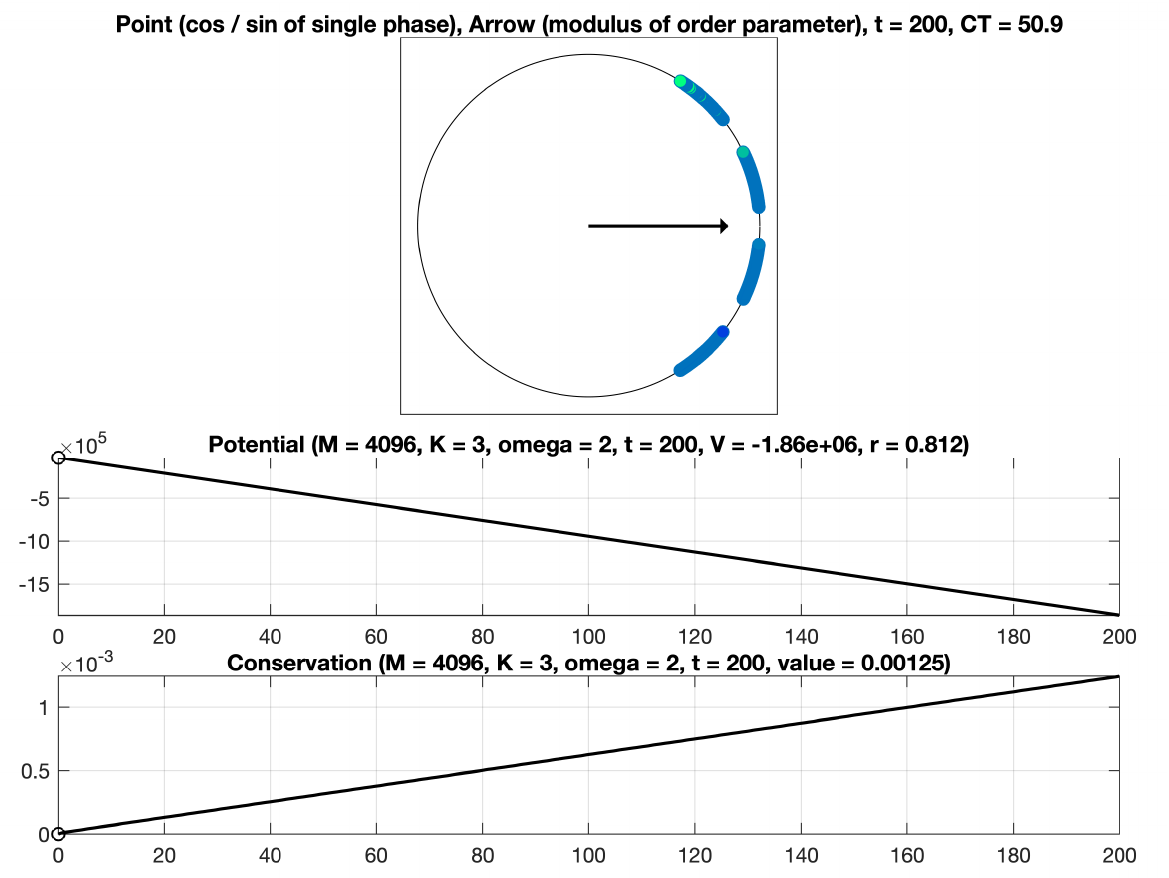} 
\end{center}
\caption{\colred{Corresponding results for the non-uniform scaling. Synchronisation within four communities is observed. Due to the lack of symmetry of the system, the conservation property does not hold.}}
\label{fig:MyFigure19}
\end{figure}

\MyParagraphAbstand{Numerical integration}
\colred{
Following up the numerical integration of the classical Kuramoto model~\eqref{eq:Kuramoto}--\eqref{eq:KuramotoReformulation}, we finally study extended Kuramoto models on graphs~\eqref{eq:KuramotoGraph}.
On the one hand, we consider the situation, where a separation of the oscillators into four communities of the same size is evident.
The structure of the associated symmetric adjacency matrix is illustrated in Figure~\ref{fig:MyFigure14} (right).
Here, our approach based on the precomputation of sums applies.
The numerical results obtained for the common uniform scaling, a total number of $M = 8^4 = 4096$ oscillators, coupling constant $K = 3$, final time $T = 200$, and intrinsic frequencies as well as initial phases of the form 
\begin{equation*}
\omega_0 = 2\,, \quad \omega_m = 1 + \omega_0 \, \tfrac{(2 m - M - 1)}{M - 1}\,, \quad 
\vartheta_m(0) = \tfrac{2 \pi m}{M}\,, \quad m \in \{1, 2, \dots, M\}\,, 
\end{equation*}
are displayed in Figure~\ref{fig:MyFigure18} (right).
On the other hand, we consider the equivalent system without recognisable block structure, see Figure~\ref{fig:MyFigure14} (left), using straightforward summation for the evaluation of the right-hand side. 
In order to achieve consistency with the previous case, we permute internal frequencies and initial phases accordingly, perform the time integration, and reorder the solution values subsequently.
By comparison of the numerical results shown in Figure~\ref{fig:MyFigure18}, it is evident that the efficient evaluation of the decisive sums based on the block structure of the adjacency matrix is beneficial and permits a significant reduction of the computation time from approximately 890 seconds to about 54 seconds.
For the sake of comparability with the classical Kuramoto, we display the values of the corresponding potentials, which decrease when time evolves, and verify the conservation property, see also~\eqref{eq:ExtendedPotentialConservation}.
The analogous results for the non-uniform scaling are given in Figure~\ref{fig:MyFigure19}. 
As expected, due to the lack of symmetry, the constraint~\eqref{eq:KuramotoGraphConservation} is not fulfilled.
}
\section{Conclusions and outlook}
\colred{
In summary, we have presented results which confirm that the combination of detecting communities and using localisations of order parameters permits significant improvements regarding the reliable and efficient numerical integration of dynamical systems on graphs.
Our approach provides a natural way to exploit the underlying graph structure and local mean-field variables.

Popular models, where the employed concepts apply at once, are given by consensus problems~\cite{OlfatiSaberMurray}.
However, as these systems are linear with known analytical solution representations, we have considered their numerical simulation to be of less importance and have focused on Kuramoto models as relevant instances.
Generalisations to other nonlinear dynamical systems on graphs occuring in applications will be the objective of future investigations. 

As mentioned in the introduction, it is straightforward to extend our approach to oscillator networks with higher-order interactions included~\cite{BickBoehleKuehn,SkardalArenas}.
Suitable modifications have to be developed for Cucker--Smale models describing flocking behaviour~\cite{CuckerSmale2007}. 
Other types of models, commonly encountered in neuroscience applications, incorporate individual excitable oscillators beyond a phase reduction~\cite{Izhikevich2000,Kopell1995}. 
A further example is the strategy of combining our approach with micro-macro numerical integration schemes such as projective integration methods~\cite{KevrekidisEtAl}, where one makes use in the time integration of a macroscopic evolution, e.g., for the order parameter or the probability density of a typical oscillator, in combination with direct microscopic simulation to improve the numerical simulation. In such a context, our improvements directly accelerate the microscopic integrator. 
Yet, some natural-looking steps are bound to be far more involved.
For example, for time-dependent or adaptive network dynamics~\cite{GrossSayama,HolmeSaramaki}, the recomputation of community structures at every iteration is computationally inefficient and requires novel perspectives. 
Besides, providing a rigorous numerical analysis to assess the quality of community detection algorithms remains an open question.
Furthermore, it is of relevance to study dynamical systems on graphs that incorporate stochastic perturbations.}
\newcommand{\Reference}[3]{\textsc{#1} \newblock \emph{#2} \newblock #3}
\bibliographystyle{unsrt}


\clearpage

\clearpage

\end{document}